\documentclass[a4paper]{amsart}
\usepackage{amsthm} 

\usepackage[leqno]{amsmath}
\usepackage{amsfonts}
\usepackage{amssymb}
\usepackage{amsxtra}
\usepackage [matrix,arrow]{xy}

\title[Characteristic cycles of standard modules]
{Characteristic cycles of standard modules for the rational
Cherednik algebra of type $\Z/l\Z$}
\author{Toshiro Kuwabara}
\address{Department of Mathematics, Graduate School
 of Science, Kyoto University, Kyoto 606-8502, Japan}
\email{toshiro@math.kyoto-u.ac.jp}

\newtheorem{definition}{Definition}[section]
\newtheorem{proposition}[definition]{Proposition}
\newtheorem{theorem}[definition]{Theorem}
\newtheorem{corollary}[definition]{Corollary}
\newtheorem{lemma}[definition]{Lemma}
\newtheorem{remark}[definition]{Remark}

\newcommand{\refprop}[1]{Proposition~\ref{#1}}
\newcommand{\refthm}[1]{Theorem~\ref{#1}}
\newcommand{\refcor}[1]{Corollary~\ref{#1}}
\newcommand{\reflemma}[1]{Lemma~\ref{#1}}
\newcommand{\refeq}[1]{(\ref{#1})}
\newcommand{\refsec}[1]{Section~\ref{#1}}

\renewcommand{\theenumi}{(\arabic{enumi})}

\newcommand{\C}{{\mathbb C}}
\newcommand{\R}{{\mathbb R}}
\newcommand{\Q}{{\mathbb Q}}

\newcommand{\Z}{{\mathbb Z}}

\newcommand{\Hom}{\operatorname{Hom}}
\newcommand{\Dim}{\mathrm{dim}}

\newcommand{\gr}{\operatorname{gr}}

\newcommand{\Ext}{\mathrm{Ext}}

\newcommand{\smashprod}{\#}

\newcommand{\repgeq}[1]{\unrhd_{rep, #1}}

\newcommand{\repgo}[1]{\rhd_{rep, #1}}

\newcommand{\geogeq}[1]{\unrhd_{geom, #1}}
\newcommand{\geoleq}[1]{\unlhd_{geom, #1}}
\newcommand{\geogo}[1]{\rhd_{geom, #1}}

\newcommand{\Rep}{\mathrm{Rep}}
\newcommand{\gl}{\mathfrak{gl}}
\newcommand{\A}{\mathcal{A}}
\newcommand{\B}{\mathcal{B}}
\newcommand{\calF}{\mathcal{F}}
\newcommand{\calM}{\mathcal{M}}
\newcommand{\calP}{\mathcal{P}}
\newcommand{\calL}{\mathcal{L}}
\newcommand{\T}{\mathcal{T}}
\newcommand{\calS}{\mathcal{S}}
\newcommand{\catO}{\mathcal{O}}
\newcommand{\calO}{\mathcal{O}}
\newcommand{\M}{\mathfrak{M}}
\newcommand{\bbT}{\mathbb{T}}
\newcommand{\calV}{\mathcal{V}}
\DeclareMathOperator{\Spec}{Spec}
\DeclareMathOperator{\Proj}{Proj}

\newcommand{\U}{\mathcal{U}}

\newcommand{\Char}{\mathrm{Char}}
\newcommand{\Ch}{\mathbf{\mathrm{Ch}}}
\newcommand{\rCh}{\mathbf{\mathrm{rCh}}}
\newcommand{\Hilb}{\mathrm{Hilb}}
\newcommand{\supp}{\mathrm{Supp}}
\newcommand{\Min}{\mathrm{Min}}
\newcommand{\ann}{\mathrm{ann}}

\newcommand{\modu}[1]{#1\text{-$\mathrm{Mod}$}}
\newcommand{\fmod}[1]{#1\text{-$\mathrm{mod}$}}
\newcommand{\Filt}[1]{#1\text{-$\mathrm{Filt}$}}
\newcommand{\filt}[1]{#1\text{-$\mathrm{filt}$}}

\newcommand{\End}{\mathrm{End}}
\newcommand{\Coh}{\mathrm{Coh}}
\newcommand{\Qcoh}{\mathrm{Qcoh}}
\newcommand{\Pic}{\mathrm{Pic}}

\newcommand{\Qgr}[1]{#1\text{-$\mathrm{Qgr}$}}
\newcommand{\fqgr}[1]{#1\text{-$\mathrm{qgr}$}}
\newcommand{\Grmod}[1]{#1\text{-$\mathrm{Grmod}$}}
\newcommand{\fgrmod}[1]{#1\text{-$\mathrm{grmod}$}}
\newcommand{\Tor}[1]{#1\text{-$\mathrm{Tor}$}}
\newcommand{\ftor}[1]{#1\text{-$\mathrm{tor}$}}

\newcommand{\hwvec}{\mathbf{1}}
\newcommand{\ii}{\eta}
\newcommand{\eu}{\mathrm{eu}}

\newcommand{\frS}{\mathfrak{S}}

\begin{document}
\maketitle

\begin{abstract}
 We study the representation theory of the rational Cherednik
 algebra $H_\kappa = H_\kappa(\Z_l)$ for 
 the cyclic group $\Z_l = \Z / l \Z$ and its connection with
 the geometry of the quiver variety $\M_\theta(\delta)$
 of type $A_{l-1}^{(1)}$.

 We consider a functor between the categories of 
 $H_\kappa$-modules with different parameters, called the
 shift functor, and give the condition when it is an equivalence of
 categories.

 We also consider a functor from the category of 
 $H_\kappa$-modules with good filtration to the category of
 coherent sheaves on $\M_\theta(\delta)$. We prove that 
 the image of the regular representation of $H_\kappa$ by
 this functor is
 the tautological bundle on $\M_\theta(\delta)$. As a
 corollary, we determine the characteristic cycles
 of the standard modules. 
 It gives an affirmative 
 answer to a conjecture given
 in \cite{Go} in the case of $\Z_l$.
\end{abstract}

\section{Introduction}
\label{sec:introduction}

\subsection{Background}
The rational Cherednik algebra for the wreath product
$\Z_l \wr \frS_n$ of the cyclic group $\Z_l = \Z / l \Z$ and
the symmetric group $\frS_n$ is defined by \cite{EG}. 
Let $D(\C^{n}_{reg})$ be the algebra of algebraic differential
operators on $\C^{n}_{reg} = \{(x_1, \dots, x_n)\: | \: x_i \ne 0, 
x_i^l \ne x_j^l \}$. The rational Cherednik algebra is
a subalgebra of the smash product
$D(\C^{n}_{reg}) \smashprod (\Z_l \wr \frS_n)$ which is
generated by the multiplication of functions, 
$\Z_l \wr \frS_n$ and
the Dunkl operators. 
The category of modules over
the rational Cherednik algebra contains an interesting
subcategory called the category $\catO$.
The category $\calO$ is the
subcategory of modules on which the Dunkl 
operators act locally nilpotently.
The category $\calO$ is a highest weight category 
in the sense of \cite{CPS}.

In this paper, we consider the case of $n=1$. Our work is
motivated by the papers \cite{GS1} and \cite{GS2}, in which
the case $l=1$ was considered.
We first review this case.

The rational Cherednik algebra $H_c(\frS_n)$ is the
algebra with a parameter $c \in \R$. We denote the
category $\catO$ of $H_c(\frS_n)$ by $\catO_c(\frS_n)$.
By results of \cite{Op}, \cite{He} and \cite{BEG}, 
we have 
a functor called a shift functor (or a Heckman-Opdam shift functor),
\begin{equation}
\label{eq:39} 
 \widehat{S}_c : \modu{H_c(\frS_n)} \longrightarrow \modu{H_{c+1}(\frS_n)}.
\end{equation}
If the shift functor $\widehat{S}_c$ is an equivalence of categories,
we can construct
a functor from the category of filtered modules $\filt{H_c(\frS_n)}$ to the
category of coherent sheaves $\Coh(\Hilb^n(\C^2))$ on the
Hilbert scheme of $n$ points on $\C^2$,
\begin{equation}
\label{eq:38} 
 \widehat{\Phi}_c : \filt{H_c(\frS_n)} \longrightarrow 
 \Coh(\Hilb^n(\C^2)).
\end{equation}
We recall that 
$\Hilb^n(\C^2)$ is a symplectic resolution
of the singularity $\C^{2n} / \frS_n$,
\[
 \pi : \Hilb^n(\C^2) \longrightarrow \C^{2n} / \frS_n.
\]
These functors $\widehat{S}_c$ and $\widehat{\Phi}_c$ are generalized to
the other cases by \cite{Mu}, \cite{Bo} and \cite{Va}.

In \cite{GS1} and \cite{GS2}, Gordon and Stafford
 also considered the images of
certain modules by $\widehat{\Phi}_c$. 
Consider the rational Cherednik algebra 
$H_c(\frS_n)$ itself as a left $H_c(\frS_n)$-module.  Then the 
corresponding coherent sheaf $\widehat{\Phi}_c(H_c(\frS_n))$
coincides with the Procesi bundle on $\Hilb^n(\C^2)$. 
The Procesi bundle, which was defined in \cite{Ha1}, is a
vector bundle 
whose fiber is isomorphic to the regular representation of 
$\frS_n$. As a corollary of the above result, they described
the images of the standard modules by the functor $\widehat{\Phi}_c$
and determined their characteristic cycles. The 
characteristic cycle $\Ch(M)$ is an invariant of a module
$M$ in $\catO_c(\frS_n)$, which is the sum of 
irreducible components of $\supp \widehat{\Phi}_c(M)$ with multiplicities.
The standard
modules of $H_c(\frS_n)$ are indexed by partitions $\lambda$ of
$n$. Denote them by $\Delta_c(\lambda)$. Let 
$(x_1, \dots, x_n, y_1, \dots, y_n)$ be a coordinate system of $\C^{2n}$.
The irreducible
components of $\pi^{-1}(\{y_1 = \dots = y_n = 0\})$ are
indexed by partitions $\mu$ of $n$. Denote them by
 $\mathcal{Z}_\mu$. Let $[\: \mathcal{Z}_\mu \:]$ be the
homology class given by $\mathcal{Z}_\mu$. 
One of the main results in \cite{GS2} is 
\begin{equation}
\label{eq:2}
 \Ch(\Delta_c(\lambda)) = \sum_{\mu} K_{\lambda \mu} 
[\: \mathcal{Z}_\mu \:]
\end{equation}
for each partition $\lambda$ of $n$. Here 
$K_{\lambda \mu} \in \Z_{\geq 0}$ is the Kostka number.

In the general case, the rational Cherednik algebra
$H_h(\Z_l \wr \frS_n)$ has an $l$-dimensional parameter
$h \in \R^l$. The $l$-multipartitions of $n$ parametrize the
standard modules of the category $\catO_h(\Z_l \wr \frS_n)$.
We have the partial ordering $\repgeq{h}$ on
the set of $l$-multipartitions of $n$, which arises from
the structure of the highest weight category 
$\catO_h(\Z_l \wr \frS_n)$. This ordering $\repgeq{h}$ 
depends on $h$.

Let $\M_h(n \delta)$ be the quiver variety  of type $A_{l-1}^{(1)}$
with the stability 
parameter $h \in \Q^l$ and the dimension vector 
$n \delta = (n, \dots, n)$. 
Similarly to the case of $\frS_n$, the $l$-multipartitions of 
$n$ parametrize the components of 
a certain subvariety of $\M_h(n \delta)$.
The action of $\bbT = (\C^*)^2$ on $\M_h(n \delta)$ induces
a partial ordering $\unrhd_{geom, h}$ on the set of 
$l$-multipartitions of $n$. The ordering $\unrhd_{geom, h}$
also depends on $h$. 

Consider an analogue of the functor \refeq{eq:38} in the
general case. It has been conjectured that the representation theory
of $H_h(\Z_l \wr \frS_n)$ is deeply connected with the geometry of
$\M_h(n \delta)$ (\cite{Go}, \cite{Va}). The functor 
\[
 \widehat{\Phi}_h :
 \filt{H_h(\Z_l \wr \frS_n)} \longrightarrow \Coh(\M_h(n \delta))
\]
is defined for generic $h$ in \cite{Va}.
In \cite{Go}, Gordon compared the ordering $\repgeq{h}$ and
the ordering $\unrhd_{geom, h}$, and proved that
$\unrhd_{geom, h}$ refines $\repgeq{h}$. He conjectured an
analogue of the identity \refeq{eq:2} in the general case
(Question~10.2 in \cite{Go}).

\subsection{Shift functors and their equivalences}
In this paper, we consider the rational Cherednik algebra 
$H_\kappa = H_\kappa(\Z_l)$ in the case of $n=1$. Here
 $\kappa$
is an $(l-1)$-dimensional 
parameter $\kappa = (\kappa_i)_{i=1, \dots, l-1} \in \R^{l-1}$.
Let $\gamma$ be the element of $\Z_l$ which acts  on $\C^*$
by the multiplication of
$\zeta = \exp(2\pi\sqrt{-1}/l)$.
The algebra $H_\kappa$ is the subalgebra of 
$D(\C^*) \smashprod \Z_l$ generated by the
coordinate function $x$, $\gamma \in \Z_l$ and
the Dunkl operator 
$y = (d/dx) + (l/x) \sum_{i=0}^{l-1} \kappa_i \bar{e}_i$ where
$\bar{e}_i = (1/l) \sum_{j=0}^{l-1} \zeta^{ij} \gamma^j \in \C\Z_l$.
As a vector space, we have
\[
 H_\kappa = \C[x] \otimes_{\C} \C \Z_l \otimes_{\C} \C[y].
\]
The algebra $H_\kappa$ is isomorphic to another algebra called
the deformed preprojective algebra defined by 
Crawley-Boevey and Holland in \cite{CBH} and \cite{Ho}. 
Let $Q = (I, E)$ be the Dynkin quiver of type $A^{(1)}_{l-1}$
such that $I = \{I_0, \dots, I_{l-1}\}$ is the set of 
vertices and $E = \{ F_i : I_{i-1} \rightarrow I_{i},
\; i=0, \dots, l-1\}$ is the set of arrows. We regard
indices for vertices and arrows as integers modulo $l$.
For $\lambda = (\lambda_i)_{i=0, \dots, l-1} \in \R^l$
(or $\Z^l$), we regard the sum 
$\lambda_i + \lambda_{i+1} + \dots + \lambda_{j-1}$ as cyclic,
i.e., 
\[
 \lambda_i + \lambda_{i+1} + \dots + \lambda_{j-1}
 = \lambda_i + \dots + \lambda_{l-1} + \lambda_0 +
 \lambda_1 + \dots + \lambda_{j-1}
\]
if $j < i$.
Let 
$\Rep(Q, \delta) \simeq \C^l$ be the space of representations
with the dimension vector $\delta = (1, \dots, 1)$. Set
$GL(\delta) = \prod_{i=0}^{l-1} \C^*$ and set
$\gl(\delta) = \bigoplus_{i=0}^{l-1} \C e^{(i)} = Lie(GL(\delta))$.
Let $t_0$, $\dots$, $t_{l-1} \in \C[\Rep(Q, \delta)]$ be the coordinate
functions, and let $\partial_0$, $\dots$, $\partial_{l-1} \in 
D(\Rep(Q, \delta))$ be the corresponding differential operators.
Let $\R^l_1$ be the set of $\lambda = (\lambda_i)_{i=0, \dots, l-1}
\in \R^l$ such that $\lambda_0 + \cdots + \lambda_{l-1} = 1$.
For a parameter $\lambda = (\lambda_i)_{i=0, \dots, l-1} \in \R^l_1$,
an algebra $\T_\lambda$ is
\[
 \T_\lambda = M_l(D(\Rep(Q, \delta)))^{GL(\delta)}
 \Bigm/ \sum_{i=0}^{l-1} M_l(D(\Rep(Q, \delta)))^{GL(\delta)}
 (\tau(e^{(i)}) - \lambda_i),
\]
where $\tau(e^{(i)}) = E_{ii} \otimes 1 + 1 \otimes 
(t_{i+1} \partial_{i+1} - t_{i} \partial_{i}) \in 
M_l(D(\Rep(Q, \delta)))^{GL(\delta)}$. By \cite[Cor 4.6]{Ho},
this algebra is isomorphic to the deformed preprojective algebra
$\Pi_\lambda$ defined by \cite{CBH}. Moreover, 
the algebra $\T_\lambda$ is isomorphic to $H_\kappa$ for
$\lambda_i = \kappa_{i+1} - \kappa_i + (1/l)$. Set
$e_i = E_{ii}$. Denote the 
spherical subalgebra of $H_\kappa$ by 
$U_\kappa = \bar{e}_0 H_\kappa \bar{e}_0$. The algebra
$U_\kappa$ is isomorphic to the following subalgebra of $\T_\lambda$
\[
 \A_\lambda = e_{0} \T_\lambda e_{0}
\simeq D(\Rep(Q, \delta))^{GL(\delta)}
 \Bigm/ \sum_{i=0}^{l-1} D(\Rep(Q, \delta))^{GL(\delta)}
 (\iota(e^{(i)}) -  \bar{\lambda}_i)
\]
where $\iota(e^{(i)}) = t_{i+1} \partial_{i+1} - t_{i} \partial_{i}$
and $\bar{\lambda}_{i} = \lambda_i - \delta_{i0}$. 
A Dynkin root is a root $\beta = (\beta_i)_{i=0, \dots, l-1}\in \Z^l$ 
such that $\beta_0 = 0$.
Then,  
the $(\A_\lambda, \T_\lambda)$-bimodule $e_0 \T_\lambda$ yields
a Morita equivalence between $\T_\lambda$ and $\A_\lambda$ if
$\lambda = (\lambda_i)_{i=0, \dots, l-1}$ satisfies
$\langle \lambda, \beta \rangle = \sum_{i=0}^{l-1} 
\lambda_i \beta_i \ne 0$ for all Dynkin roots
$\beta = (\beta_i)_{i=0, \dots, l-1} \in \Z^l$.

For $i=0$, $\dots$, $l-1$,
the standard module $\Delta_\lambda(i)$ is the following
$\T_\lambda$-module:
\[
 \Delta_\lambda(i) = (\T_\lambda / \T_\lambda A^*) e_i
\]
where $A^* = \sum_{i=0}^{l-1} E_{i-1,i} \otimes \partial_i \in
\T_\lambda$. We have the isomorphism of vector spaces
\[
 \Delta_\lambda(i) = \C[A] \hwvec_i
\]
where $A = \sum_{i=0}^{l-1} E_{i,i-1} \otimes t_i \in
\T_\lambda$ and $\hwvec_i$ is the image of $e_i$. 

Let $\Z^l_0$ be the set of $\theta = (\theta_i)_{i=0, \dots, l-1} \in
\Z^l$ such that $\theta_0 + \dots + \theta_{l-1} = 0$.
For $\theta \in \Z^l_0$, 
let $\chi_\theta$ be the character of $GL(\delta)$
\[
 \chi_\theta(g) = \prod_{i=0}^{l-1} (g_i)^{\theta_i}
\]
for $g = (g_i)_{i=0, \dots, l-1} \in GL(\delta)$.
We define the shift functor $\calS_\lambda^\theta$ for
$\lambda = (\lambda_i)_{i=0, \dots, l-1} \in \R^l_1$ and 
$\theta = (\theta_i)_{i=0, \dots, l-1} \in \Z^l_0$,
\begin{align*}
 \calS_\lambda^\theta : \fmod{\A_\lambda} &\longrightarrow
 \fmod{\A_{\lambda+\theta}} \\
 N & \mapsto \B_\lambda^\theta \otimes_{\A_\lambda} N
\end{align*}
where $\B_\lambda^\theta$ is the following 
$(\A_{\lambda+\theta}, \A_\lambda)$-bimodule of semi-invariants
\[
 \B_\lambda^\theta = \left[D(\Rep(Q, \delta))
 \Bigm/ \sum_{i=0}^{l-1} D(\Rep(Q, \delta))
 (\iota(e^{(i)}) -  \bar{\lambda}_i)\right]^{GL(\delta), \chi_\theta}.
\]
We also define
the functor,
\begin{align*}
 \widehat{\calS}_\lambda^\theta : \fmod{\T_\lambda} &\longrightarrow
 \fmod{\T_{\lambda+\theta}}, \\
 M & \mapsto \T_\lambda e_0 \otimes_{\A_{\lambda+\theta}}
 \B_\lambda^\theta \otimes_{\A_\lambda} e_0 M.
\end{align*}
Since the algebra $\T_\lambda$ has $(l-1)$-dimensional parameter,
we have the $(l-1)$-dimensional parameter $\theta \in \Z^l$ for
shifting the parameter $\lambda$. Thus we have many shift functors
for the same $\T_\lambda$ while we have only one shift functor $S_c$
for $H_c(\frS_n)$.

We study the case when the shift functors $\calS_\lambda^\theta$ and
$\widehat{\calS}_\lambda^\theta$ are equivalences of categories.
The main difficulty of this question is that we must consider 
complicated combinatorics which depends on the $(l-1)$-dimensional
parameters $\lambda$ and $\theta$.

Define the following sets of parameters
\begin{align*}
 \R^l_{reg} &= \{\lambda = (\lambda_i)_{i=0, \dots, l-1}
\in \R^l_1\; | \; \bar{\lambda}_{i} + \dots + \bar{\lambda}_{j-1} \ne 0
\quad \text{for all $i \ne j$}\}, \\
 \Z^l_{reg} &= \left\{
 \theta \in \Z^l_0 \bigm|
 \theta_i + \dots + \theta_{j-1} \ne 0 \quad
  \text{for all $i \ne j$}
\right\}, \\
 \Z^l_{\lambda} &= \{ \theta = (\theta_i)_{i=0, \dots, l-1}
\in \Z^l_{reg}\; | \;\theta_{i} + \dots + \theta_{j-1} < 0
\quad \text{if $\lambda_{i} + \dots + \lambda_{j-1} \in \Z_{\leq 0}$}\}.
\end{align*}
The set of parameters $\Z^l_{reg}$ is decomposed into 
$(l-1)!$ alcoves. If $\lambda \in \R^l_{reg}$ is generic, we have
 $\Z^l_\lambda = \Z^l_{reg}$. If $\lambda$ belongs to 
$\R^l_{reg} \cap \Z^l$,
$\Z^l_\lambda$ is one of $(l-1)!$ alcoves in $\Z^l_{reg}$.

Then the following theorem is the first main result of this
paper.
\begin{theorem}
\label{thm:7}
 For $\lambda \in \R^l_{reg}$ and $\theta \in \Z^l_\lambda$,
 the shift functors $\calS_\lambda^\theta$ and 
 $\widehat{\calS}_\lambda^\theta$ are equivalences of 
 categories.
\end{theorem}

Moreover, we explicitly determine the images of the standard modules by 
$\calS_\lambda^\theta$.

The parameter $\theta = (\theta_i)_{i=0, \dots, l-1}$ 
defines a total ordering $\unrhd_\theta$ on 
the set of indices $\Lambda = \{0, 1, \dots, l-1\}$,
\[
 i \rhd_{\theta} j \Leftrightarrow
 \theta_{i} + \dots + \theta_{j-1} < 0.
\]
If $\theta$ and $\theta'$ belong to the same alcove in $\Z^l_{reg}$,
then $\rhd_\theta$ is equal to $\rhd_{\theta'}$. 
Let $\repgeq{\lambda}$ be the partial ordering on $\Lambda$ 
defined as $i \repgeq{\lambda} j$ if and only if 
$\Hom_{\T_\lambda}(\Delta_\lambda(j), \Delta_\lambda(i)) \ne 0$.
Our total ordering $\unrhd_\theta$
refines the partial ordering $\repgeq{\lambda}$ when
$\theta \in \Z^l_{\lambda}$, i.e., 
$\Hom_{\T_\lambda}(\Delta_\lambda(j), \Delta_\lambda(i)) \ne 0$
implies $i \unrhd_\theta j$.
Let $\ii_1$, $\dots$, $\ii_l$ be the elements 
of $\Lambda$ such that 
\begin{equation}
\label{eq:48} 
 \ii_l \unrhd_\theta \ii_{l-1} \unrhd_\theta \dots
 \unrhd_\theta \ii_1.
\end{equation}

\begin{proposition}
 \label{prop:19}
 For $i=1$, $\dots$, $l$,
 we have an isomorphism of $\A_{\lambda+\theta}$-modules
\begin{align*}
 e_0 \Delta_{\lambda+\theta}(\ii_i)
 &\longrightarrow
 \calS_\lambda^\theta(e_0 \Delta_\lambda(\ii_i)) = \B_\lambda^\theta
 \otimes_{\A_\lambda} e_0 \Delta_{\lambda}(\ii_i), \\
 e_0 t_{0} t_{l-1} \cdots t_{\ii_i} \hwvec_{\ii_i}
 &\mapsto
 \tilde{f}_i \otimes e_0 t_0 t_{l-1} \cdots t_{\ii_i} \hwvec_{\ii_i}.
\end{align*}
 where
 \[
  \tilde{f}_i = 
 \prod_{j=i+1}^{l-1} (t_{\ii_{j+1}} \cdots t_{\ii_l})^{\theta_{\ii_j} +
 \theta_{\ii_j+1} + \dots + \theta_{\ii_{j+1}-1}}
 \prod_{j=1}^{i-1} (\partial_{\ii_1} \cdots \partial_{\ii_j})^{\theta_{\ii_j} +
 \theta_{\ii_j+1} + \dots + \theta_{\ii_{j+1}-1}}.
\]
\end{proposition}

\subsection{Construction of a tautological bundle}
Next, we consider analogues of the functor \refeq{eq:38} and
determine the image of the regular representation $H_\kappa$
by this functor. 

For $\theta \in \Z_{reg}^l$, the quiver variety
$\M_\theta(\delta)$ with the stability parameter $\theta$ can
be described as follows.
\begin{gather*}
 \M_\theta(\delta) = \Proj S, \\
 S = \bigoplus_{m \in \Z_{\geq 0}} S_m, \qquad
 S_m = \C[\mu^{-1}(0)]^{GL(\delta), \chi_\theta^m}.
\end{gather*}
where $\mu : \Rep(\bar{Q}, \delta) = T^* \Rep(Q, \delta) 
\longrightarrow \gl(\delta)^*$
is the moment map. The variety $\M_\theta(\delta)$ gives a minimal
resolution of the Kleinian singularity $\C^2 / \Z_l$,
\[
 \pi_\theta : \M_\theta(\delta) \longrightarrow \C^2 / \Z_l.
\]
For any $\theta \in \Z^l_{reg}$,
$\M_\theta(\delta)$ is isomorphic to the
toric variety $X(\Delta)$ defined in \refsec{sec:quiver-varieties-vs}.

For $\lambda \in \R^l_{reg}$ and $\theta \in \Z^l_\lambda$, we 
define the functor 
\begin{equation}
\label{eq:13}
\begin{split}
 \widehat{\Phi}_\lambda^\theta : \filt{\T_\lambda} &\longrightarrow
 \Coh(\M_\theta(\delta)) \\
 M &\mapsto \Bigl(\bigoplus_{m \in \Z_{\geq 0}} \gr \bigl(
 \B_\lambda^{m \theta} \otimes_{\A_\lambda} e_0 M \bigr)
\Bigr)\sptilde.
\end{split}
\end{equation}
as in \cite{Bo}. 

We define a locally free sheaf $\widetilde{\calP}_\theta$ on
$\M_\theta(\delta)$ as follows
\[
 \widetilde{\calP}_\theta = \Bigl(\bigoplus_{m \in \Z_{\geq 0}} 
 e_0 M_l(\C[\mu^{-1}(0)])^{GL(\delta), \chi_\theta^m} \Bigr)\sptilde.
\]
Then $\widetilde{\calP}_\theta$ is a tautological bundle of the quiver 
variety $\M_\theta(\delta)$. It is an analogue of the Procesi 
bundle on $\Hilb^n(\C^2)$. 
Although the structure as an algebraic variety of $\M_\theta(\delta)$ is
independent of $\theta$, the tautological bundle 
$\widetilde{\calP}_\theta$ depends on
$\theta$.

We have a construction of the minimal resolution of $\C^2 / \Z_l$
by the toric variety $X(\Delta)$ (\refsec{sec:quiver-varieties-vs}).
We prove the following proposition 
using this construction. This is an analogue of a result
of \cite{Ha2} for the Procesi bundle. 
\begin{proposition}
\label{prop:15}
 For $m \in \Z_{> 0}$, we have
\[
 H^p(\M_\theta(\delta), \widetilde{\calP}_\theta \otimes \calO(m))
 = 
 \left\{
 \begin{array}{ll}
  e_0 M_l(\C[\mu^{-1}(0)])^{GL(\delta), \chi_\theta^m} &
   (p = 0),\\
  0 &  (p \ne 0).
 \end{array}
 \right.
\]
 where $\calO(1)$ is the twisting sheaf of of $\M_\theta(\delta)$
 associated to the homogeneous coordinate ring $S$
 (see \cite[p.117]{Har}).
\end{proposition}
We make use of this
proposition to calculate the $(q,t)$-dimension of the module
$e_0 M_l(\C[\mu^{-1}(0)])^{GL(\delta), \chi_\theta^m}$, i.e.,
the character with respect to the $\bbT$-action.

Using this result, we obtain the following second main result.
Set 
\[
 \widetilde{\R}^l_{reg} = \{(\lambda_i)_{i=0, \dots, l-1} \in \R^l_{reg}
\; | \; \lambda_i + \dots + \lambda_{j-1} \ne 0 \quad
\text{for all $i \ne j$}\}.
\]
\begin{theorem}
\label{thm:6}
 For $\lambda \in \widetilde{\R}^l_{reg}$ and $\theta \in \Z^l_\lambda$,
we have an
 isomorphism  of coherent sheaves on $\M_\theta(\delta)$
 \[
  \widehat{\Phi}_\lambda^{\theta}(\T_\lambda) \simeq \widetilde{\calP}_\theta.
 \]
\end{theorem}

As a corollary of \refthm{thm:6}, we have
\begin{corollary}
 \label{cor:7}
 For $\lambda \in \widetilde{\R}^l_{reg}$, $\theta \in \Z^l_\lambda$,
 we have the 
 isomorphism
 \[
  \widehat{\Phi}_\lambda^{\theta}(\Delta_\lambda(i)) \simeq 
 \bigl(\widetilde{\calP}_\theta / \widetilde{\calP}_\theta \bar{A}^* \bigr) e_i
 \]
 where $\bar{A}^* = \sum_{i=0}^{l-1} E_{i-1,i} \otimes \xi_i
 \in M_l(\C[\mu^{-1}(0)])$ and $\xi_i = \overline{\partial_i}$ is the
 image of $\partial_i$ 
in $\gr D(\Rep(Q, \delta))^{GL(\delta)} \simeq \C[\mu^{-1}(0)]^{GL(\delta)}$.
\end{corollary}

\subsection{The characteristic cycles of the standard modules}
Finally, we determine the characteristic cycles of the standard
modules. 

The structure of the subvariety $\pi_\theta^{-1}(\{y=0\})$ 
is well-known: we have
\[
 \pi_\theta^{-1}(\{y=0\}) = \bigsqcup_{i=0}^{l-1} \U^{0}_i
\]
where $\U_i^0$ is one-dimensional affine subvariety depending on
$\theta$
defined by
\refeq{eq:41} in \refsec{sec:defin-quiv-vari}.
We denote by $\U_i$ the closure of $\U^0_i$. Then, the irreducible
components of $\pi_\theta^{-1}(\{y=0\})$ are $\U_0$, $\dots$, $\U_{l-1}$.

The following proposition is the third main result.
\begin{proposition}
 \label{prop:16}
 For $i=1$, $\dots$, $l$,  we have
 \[
  \Ch(\Delta_\lambda(\ii_i)) = 
 \sum_{j=1}^{i} [\: \U_{\ii_j} \:].
 \]
 where $\ii_i$ is the index defined on \refeq{eq:48}.
\end{proposition}

This proposition answers a conjecture in \cite{Go} in the case
of $\Z_l$. 
Consider the geometric ordering defined in Section~5.4 of
\cite{Go},
\[
 i \geogo{\theta} j \quad \text{if $i \ne j$ and $\U_i \cap \U^0_j \ne \emptyset$.}
\]
By the general theory of quiver varieties, if $\theta$ and $\theta'$
belong to the same alcove in $\Z^l_{reg}$, $\rhd_{geom, \theta}$ 
is equal to $\rhd_{geom, \theta'}$.
Moreover, we show that the geometric ordering $\geogo{\theta}$ is equal to the
ordering $\rhd_{\theta}$. Thus, we have
\[
 \ii_l \geogo{\theta} \ii_{l-1} \geogo{\theta} \dots \geogo{\theta} \ii_1. 
\]

Therefore, \refprop{prop:16} is written as 
\[
  \Ch(\Delta_\lambda(i)) =  \sum_{j \geoleq{\theta} i}
 [\: \U_{j} \:].
\]
This is an affirmative answer to Question 10.2 of \cite{Go}.

The functor $\widehat{\Phi}_\lambda^\theta$ of \refeq{eq:13} confirms 
a deep connection between the representation theory of 
$\T_\lambda$ and the geometry of $\M_\theta(\delta)$. In fact, 
to prove \refthm{thm:7} and \refthm{thm:6}, we make use of this
connection. Although \refprop{prop:19} is purely 
representation theoretical, the elements $\tilde{f}_i$ are
obtained from the information on the geometry of $\M_\theta(\delta)$.

\subsection{Plan of paper}
The paper is organized as follows. In \refsec{sec:defin-quiv-vari}, 
we recall the definition and basic facts about the 
quiver variety $\M_\theta(\delta)$. In 
\refsec{sec:quiver-varieties-vs}, we recall the construction 
of the minimal resolution of the singularity $\C^2 / \Z_l$ as
 a toric variety,
and compare it with $\M_\theta(\delta)$. In 
\refsec{sec:tautological-bundles}, we construct the tautological
bundle $\widetilde{\calP}_\theta$. In \refsec{sec:line-bundles}, we prove
\refprop{prop:15} by using well-known facts about line bundles on
toric varieties. In \refsec{sec:q-t-dimension} and 
\refsec{sec:proof-theorem}, we calculate the $(q, t)$-dimension
of $\C \otimes_{\C[t_0 \cdots t_{l-1}]} 
e_0 M_l(\C[\mu^{-1}(0)])^{GL(\delta), \chi_\theta}$. 
This result is used when we prove \refthm{thm:6}.
In \refsec{sec:rati-cher-algebr} and \refsec{sec:deform-prepr-algebr},
we define the rational Cherednik algebra $H_\kappa$ and
the deformed preprojective algebra $\T_\lambda$, and recall 
their fundamental properties. In \refsec{sec:parameters-orderings}, we
prepare some conditions for the parameters $\lambda$ and $\theta$,
and define the ordering
$\unrhd_\theta$. In \refsec{sec:shift-functors}, we define
the shift functor $\calS_\lambda^\theta$ and prove \refthm{thm:7}.
In \refsec{sec:q-dimens-repr}, we calculate the $q$-dimension
of $\C \otimes_{\C[t_0 \cdots t_{l-1}]} \B_\lambda^\theta
\otimes_{\A_\lambda} e_0 \T_\lambda$. We use this result to prove
\refthm{thm:6}. 
In \refsec{sec:gord-staff-funct},
we recall the definition of the functor $\widehat{\Phi}_\lambda^\theta$.
In \refsec{sec:proof-main-theorem},
 we prove \refthm{thm:6}.
In \refsec{sec:char-cycl}, we 
prove \refcor{cor:7} and \refprop{prop:16}.

\subsection*{Acknowledgments.}
I am deeply grateful to Iain Gordon for useful advices and
discussions during the study of this subject.
I also thank him for reading manuscript and giving
some valuable comments. 

I am deeply grateful to Susumu Ariki, Masaki Kashiwara,
 and Raphael Rouquier for
useful discussions, Richard Vale for showing me a draft of
\cite{Va}.
I also thank my adviser Tetsuji
Miwa for reading the manuscript and for his kind encouragement.

\section{Preliminaries}
\label{sec:preliminaries}
\subsection{Basic notations}
\label{sec:basic-notations}

Fix an integer $l \in \Z_{> 0}$. We define two sets of
parameters. 
Let $\Z^l_0$ be the set of
$\theta = (\theta_i)_{i=0, \dots, l-1} \in \Z^l$ such that 
$\theta_0 + \dots + \theta_{l-1} = 0$. Let $\R^l_1$ be the
set of $\lambda = (\lambda_i)_{i=0, \dots, l-1} \in \R^l$ such that
$\lambda_0 + \dots + \lambda_{l-1} = 1$.

Set
$\gamma = \begin{pmatrix}
	   \zeta & 0 \\
	   0 & \zeta^{-1}
	  \end{pmatrix}$
to be the element of $SL_2(\C)$ where 
$\zeta = \mathrm{exp}(2\pi \sqrt{-1} / l)$.
Let $\Z_l = \Z / l \Z$ be the finite subgroup of $SL_2(\C)$ 
generated by the element $\gamma$. 
Denote the group ring of $\Z_l$ over the field $\C$ by $\C \Z_l$.
For $i = 0$, 
$\dots$, $l-1$, let $\bar{e}_i$ be the idempotent
$\bar{e}_i = (1/l) \sum_{j=0}^{l-1} \zeta^{ij} \gamma^j \in \C\Z_l$.
Then we have $\bar{e}_i \bar{e}_j = \delta_{ij} \bar{e}_i$, and
$\C \Z_l = \bigoplus_{i=0}^{l-1} \C \bar{e}_i$.
For $i = 0$, $\dots$, $l-1$, let $L_i = \C \hwvec_i$ be
the one-dimensional irreducible representation of $\Z_l$
on which $\bar{e}_j$ acts by $\bar{e}_j \hwvec_i = \delta_{ij} \hwvec_i$.

For a group $G$ and a $G$-module $M$, we denote by $M^{G}$ the 
$G$-invariant subspace of $M$. 
For a character $\chi$ of the group $G$, we denote by $M^{G, \chi}$
the semi-invariant subspace of $M$ belonging to the 
character $\chi$, i.e., 
\[
 M^{G, \chi} = \{v \in M\; | \;g v = \chi(g) v \quad (g \in G)\}.
\]

For a $\C$-algebra $R$, let $M_l(R)$ be the $l \times l$ matrix algebra 
whose elements have coefficients in $R$, i.e.
\[
 M_l(R) \simeq M_l(\C) \otimes_\C R.
\]
Let $E_{ij}$ be the matrix in $M_l(R)$ such that the 
$(i', j')$-entry of $E_{ij}$ is given by $\delta_{ii'} \delta_{jj'}$.

For an algebra $R$, denote by $\modu{R}$  the category of $R$-modules.
Let $\fmod{R}$ be the full subcategory of $\modu{R}$ whose
objects are finitely generated over $R$. 

Fix an affine variety $X$. Let $\C[X]$ be the ring of
polynomial functions on $X$. Let $D(X)$ be the ring of 
algebraic differential operators on $X$. For an $\calO_X$-module
$\calF$ and a point $x \in X$, we denote
the stalk of $\calF$ at $x$ by $\calF_x$. We define
its fiber at $x$ by $\calF(x) = \calF_x \otimes_{\calO_{X,x}} \C$
where $\C = \calO_{X,x} / \mathfrak{m}_x$ and $\mathfrak{m}_x$ is 
a unique maximal ideal of $\calO_{X,x}$.
We denote by $\Qcoh(X)$ the category of
quasi-coherent sheaves on $X$ and by $\Coh(X)$ the category of
coherent sheaves on $X$.

\subsection{Quivers}
\label{sec:quivers}

A quiver $Q = (I, E)$ is a pair of a set of vertices $I$ and
a set of arrows $E$ equipped with two maps $in$, 
$out : E \rightarrow I$. We assume $I$ and $E$ are finite
sets.
Let $Q^* = (I, E^*)$ be
the quiver with the same set of vertices $I$ and the set
of arrows 
$E^* = \{\alpha^* \; | \; \alpha \in E\}$ where 
$\alpha^{*}$ is an arrow such that $in(\alpha^{*}) = out(\alpha)$
and $out(\alpha^{*}) = in(\alpha)$. 
Let $\bar{Q}$ be 
the quiver $(I, E \sqcup E^*)$. 

A representation of a quiver $Q = (I, E)$  is a pair 
$(V, (\phi_\alpha)_{\alpha \in E})$ of an $I$-graded vector space
$V = \bigoplus_{i \in I} V_i$ and a set of linear maps 
$\phi_\alpha : V_{out(\alpha)} \rightarrow V_{in(\alpha)}$. For 
an $I$-graded vector space $V = \bigoplus_{i \in I} V_i$, its
dimension vector is $\underline{\dim} V = ( \dim V_i )_{i \in I}
\in (\Z_{\geq 0})^I$.
For a dimension vector $v = (v_i)_{i \in I}$, the space of 
representation $\Rep(Q, v)$ is the following space: 
\[
 \Rep(Q, v) = \bigoplus_{\alpha \in E} 
\Hom(\C^{v_{out(\alpha)}}, \C^{v_{in(\alpha)}}).
\]
We identify a point $(\phi_{\alpha})_{\alpha \in E}$ of 
$\Rep(Q, v)$ and a representation 
$(\bigoplus_{i \in I} \C^{v_i}, (\phi_{\alpha})_{\alpha \in E})$.

Fix a dimension vector $v = (v_i)_{i \in I}$. Let $GL(v)$ be the 
Lie group $\prod_{i \in I} GL_{v_i}(\C)$ and let 
$\gl(v)$ be the Lie algebra of $GL(v)$: 
$\gl(v) = \bigoplus_{i \in I} \gl_{v_i}(\C)$. The group $GL(v)$
acts naturally on $\Rep(Q, v)$.

Let $\langle \; , \; \rangle: \R^I \times \R^I \longrightarrow
\R$ be the bilinear form 
$\langle \lambda, \mu \rangle = \sum_{i \in I} \lambda_i \mu_i$
for $\lambda = (\lambda_i)_{i \in I}$, 
$\mu = (\mu_i)_{i \in I}$.

\subsection{Quiver of type $A_{l-1}^{(1)}$}
\label{sec:quiver-type-a_l}
In this paper, we consider the McKay quiver associated to a
group $\Z_l = \Z / l \Z$ with cyclic orientation. 
In other words, it is a Dynkin quiver of type 
$A_{l-1}^{(1)}$ with cyclic orientation. 
Let $Q = (I, E)$ be a quiver with 
$I = \{I_0, \dots, I_{l-1}\}$ as the set of vertices and 
$E = \{ F_i : I_{i-1} \rightarrow I_{i}\:|\: i = 0, \dots, l-1\}$
as the set of arrows. We regard indices for vertices and arrows as
integers modulo $l$, i.e., we consider $I_{-1} = I_{l-1}$,
$F_l = F_0$ and so on.
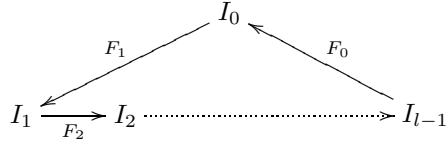
\[
 \xymatrix{
& & I_0 \ar[dll]_{F_1} & & \\
I_1 \ar[r]_{F_2} & I_2 \ar@{.>}[rrr] & & &
I_{l-1} \ar[ull]_{F_{0}}
}
\]

We identify $\Z^I = \Z^l$ and the 
root lattice of type $A^{(1)}_{l-1}$, and identify 
$\R^I = \R^l$ and the dual of the Cartan subalgebra. 
For $i=0$, $\dots$, $l-1$, let 
$\epsilon_i \in \Z^l$ be the standard coordinate vector corresponding to
the vertex $I_i$. Under the above identification, 
we regard $\epsilon_0$, $\dots$, $\epsilon_{l-1}$
as simple roots. Let $\delta = (1, \dots, 1)$ be the
minimal positive imaginary root. A root $\beta \in \Z^l$ is called
a Dynkin root when $\langle \beta, \epsilon_0 \rangle = 0$.
When a positive root $\beta = (\beta_i)_{i=0, \dots, l-1} \in \Z^l$
satisfies $\beta \ne \delta$ and $\beta_i \leq \delta_i = 1$
for all $i = 0$, $\dots$, $l-1$, we write $\beta < \delta$.
For $\lambda = (\lambda_i)_{i=0, \dots, l-1} \in \R^I = \R^l$
(or $\Z^I = \Z^l$), we regard the sum 
$\lambda_i + \lambda_{i+1} + \dots + \lambda_{j-1}$ as cyclic,
i.e., 
\[
 \lambda_i + \lambda_{i+1} + \dots + \lambda_{j-1}
 = \lambda_i + \dots + \lambda_{l-1} + \lambda_0 +
 \lambda_1 + \dots + \lambda_{j-1}
\]
if $j < i$.

\section{Quiver varieties}
\label{sec:quiver-varieties}
\subsection{Definition of quiver varieties}
\label{sec:defin-quiv-vari}

In this subsection we review the definition and fundamental 
properties of quiver varieties which were introduced by Nakajima
in \cite{Na1}. 

Define the quiver $Q = (I, E)$
as in \refsec{sec:quivers}.
The space of representations
\[
 \Rep(\bar{Q}, \delta) 
 = \{ (a_i, b_i)_{i=0, \dots, l-1} \; | \; a_i, b_i \in \C\} 
 \simeq \C^{2l}
\]
is a symplectic manifold with the symplectic form 
$\sum_{i=0}^{l-1} d b_i \wedge d a_i$. 
We have the symplectic action of $GL(\delta)$ on $\Rep(\bar{Q}, \delta)$,
and the corresponding moment map is
\begin{align*}
 \mu : \Rep(\bar{Q}, \delta) &\longrightarrow \gl(\delta)^{*} \simeq
 \C^l, \\
 (a_i, b_i)_{i=0, \dots, l-1} &\mapsto (a_{i} b_{i} - a_{i+1} b_{i+1})_{i=0,
 \dots, l-1}.
\end{align*}

Let $t_0$, $\dots$, $t_{l-1}$ and $\xi_0$, $\dots$, $\xi_{l-1}
\in \C[\Rep(\bar{Q}, \delta)]$ be
the coordinate functions such that $t_i((a_j, b_j)_{j=0, \dots, l-1}) = a_i$
and $\xi_i((a_j, b_j)_{j=0, \dots, l-1}) = b_i$ for 
$(a_j, b_j)_{j=0, \dots, l-1} \in \Rep(\bar{Q}, \delta)$. 
Then we have 
\[
\C[\mu^{-1}(0)] = \C[t_0, \dots, t_{l-1}, \xi_0, \dots, \xi_{l-1}]
\bigm/ (t_i \xi_i - t_{i+1} \xi_{i+1})_{i=0, \dots, l-1}.
\]
The group $GL(\delta)$ acts on $\C[\mu^{-1}(0)]$ as follows.
\begin{align*}
 g \cdot t_i &= g_i^{-1} g_{i-1} t_i, \\
 g \cdot \xi_i &= g_i g_{i-1}^{-1} \xi_i 
\end{align*}
for $i=0$, $\dots$, $l-1$ and 
$g = (g_k)_{k=0, \dots, l-1} \in GL(\delta)$.

Fix a parameter $\theta = (\theta_i)_{i=0, \dots, l-1} \in \Z^l_0$ 
called a stability parameter. For a representation 
$(V, (a_i, b_i)_{i=0, \dots, l-1})$ of
$\bar{Q}$ with the dimension vector $\delta$, 
we call it $\theta$-semistable if 
$\langle \underline{\Dim} W, \theta \rangle \leq 0$ for any 
subrepresentation $W$ of $V$. Define the subset $\mu^{-1}(0)_\theta$
of $\mu^{-1}(0) \subset \Rep(\bar{Q}, \delta)$:
\begin{multline*}
 \mu^{-1}(0)_\theta = 
 \bigl\{ (a_i, b_i)_{i=0, \dots, l-1} \in \mu^{-1}(0) \bigm| \\
 \text{$(a_i, b_i)_{i=0, \dots, l-1}$ is a $\theta$-semistable
 representation.}
 \bigr\}.
\end{multline*}
It is a Zariski open subset of $\mu^{-1}(0)$. 
For $p$, $q \in \mu^{-1}(0)_\theta$, we denote $p \sim q$ when
the closures of $GL(\delta)$-orbits intersect in $\mu^{-1}(0)_\theta$.
Then $\sim$ is an equivalence relation. 
Then we define the quiver variety $\M_\theta(\delta)$ as follows:
\[
 \M_\theta(\delta) = \mu^{-1}(0)_\theta / \sim.
\]
For a point $(a_i, b_i)_{i=0, \dots, l-1} \in \mu^{-1}(0)_\theta$, 
we denote by
$[a_i, b_i]_{i=0, \dots, l-1}$ the corresponding point of
$\M_\theta(\delta)$.

\begin{remark}
 Although our definition of the quiver variety is different
 from one in \cite{Na1}, $\M_\theta(\delta)$ coincides with
 $\M_{(\zeta_\C, \zeta_\R)}(\mathbf{v}, \mathbf{w})$ with
 $\zeta_\C = 0$, $\zeta_\R = \theta$, $\mathbf{v} = \delta$ and
 $\mathbf{w} = \epsilon_0$. This definition is the same as
 one of \cite[Section 4]{Na2}.
\end{remark}

For $\theta = 0 = (0, \dots, 0)$, we have
\[
 \M_0(\delta) \simeq \C^2 / \Z_l
\]
(see \cite{CS}).

\begin{proposition}[\cite{Kr}, \cite{Na2}]
\label{prop:1}
 If a stability parameter $\theta = (\theta_i)_{i=0, \dots, l-1} \in
 \Z^l_0$ satisfies $\langle \theta, \beta \rangle \ne 0$ for 
all positive roots $\beta$ which satisfy $\beta < \delta$, 
$\M_\theta(\delta)$ is
nonsingular and we have a minimal resolution of Kleinian
singularities of type $A_{l-1}$:
\[
\label{eq:15} 
\pi_\theta: \M_\theta(\delta) \longrightarrow \M_0(\delta) \simeq \C^2 / \Z_l.
\]
\end{proposition}

In this paper, we always consider the case when the stability
parameter $\theta \in \Z^l_0$ satisfies 
$\langle \theta, \beta \rangle \ne 0$ for all positive roots
$\beta$ which satisfy $\beta < \delta$. Set
\begin{equation}
 \label{eq:16}
 \Z^l_{reg} = \left\{
 \theta \in \Z^l_0 \bigm|
 \langle \theta, \beta \rangle \ne 0 \quad
 \text{for all positive roots $\beta$ which satisfy $\beta < \delta$}
\right\}.
\end{equation}

For a stability parameter  
$\theta = (\theta_i)_{i=0, \dots, l-1} \in \Z^l_{reg}$,
define the following graded commutative algebra
\begin{align*}
 S &= \bigoplus_{m \in \Z_{\geq 0}} S_m, \\
 S_m &= \C[\mu^{-1}(0)]^{GL(\delta), \chi_\theta^m},
\end{align*}
where $\chi_\theta$ is the character of $GL(\delta)$ given by 
$\chi_\theta(g) = \prod_{i=0}^{l-1} (g_i)^{\theta_i}$ for
$g = (g_i)_{i=0, \dots, l-1} \in GL(\delta)$.
The injective homomorphism $S_0 \rightarrow S$ induces the
morphism of schemes
\[
 \Proj S \longrightarrow \Spec S_0 \simeq \C^2 / \Z_l.
\]

We have the following construction of quiver varieties.
\begin{proposition}[\cite{CS}, \cite{Na2}]
\label{prop:2}
As schemes over $\C^2 / \Z_l$, we have the following isomorphism:
 \[
  \M_\theta(\delta) \simeq \Proj S.
 \]
\end{proposition}
The above construction induces the twisting sheaf on 
the scheme $\M_\theta(\delta) \simeq \Proj S$ which we 
denote by $\calO(1)$ (see \cite[p.117]{Har}).

Fix a stability parameter $\theta = (\theta_i)_{i=0, \dots, l-1} 
\in \Z^l_{reg}$. The
two-dimensional torus $\bbT = (\C^*)^2$ acts on the quiver
variety $\M_\theta(\delta)$ as follows: for 
$[a_i, b_i]_{i=0, \dots, l-1} \in \M_\theta(\delta)$ and
$(m_1, m_2) \in \bbT$,
\[
 (m_1, m_2) [a_i, b_i]_{i=0, \dots, l-1}
 = [m_1 a_i, m_2 b_i]_{i=0, \dots, l-1}.
\]
The group $\bbT$ acts on $\C^2$ by 
$(m_1, m_2) (a, b) \mapsto (m_1 a, m_2 b)$ and it induces the
$\bbT$-action on $\C^2 / \Z_l$. Moreover $\pi_\theta$ is 
$\bbT$-equivariant. 
The variety $\M_\theta(\delta)$ has $l$ $\bbT$-fixed points 
$p_0$, \dots, $p_{l-1}$ where $p_i = [a_j, b_j]_{j=0, \dots, l-1}$ 
is given as follows:
\begin{gather}
  a_i = 0,\quad b_i = 0, \nonumber\\
  a_j = 0,\quad b_j \ne 0 \quad \text{if $\theta_{i} + \theta_{i+1} +
 \dots + \theta_{j-1} < 0$}, \nonumber\\
\label{eq:34}
  a_j \ne 0,\quad b_j = 0 \quad \text{if $\theta_{i} + \theta_{i+1} +
 \dots + \theta_{j-1} > 0$}.
\end{gather}
Note that we have $\theta_i + \theta_{i+1}
\dots + \theta_{j-1} \ne 0$ for all $i \ne j$ 
by the assumption $\theta \in \Z_{reg}^l$.

Define the ordering $\geogeq{\theta}$ on the set $\Lambda = \{0, \dots l-1\}$ by 
\begin{equation}
 \label{eq:35}
 i \geogo{\theta} j \Longleftrightarrow 
 \theta_{i} + \dots + \theta_{j-1} < 0.
\end{equation}
Since we take $\theta \in \Z^l_{reg}$, the ordering $\geogeq{\theta}$
is a total ordering. 

Set $\ii_1$, $\dots$, $\ii_l$ be the indices in $\Lambda$ such that
\begin{equation}
\label{eq:36} 
 \ii_l \geogo{\theta} \ii_{l-1} \geogo{\theta} \dots \geogo{\theta} \ii_1.
\end{equation}
By \refeq{eq:34} and \refeq{eq:35}, for 
$p_{\ii_i} = [a_j, b_j]_{j=0, \dots, l-1}$, we have
\begin{equation}
\label{eq:37}
\#\{ j \in \Lambda \; | \; b_j \ne 0 \} = i-1. 
\end{equation}

\begin{proposition}
 \label{prop:18}
 For $i=1$, $\dots$, $l$, the fixed point 
$p_{\ii_i} = [a_j, b_j]_{j=0, \dots, l-1}$ is given by
 \begin{gather*}
  a_{\ii_i} = 0, \quad b_{\ii_i} = 0, \\
  a_{\ii_j} = 0, \quad b_{\ii_j} \ne 0 \qquad \text{for $j < i$}, \\
  a_{\ii_j} \ne 0, \quad b_{\ii_j} = 0 \qquad \text{for $j \geq i$}.
 \end{gather*}
\begin{proof}
 By \refeq{eq:35} and \refeq{eq:36}, we have 
 \[
 \theta_{\ii_i} + \dots + \theta_{\ii_j - 1} < 0
 \qquad \text{for $j=1$, $\dots$, $i-1$}.
 \]
Thus we have
\[
  a_{\ii_j} = 0, \quad b_{\ii_j} \ne 0 \qquad \text{for $j < i$}.
\]
By \refeq{eq:37}, we have
\[
  a_{\ii_j} \ne 0, \quad b_{\ii_j} = 0 \qquad \text{for $j \geq i$}. 
\] 
\end{proof}
\end{proposition}

For $i=0$, $\dots$, $l-1$, we define the following 
one-dimensional 
affine subvariety of $\M_\theta(\delta)$
\begin{equation}
\label{eq:41} 
  \U^0_i = 
 {\left\{ [a_j,b_j]_{j=0,\dots,l-1} 
 \biggm|
 \begin{array}{l}
  b_i = 0, \\
  a_j = 0, b_j \ne 0 \quad \text{if $\theta_{i} + \dots + \theta_{j-1} < 0$}, \\
  a_j \ne 0, b_j = 0 \quad \text{if $\theta_{i} + \dots + \theta_{j-1} > 0$}.
 \end{array}
 \right\}} \subset \M_{\theta}(\delta).
\end{equation}
Clearly $p_i \in \U^0_i$ for all $i=0$, $\dots$, $l-1$. We denote
 by $\U_i$ the closure of $\U^0_i$.
By
\refprop{prop:18}, $\U_{\ii_i}$ is given as follows.
\begin{equation}
\label{eq:40} 
  \U^0_{\ii_i} = 
 {\left\{ [a_j,b_j]_{j=0,\dots,l-1} 
 \biggm|
 \begin{array}{l}
  b_{\ii_i} = 0, \\
  a_{\ii_j} = 0, \; b_{\ii_j} \ne 0 \qquad \text{for $j < i$}, \\
  a_{\ii_j} \ne 0, \; b_{\ii_j} = 0 \qquad \text{for $j \geq i$}.
 \end{array}
 \right\}}.
\end{equation}

The ordering $\geogeq{\theta}$ is related with the $\bbT$-action 
on $\M_\theta(\delta)$ as follows. We denote $p_i \rightarrow p_j$ when
there is a point $[a_k, b_k]_{k=0, \dots, l-1}
\in \pi_\theta^{-1}(0)$ such that
\[
 \lim_{m \rightarrow 0} (m^{-1}, m) [a_k, b_k]_{k=0, \dots, l-1}
  = p_{i},
\]
and
\[
 \lim_{m \rightarrow 0} (m, m^{-1}) [a_k, b_k]_{k=0, \dots, l-1}
  = p_{j}.
\]
By \refprop{prop:18} and \refeq{eq:40}, taking 
$[a_k, b_k]_{k=0, \dots, l-1}$ on $\U_{\ii_i}$,
 we have $p_{\ii_{i}} \rightarrow p_{\ii_{i-1}}$ for each
$i=2$, $\dots$, $l$. 

The structure of the subvariety $\pi_\theta^{-1}(\{y=0\})$ is 
well-known (see \cite[Lecture 1]{Sl}).
The subvariety $\pi_\theta^{-1}(\{y=0\})$ is 
the disjoint union of $\U_i^0$:
\[
 \pi_\theta^{-1}(\{y=0\}) = \bigsqcup_{i=0}^{l-1} \U^0_{i}
\]
The irreducible components of $\pi_\theta^{-1}(\{y=0\})$ are
$\U_0$, $\dots$, $\U_{l-1}$.

For $i=1$, $\dots$, $l-1$, we have
\[
 \U_{\ii_{i+1}} \cap \U^0_{\ii_i}  = \{p_{\ii_i} \} 
\]
and $\U_{\ii_j}$ does not intersect with $\U_{\ii_i}$ unless
$j=i+1$, $i$, $i-1$.

\subsection{Quiver varieties vs. toric varieties}
\label{sec:quiver-varieties-vs}

In this subsection, we compare two constructions of 
the minimal resolution of Kleinian singularity $\C^2 / \Z_l$:
i.e. as a quiver variety and as a toric variety.

Let $N = \Z^2$, $M = \Hom(N, \Z) \simeq \Z^2$ and let
\[
 \langle\; , \;\rangle : M \times N \longrightarrow \Z
\]
be the natural pairing. Set $v_i = (1, i) \in N$  for 
$0 \leq i \leq l$. Let 
$\sigma_i = \R_{\geq 0} v_i + \R_{\geq 0} v_{i-1}$  be a
$2$-dimensional cone for $i=1$, $\dots$, $l$. Let
$\Delta$ be the fan 
with the $2$-dimensional cones $\sigma_i$ for $i=1$, $\dots$,
$l$ and the $1$-dimensional cones $\R_{\geq 0} v_i$ for $i=0$, $\dots$, $l$.
Then the toric variety $X(\Delta)$ associated to the fan 
$\Delta$ gives the minimal resolution of $\C^2 / \Z_l$:
\begin{equation}
\label{eq:17} 
 X(\Delta) \longrightarrow \C^2 / \Z_l.
\end{equation}

For $i = 1$, $\dots$, $l$,
let $M_i = M \cap \check{\sigma}_i$ be the semigroup 
where $\check{\sigma}_i = \R_{\geq 0} (i, -1) + \R_{\geq 0} (1-i, 1)$
is the dual cone of $\sigma_i$.
Let $R_i = \C M_i$ be the group ring of $M_i$ and
let $X_i = \Spec R_i$.
The toric variety $X(\Delta)$ has the open covering 
$X(\Delta) = \bigcup_{i=1}^{l} X_i$.

Let $u = (1, 0)$, $v = (0, 1)$ be the basis of the lattice
$M$. Then, $R_i = \C[u^i v^{-1}, u^{1-i} v]$. Let
$xy = u$, $y^l = v$. 
Then we have
\[
 R_i = \C[x^i y^{i-l}, x^{1-i} y^{l+1-i}].
\]
The natural embedding $\C[x,y]^{\Z_l} \rightarrow R_i$ induces
the morphism $X_i \rightarrow \C^2 / \Z_l$ of \refeq{eq:17}.

Let $\mathfrak{m}_i = (x^i y^{i-l}, x^{1-i} y^{l+1-i}) \subset R_i$
be the maximal ideal of $R_i$. It is the maximal ideal 
corresponding to the unique $\bbT$-fixed point in $X_i$.

Fix a stability parameter 
$\theta = (\theta_i)_{i=0, \dots, l-1} \in \Z^l_{reg}$. We consider
the quiver variety $\M_\theta(\delta)$ defined in 
\refsec{sec:defin-quiv-vari}. 
By \refprop{prop:1}, we also have the minimal resolution
\[
 \M_\theta(\delta) \longrightarrow \C^2 / \Z_l.
\]
Thus we have an isomorphism of algebraic varieties
\[
\M_\theta(\delta) \simeq X(\Delta).
\]
We construct this
isomorphism explicitly. 

For $i=1$, $\dots$, $l$, let
\[
 R'_i = 
 \C\left[
 \frac{t_{\ii_{l-i+1}} \cdots t_{\ii_l}}
 {\xi_{\ii_{1}} \cdots \xi_{\ii_{l-i}}},
 \frac{\xi_{\ii_1} \cdots \xi_{\ii_{l-i+1}}}
 {t_{\ii_{l-i+2}} \cdots t_{\ii_{l}}}
 \right]
\]
where $\ii_1$, $\dots$, $\ii_l$ are the indices defined on
\refeq{eq:36}. Note that the polynomials in $R'_i$ have no
poles at the fixed point $p_{\ii_{l-i+1}}$ by \refprop{prop:18}.

Then we have an open covering
\[
 \M_\theta(\delta) = \bigcup_{i=1}^l X'_i, \qquad
 X'_i = \Spec R'_i.
\]
Thus, for $i=1$, $\dots$, $l$, we define an isomorphism,
\begin{align*}
 R'_i &\longrightarrow R_i, \\
 t_j &\mapsto x \qquad (j=0, \dots, l-1), \\
 \xi_j &\mapsto y \qquad (j=0, \dots, l-1).
\end{align*}
This induces the isomorphism of algebraic varieties
\begin{equation}
\label{eq:18} 
\M_\theta(\delta) \longrightarrow X(\Delta).
\end{equation}

Set
\begin{equation}
\label{eq:32} 
 \mathfrak{m}'_i =
 \left(
 \frac{t_{\ii_{l-i+1}} \cdots t_{\ii_l}}
 {\xi_{\ii_{1}} \cdots \xi_{\ii_{l-i}}},
 \frac{\xi_{\ii_1} \cdots \xi_{\ii_{l-i+1}}}
 {t_{\ii_{l-i+2}} \cdots t_{\ii_{l}}}
\right).
\end{equation}
This is a maximal ideal of $R'_i$.
This ideal corresponds to the $\bbT$-fixed point 
$p_{\ii_{l-i+1}} \in X'_i$ .

\subsection{Tautological bundles}
\label{sec:tautological-bundles}

In this subsection, we define tautological bundles on the quiver 
variety $\M_\theta(\delta)$ and give their explicit construction.

A tautological bundle is defined as follows. Consider a
vector bundle of rank $l$ on $\mu^{-1}(0)_\theta$ whose fiber
is isomorphic to the representation of $\bar{Q}$ given by 
$(a_i, b_i)_{i=0, \dots, l-1}$ for each point 
$(a_i, b_i)_{i=0, \dots, l-1} \in \mu^{-1}(0)_\theta$.
If the vector bundle descends to a vector bundle 
$\widetilde{\calP}_\theta$ on $\M_\theta(\delta)$, we call
$\widetilde{\calP}_\theta$ a tautological bundle. 

To construct a tautological bundle, consider the matrix algebra 
\[
 M_l(\C[\mu^{-1}(0)]) \simeq M_l(\C) \otimes_\C \C[\mu^{-1}(0)].
\]
The group $GL(\delta)$ acts on $M_l(\C)$ by
\begin{equation}
\label{eq:42} 
 g \cdot E_{ij} = g_i g_j^{-1} E_{ij}
\end{equation}
for $g = (g_k)_{k=0, \dots, l-1} \in GL(\delta)$ and 
$0 \leq i$, $j \leq l-1$.
The group $GL(\delta)$ acts on  $\C[\mu^{-1}(0)]$. 
Thus $GL(\delta)$ acts on $M_l(\C[\mu^{-1}(0)])$.
We define a graded $S$-module
\begin{equation}
\label{eq:47} 
 \calP_\theta = \bigoplus_{m \in \Z_{\geq 0}} e_0 
 M_l(\C[\mu^{-1}(0)])^{GL(\delta), \chi_\theta^m}.
\end{equation}
where $e_0 = E_{00}$.
Let $\widetilde{\calP}_\theta$  be the sheaf
associated to $\calP_\theta$.
We show that $\widetilde{\calP}_\theta$ is 
a direct sum of $l$ line bundles. 
For $i=0$, $\dots$, $l-1$
let $\widetilde{\calL}_i$ be
the sheaf associated to the following graded $S$-module $\calL_i$,
\begin{equation}
\label{eq:19} 
 \calL_i = \bigoplus_{m \in \Z_{\geq 0}} 
 \C[\mu^{-1}(0)]^{GL(\delta), \chi_\theta^m \chi_{\tau_i}}
\end{equation}
where $\tau_i = \epsilon_i - \epsilon_0 \in \Z^l_0$. Consider
the $i$-th column of $\calP_\theta$. The group $GL(\delta)$ acts on
$E_{0i} \in M_l(\C)$ by the character $\chi_{\tau_i}^{-1}$ by
\refeq{eq:42}. Thus the coefficients of the $i$-th column of $\calP_\theta$
coincide with $\calL_i$.
Then, each $\widetilde{\calL}_i$ is a line bundle
on $\M_\theta(\delta)$ and 
we have
\[
 \widetilde{\calP}_\theta = \bigoplus_{i=0}^{l-1} \widetilde{\calL}_i.
\]

Note that ${\calP_\theta}$ is a right module of the matrix algebra
$M_l(\C[\mu^{-1}(0)])^{GL(\delta)}$. Let $A_i$ and $\bar{A}^*_i$ be
the following elements of $M_l(\C[\mu^{-1}(0)])^{GL(\delta)}$:
$A_i = E_{i,i-1} \otimes t_i$ and $\bar{A}^*_i = E_{i-1,i} \otimes \xi_i$.
Then the collection of maps $(A_i, \bar{A}^*_i)_{i=0, \dots, l-1}$ 
gives an action of $\bar{Q}$ on $\widetilde{\calP}_\theta$.
Thus we have the following 
proposition.
\begin{proposition}
 The vector bundle $\widetilde{\calP}_\theta$ is a tautological bundle on
 $\M_\theta(\delta)$.
\end{proposition}

\begin{lemma}
\label{lemma:10}
 The module 
 $\C[\mu^{-1}(0)]^{GL(\delta), \chi_\theta^m \chi_{\tau_i}}$ is a
 $\C[t_0 \cdots t_{l-1}]$-free module.
 \begin{proof}
  Consider the grading of the $\C[t_0 \cdots t_{l-1}]$-module
  $\C[\mu^{-1}(0)]^{GL(\delta), \chi_\theta^m \chi_{\tau_i}}$ defined
  by the degree
  \[
   \deg t_k = 0, \qquad \deg \xi_k = 1.
  \]
  Since $\mu^{-1}(0) \not\subseteq \bigcup_{i=0}^{l-1} \{t_i = 0\}$,
  the $\C[t_0 \cdots t_{l-1}]$-module $\C[\mu^{-1}(0)]$ is torsion free.
  Thus, each component with respect to the above grading is a finitely 
  generated torsion free $\C[t_0 \cdots t_{l-1}]$-module. 
  The algebra $\C[t_0 \cdots t_{l-1}]$ is a one-dimensional polynomial
  algebra. Therefore a finitely generated
  $\C[t_0 \cdots t_{l-1}]$-torsion-free module
  is automatically $\C[t_0 \cdots t_{l-1}]$-free.
 \end{proof}
\end{lemma}

\subsection{Vanishing of higher cohomologies}
\label{sec:line-bundles}

In the previous subsection, we constructed the tautological bundle
$\widetilde{\calP}_\theta$. To calculate the higher cohomologies of 
$\widetilde{\calP}_\theta$, we recall well-known facts about line bundles on
the toric variety $X(\Delta) \simeq \M_\theta(\delta)$.
Let $\Pic(X(\Delta))$
be the Picard group of $X(\Delta)$. Let $D_i$ be the divisor
of $X(\Delta)$ corresponding to $v_i \in N$ for $i=0$, 
$\dots$, $l$ as in \cite[Sec 3.3]{Fu}. Under the isomorphism \refeq{eq:18},
$D_i$ corresponds to $\U_{\ii_i}$ defined by \refeq{eq:41} 
for $i=1$, $\dots$, $l$.
By the general theory of toric varieties, we have
the following lemma.

\begin{lemma}[\cite{Mu} (2.3), \cite{Fu} Prop 3.4]
The Picard group $\Pic(X(\Delta))$ is generated by the 
divisors $D_0$, $\dots$, $D_l$. Moreover their relations 
in $\Pic(X(\Delta))$ are given by:
\begin{align*}
 D_0 + D_1 + \dots + D_l &= 0,  \\
 \sum_{i=1}^{l} i D_i &= 0.
\end{align*}
\end{lemma}
For $i=1$, $\dots$, $l-1$,  we define the cycle
$D(i) = \sum_{j=0}^{i-1} (i-j) D_{l-j} \in \Pic(X(\Delta))$.
We have $\Pic(X(\Delta)) = \bigoplus_{i=1}^{l-1} \Z D(i)$.

For $b = (b_i)_{i=1, \dots, l-1} \in \Z^{l-1}$, let
$D(b) = \sum_{i=1}^{l-1} b_i D(i) \in \Pic(X(\Delta))$ as
in \cite{Mu}. For each divisor $D \in \Pic(X(\Delta))$, 
we have the $\bbT$-invariant line bundle $\calO(D)$ on $X(\Delta)$.
The following two lemmas are proved in \cite{Mu}.

\begin{lemma}[\cite{Mu}, Lemma 2.4]
 When $b = (b_k)_{k=1, \dots, l-1} \in \Z^{l-1}_{\geq 0}$, 
 the space of local sections $H^0(X_j, \calO(D(b)))$ is a
 free $R_j$-module generated by the element
 $x^{\sum_{k=l-j+1}^{l-1} (l-k)b_k} y^{-\sum_{k=l-j+1}^{l-1} k b_k}$.
\end{lemma}

\begin{proposition}[\cite{Mu}, Lemma 2.1]
 \label{prop:4}
 If $b = (b_k)_{k=1, \dots, l-1} \in \Z^{l-1}_{\geq 0}$,
 then we have
 \[
 H^p(X(\Delta), \calO(D(b)))
 = 0  \]
for $p \ne 0$.
\end{proposition}

For $b = (b_k)_{k=1, \dots, l-1} \in (\Z_{\geq 0})^{l-1}$, let
$\calO'(D(b))$ be the line bundle on $X(\Delta)$ such that 
$H^0(X_j, \calO'(D(b)))$ is the free $R_j$-module generated
by the element \break
$x^{\sum_{k=l-j+1}^{l-1} (l-k)b_k} y^{\sum_{k=1}^{l-j} k b_k}$.
Namely, as a $\bbT$-equivariant line bundle
\[
 \calO'(D(b)) \simeq \calO(D(b)) \otimes_{\C}
 \C_{(0, \sum_{k=1}^{l-1} k b_k)}
\]
where $\C_{(a,b)}$ is the one-dimensional vector space with
the $\bbT$-action of weight $(a,b)$.

Fix a stability parameter 
$\theta = (\theta_i)_{i=0, \dots, l-1} \in \Z^l_{reg}$.
In \refsec{sec:tautological-bundles}, we
defined the line bundle $\widetilde{\calL}_i$ on the quiver 
variety $\M_\theta(\delta)$ for $i=0$, $\dots$, $l-1$ by 
\refeq{eq:19}. 

For $i=0$, $\dots$, $l-1$ and $m \in \Z_{\geq 0}$, set
\begin{equation}
\label{eq:20} 
 b_k^{\theta'} = \theta'_{\ii_k} + \theta'_{\ii_k+1}
 + \dots + \theta'_{\ii_{k+1}-1}.
\end{equation}
where $\theta' = (\theta'_k)_{k=0, \dots, l-1} = m \theta + \tau_i 
\in \Z^l_0$.
Note that we have $b_k^{\theta'} \in \Z_{\geq 0}$ for all
$k$.

For
$\theta' = (\theta'_k)_{k=0, \dots, l-1} 
= m \theta + \tau_i \in \Z^l$ where $m \in \Z_{\geq 0}$
and $i=0$, $\dots$, $l-1$, 
set
\begin{equation}
\label{eq:12} 
  f^{\theta'}_j = 
\prod_{k=j}^{l-1} (t_{\ii_{k+1}} \cdots t_{\ii_l})^{b_k^{\theta'}}
 \prod_{k=1}^{j-1} (\xi_{\ii_1} \cdots \xi_{\ii_k})^{b_k^{\theta'}}.
\end{equation}
Note that $f_j^{\theta'}$ does not vanish at the fixed point $p_{\ii_j}$
by \refprop{prop:18}.

We show that $f^{\theta'}_j$ belongs to 
$\C[\mu^{-1}(0)]^{GL(\delta), \chi_{\theta'}}$. 
 We calculate the weight of the function $f^{\theta'}_j$. Because 
 the weight of $\xi_{\ii_1} \cdots \xi_{\ii_k}$ is equal to
 the weight of $t_{\ii_{k+1}} \cdots t_{\ii_l}$ for all $k$,
 the weight of $f^{\theta'}_j$ is independent of 
 $j$. By $\theta'_0 + \dots + \theta'_{l-1} = 0$, we have
\begin{align*}
 \lefteqn{b_k^{\theta'} + b_{k+1}^{\theta'} + \dots + b_{l-1}^{\theta'}}
 & \\
&= (\theta'_{\ii_k} + \dots + \theta'_{\ii_{k+1} - 1}) +
 (\theta'_{\ii_{k+1}} + \dots + \theta'_{\ii_{k+2} - 1}) + \dots +
 (\theta'_{\ii_{l-1}} + \dots + \theta'_{\ii_{l} - 1}) \\
&=
 \theta'_{\ii_k} + \theta'_{\ii_k + 1} + \dots + \theta'_{\ii_l - 1}
\end{align*}
Thus, we have
\[
 f^{\theta'}_l 
 = \prod_{k=1}^{l-1} (\xi_{\ii_1} \cdots \xi_{\ii_k})^{b_k^{\theta'}}
 = \prod_{k = 0, \dots, l-1}
 \xi_k^{\theta'_{k} + \theta'_{k+1} + \dots + \theta'_{\ii_l-1}}.
\]
Thus the weight of $f^{\theta'}_l$ is
\begin{align*}
 \lefteqn{\sum_{k=0}^{l-1}\{ (\theta'_{k} + \dots + \theta'_{\ii_l-1}) 
\epsilon_{k} - 
(\theta'_{k} + \dots + \theta'_{\ii_l-1}) \epsilon_{k-1} \}} & \\
 &= \sum_{k=0}^{l-1} \{
 (\theta'_k + \dots + \theta'_{\ii_l-1}) - 
 (\theta'_{k+1} + \dots + \theta'_{\ii_l-1})\} \epsilon_k
 = \sum_{k=0}^{l-1} \theta'_k \epsilon_k.
\end{align*}
 Therefore, $f^{\theta'}_j$ belongs to 
$\C[\mu^{-1}(0)]^{GL(\delta), \chi_{\theta'}}$ for all $j$. 

\begin{lemma}
 \label{lemma:5}
 For $1 \leq j \leq l$ and $0 \leq i \leq l-1$ and
 $m \in \Z_{> 0}$, $H^0(X'_j, \widetilde{\calL}_i \otimes \catO(m))$
 is the free $R'_j$-module with the generator $f^{\theta'}_{l+1-j}$ 
 with $\theta' = m \theta + \tau_i \in \Z^l_0$.
\begin{proof}
By the definition of $\widetilde{\calL}_i$ \refeq{eq:19}, 
Laurent monomials of $t_0$, $\dots$, $t_{l-1}$, $\xi_0$, $\dots$,
 $\xi_{l-1}$ with the $GL(\delta)$-character $\chi_{\theta'}$
span $H^0(X'_j, \widetilde{\calL}_i \otimes \calO(m))$ over $\C$.

Since the unique fixed point $p_{\ii_{l-j+1}} \in X'_j$  corresponds
to the maximal ideal $\mathfrak{m}'_j \subset R'_j$ defined by
\refeq{eq:32}, the generators of $H^0(X'_j, \widetilde{\calL}_i \otimes 
\calO(m))$ must have no zero at $p_{\ii_{l-j+1}}$. 

The Laurent monomial $f^{\theta'}_j$ has no zero and no pole at
$p_{\ii_{l-j+1}}$ by \refprop{prop:18}. On the other hand, let
$g$ be a Laurent monomial with the $GL(\delta)$-character 
$\chi_{\theta'}$. Then, $g / f^{\theta'}_j$ is a product of
the following Laurent monomials:
\begin{gather*}
 (t_0 \cdots t_{l-1})^{\pm 1}, \quad  (\xi_0 \cdots \xi_{l-1})^{\pm 1}, 
 \quad (t_0 \xi_0)^{\pm 1},\\
 \left(\frac{t_p t_{p+1} \cdots t_{q-1}}{\xi_q \xi_{q+1}\cdots
 \xi_{p-1}}\right)^{\pm 1} \qquad \text{for $p \ne q$}.
\end{gather*}
Therefore, $g$ has either 
zeros or poles at $p_{\ii_{l-j+1}}$ by \refprop{prop:18}.

Therefore, the Laurent monomials other than
$f_j^{\theta'}$ have either zeros or poles at $p_{\ii_{l-j+1}}$, and
$f_j^{\theta'}$ is  the generator of
 $H^0(X'_j, \widetilde{\calL}_i \otimes \calO(m))$ over
$R'_j$.
\end{proof}
\end{lemma}

We identify $\M_\theta(\delta)$ and $X(\Delta)$ 
by the isomorphism \refeq{eq:18}.

\begin{proposition}
 \label{prop:5}
 For $i=0$, $\dots$, $l-1$ and $m \in \Z_{> 0}$, we have
 an isomorphism of $\bbT$-equivariant line bundles on 
 $\M_\theta(\delta) \simeq X(\Delta)$: 
 $\widetilde{\calL}_i \otimes \calO(m) \simeq \calO'(D(b))$
 where $b = (b_j^{\theta'})_{j=1, \dots, l-1}$ and
 $\theta' = m \theta + \tau_i$.
 \begin{proof}
  For $j=1$, $\dots$, $l$, 
  $H^0(X'_j, \widetilde{\calL}_i \otimes \calO(m))$ is a free
  $R'_j$-module with the generator $f^{\theta'}_{l-j+1}$. On the
  other hand $H^0(X_j, \calO'(D(b)))$ is a free $R_j$-module
  with the generator 
  $x^{\sum_{k=l-j+1}^{l-1} (l-k) b^{\theta'}_k}
  y^{\sum_{k=1}^{l-j} k b^{\theta'}_k}$. Thus the map given by 
  $t_k \mapsto x$, $\xi_k \mapsto y$ is a $\bbT$-equivariant
  isomorphism 
  $\widetilde{\calL}_i \otimes \calO(m) \simeq \calO'(D(b))$.
 \end{proof}
\end{proposition}

By the general theory of toric varieties, we have the 
following $\C$-basis
of \break
$H^0(X(\Delta), \calO(D(b)))$ for $b = (b_i)_{i=1, \dots, l-1}
\in \Z^{l-1}_{\geq 0}$.
As in \cite[page 66]{Fu}, set
\[
 P_{D(b)} = \{ m \in M \; | \;
 \langle m, v_i \rangle \geq - a_i \quad \text{for $i=0$, $\dots$, $l$}.\}
\]
where $a_i = b_{l-i+1} + 2 b_{l-i+2} + \dots + (i-1) b_{l-1}$.
Then we have
\[
 H^0(X(\Delta), \calO(D(b))) 
 = \bigoplus_{m \in P_{D(b)} \cap M}
 \C x^{m_1} y^{m_1 + l m_2}.
\]
The following lemma is an immediate consequence of this
fact.

\begin{lemma}
 \label{lemma:2}
 As a $\C[x^l, xy]$-module, we have 
 \[
  H^0(X(\Delta), \calO'(D(b))) =
 \sum_{i=1}^{l} \sum_{m=0}^{b_i - 1}
 \C[x^l, xy] x^{\sum_{j=i+1}^{l} (l-j) b_j + (l-i) (b_i - m)}
 y^{\sum_{j=1}^{i-1} j b_j + i m}
 \]
 where we set $b_l = \infty$.
\end{lemma}

Now we have the following proposition.

\begin{proposition}
 \label{prop:3}
 For $i=0$, $\dots$, $l-1$ and $m \in \Z_{> 0}$, 
we have
\[
 H^p(\M_\theta(\delta), \widetilde{\calL}_i \otimes \calO(m))
 = 
 \left\{
 \begin{array}{ll}
  \C[\mu^{-1}(0)]^{GL(\delta), \chi_\theta^m \chi_{\tau_i}} &
   (p = 0),\\
  0 &  (p \ne 0).
 \end{array}
 \right.
\]
\begin{proof}
 Set $b = (b^{\theta'}_j)_{j=1, \dots, l-1}$
with $\theta' = m \theta + \tau_i$.
By \refprop{prop:5}, we have an isomorphism of line bundles 
$\widetilde{\calL}_i \otimes \calO(m) \simeq \calO'(D(b))$. 
Therefore, by \refprop{prop:4}, we have the
 vanishing of the higher cohomologies
\[
  H^p(\M_\theta(\delta), \widetilde{\calL}_i \otimes \calO(m)) = 0
\]
 for $p \ne 0$. By the definition of $\widetilde{\calL}_i$
 at \refeq{eq:19}, it is clear that 
 \[
 \C[\mu^{-1}(0)]^{GL(\delta), \chi_\theta^m \chi_{\tau_i}}
 \subseteq H^0(\M_\theta(\delta), \widetilde{\calL}_i \otimes \calO(m)).
 \]
 We show the opposite inclusion. By \reflemma{lemma:2}, we have
\begin{equation}
 \label{eq:1}
  H^0(X(\Delta), \calO'(D(b))) =
 \sum_{k=1}^{l} \sum_{n=0}^{b^{\theta'}_k - 1}
 \C[x^l, xy] x^{\sum_{j=k+1}^{l} (l-j) b^{\theta'}_j + (l-k) 
 (b^{\theta'}_k - n)}
 y^{\sum_{j=1}^{k-1} j b^{\theta'}_j + k n}.
\end{equation}
On the other hand,
 we consider the elements 
 $g^{\theta'}_k(n) 
 \in \C[\mu^{-1}(0)]^{GL(\delta), \chi_\theta^m \chi_{\tau_i}}$
 \[
  g^{\theta'}_k(n) = (t_{\ii_{k+1}} \cdots t_{\ii_l})^{b_k^{\theta'} - n}
 \prod_{j=k+1}^{l-1} 
 (t_{\ii_{j+1}} \cdots t_{\ii_l})^{b_j^{\theta'}}
 (\xi_{\ii_1} \cdots \xi_{\ii_k})^n 
 \prod_{j=1}^{k-1} (\xi_{\ii_1} \cdots \xi_{\ii_j})^{b_j^{\theta'}}
 \]
 for $k = 1$, $\dots$, $l$ and $n = 0$, $\dots$, $b_k^{\theta'} - 1$.
 Here we set $b_l^{\theta'} = \infty$.
 The homomorphism given by $t_j \mapsto x$, $\xi_j \mapsto y$ maps 
 the elements $g^{\theta'}_k(n)$ to the generators in \refeq{eq:1}. 
The isomorphism $\widetilde{\calL}_i \otimes \calO(m) \simeq \calO'(D(b))$
implies that 
 $H^0(\M_\theta(\delta), \widetilde{\calL}_i \otimes \calO(m))$ 
 is isomorphic to
 $H^0(X(\Delta), \calO'(D(b)))$. Thus we have
 \[
 \C[\mu^{-1}(0)]^{GL(\delta), \chi_\theta^m \chi_{\tau_i}}
 = H^0(\M_\theta(\delta), \widetilde{\calL}_i \otimes \calO(m)).
 \]
\end{proof}
\end{proposition}

\begin{corollary}
 \label{cor:1}
 For $m \in \Z_{> 0}$, we have
\[
 H^p(\M_\theta(\delta), \widetilde{\calP}_\theta \otimes \calO(m))
 = 
 \left\{
 \begin{array}{ll}
  e_0 M_l(\C[\mu^{-1}(0)])^{GL(\delta), \chi_\theta^m} &
   (p = 0),\\
  0 &  (p \ne 0).
 \end{array}
 \right.
\]
\end{corollary}

\subsection{$(q, t)$-dimension}
\label{sec:q-t-dimension}

Let $V$ be a possibly infinite-dimensional vector space equipped
with an action of the torus $\bbT = (\C^*)^2$. Then, we have the weight space
decomposition of $V$: $V = \bigoplus_{r,s} V_{r,s}$ where
$V_{r,s}$ is the weight space which belongs to the weight
\begin{align*}
 \rho_{r,s} : \bbT &\longrightarrow \C \\
 (m_1, m_2) &\mapsto m_1^r m_2^s
\end{align*}
The $(q, t)$-dimension of $V$ is the following formal series:
\[
 \Dim_{q,t} V = \sum_{r,s} (\Dim V_{r,s}) q^r t^s
\]
The torus $\bbT$ acts on $\Rep(\bar{Q}, \delta)$. The action 
induces an action of $\bbT$ on $\C[\mu^{-1}(0)]$. The 
weight spaces with respect to this action are equal to the 
homogeneous spaces with respect to the following bi-grading 
on $\C[\mu^{-1}(0)]$. 
\[
 \deg t_i = (1,0), \qquad \deg \xi_i = (0,1),
\]
for $i=0$, $\dots$, $l-1$.
We consider the $(q, t)$-dimensions of 
$e_0 M_l(\C[\mu^{-1}(0)])^{GL(\delta), \chi_\theta^m}$ for 
$m \in \Z_{\geq 0}$
with respect to this
action. 

To calculate the $(q,t)$-dimension of 
$e_0 M_l(\C[\mu^{-1}(0)])^{GL(\delta), \chi_\theta^m}$ for
$m \in \Z_{> 0}$, we use the following
 Atiyah-Bott-Lefschetz formula
together with \refcor{cor:1}.

\begin{theorem}[\cite{Ha1} Theorem 3.1]
 \label{thm:1}
 Let $X$ be a smooth surface equipped with an action of 
 $\bbT$, and assume the fixed point set $X^\bbT$ is finite.
 Let $\calF$ be a $\bbT$-equivariant locally free sheaf on $X$.
 For $x \in X^\bbT$, $\bbT$ acts on $\calF(x)$. Suppose that
 $\bbT$ acts on the cotangent space at $x$ with weights
 $(v_1, v_2)$ and $(w_1, w_2)$. Then we have,
\[
 \sum_{p \geq 0} (-1)^p \dim_{q,t} H^p(X, \calF)
 = \sum_{x \in X^\bbT}
 \frac{\dim_{q,t} \calF(x)}{(1-q^{v_1} t^{v_2})(1 - q^{w_1}t^{w_2})}.
\]
\end{theorem}

We will apply the above theorem for $X = \M_\theta(\delta)$
and $\calF = \widetilde{\calP}_\theta \otimes \calO(m)$. 
Then we have the following theorem.

\begin{theorem}
 \label{thm:3}
 For $m \geq 0$, $i=0$, $\dots$, $l-1$, we have 
 the following identity
\begin{multline}
\label{eq:22}
  \dim_{q,t} \C \otimes_{\C[t_0 \cdots t_{l-1}]}
 e_0 M_l(\C[\mu^{-1}(0)])^{GL(\delta), \chi_\theta^m}
 \Bigr|_{t=q^{-1}} \\
 = 
 \sum_{i=1}^{l} q^{d_i^{m \theta} + (l-i)} \frac{1}{1-q^{-1}} 
\end{multline}
 where 
\begin{equation}
\label{eq:21} 
 d_i^\theta = - \theta_0 - 2\theta_1 - \dots - i \theta_{i-1}
 + (l-i-1) \theta_i + \dots + \theta_{l-2}.
\end{equation}
for $i=0$, $\dots$, $l-1$. Here we set $d_l^{\theta} = d_0^\theta$.
\end{theorem}

This will be proved in the next subsection. 

\subsection{The proof of the theorem}
\label{sec:proof-theorem}

In this subsection, we prove \refthm{thm:3}. To prove the theorem,
we apply \refthm{thm:1} for $\calF = \widetilde{\calP}_\theta \otimes \calO(m)$
together with \refcor{cor:1} for $m \geq 1$. Since we did not prove the
vanishing of the cohomologies \refcor{cor:1} for $m=0$, we cannot 
make use of \refthm{thm:1} in this case. On the other hand, 
in the case of
$m=0$, the space $\C \otimes_{\C[t_0 \cdots t_{l-1}]} 
e_0 M_l(\C[\mu^{-1}(0)])^{GL(\delta)}$ is independent of
the stability parameter $\theta = (\theta_i)_{i=0, \dots, l-1} 
\in \Z^l_{reg}$. Therefore we can easily show \refeq{eq:22}
by a direct calculation.

Set $\bar{A}^* = \bar{A}^*_0 + \bar{A}^*_1 + \dots + \bar{A}^*_{l-1}$.
By the right action on $\calP_\theta$ given in 
\refsec{sec:tautological-bundles}, $e_0
M_l(\C[\mu^{-1}(0)])^{GL(\delta)}$ is a right $\C[\bar{A}^*]$-module.
As the right $\C[\bar{A}^*]$-module, we have
\[
 \C \otimes_{\C[t_0 \cdots t_{l-1}]} e_0 M_l(\C[\mu^{-1}(0)])^{GL(\delta)}
 = \bigoplus_{i=0}^{l-1} E_{0i} \otimes t_{i+1} \cdots t_{l-1} t_0 \C[\bar{A}^*].
\]
Therefore, we have
\[
 \dim_{q,t} \C \otimes_{\C[t_0 \cdots t_{l-1}]}
 e_0 M_l(\C[\mu^{-1}(0)])^{GL(\delta)}  \Bigr|_{t=q^{-1}}
 = \sum_{i=1}^{l} q^{l-i} \frac{1}{1-q^{-1}}.
\] 

Consider the case of $m \in \Z_{> 0}$. Since, at the fixed point
$p_i \in \M_\theta(\delta)$ for $i=0$, $\dots$, $l-1$,
\[
 \dim_{q,t} \left(\widetilde{\calP}_\theta \otimes \calO(m)\right)(p_i)
 = \dim_{q,t} \widetilde{\calP}_\theta(p_i) \cdot 
 \left(\dim_{q,t} \calO(1)(p_i) \right)^m,
\]
we need to
calculate $\dim_{q,t} \calO(1)(p_i)$ and $\dim_{q,t} \widetilde{\calP}_\theta(p_i)$.

First we consider the fibers of $\calO(1)$ at the fixed 
points $p_{\ii_1}$, $\dots$, $p_{\ii_l}$.

\begin{lemma}
 \label{lemma:3}
 We have $\calO(1)_{p_{\ii_i}} = \calO_{p_{\ii_i}} f^\theta_i$
 where $f^\theta_i$ is the function defined on \refeq{eq:12}
 with $\theta' = \theta$, i.e., 
\[
 f^\theta_i = \prod_{j=i}^{l-1} (t_{\ii_{j+1}} \cdots t_{\ii_l})^{b_j^\theta}
 \prod_{j=1}^{i-1} (\xi_{\ii_1} \cdots \xi_{\ii_j})^{b_j^\theta}.
\]
\begin{proof}
 Apply \reflemma{lemma:5} for $i=0$ and $m=1$, we have
 $\widetilde{\calL}_0 \otimes \calO(1) = \calO(1)$, thus we have
 \[
  H^0(X'_i, \calO(1)) = R'_i f^\theta_i.
 \]
 The fixed point $p_{\ii_{l-i+1}} \in X'_i$ corresponds to the maximal
 ideal $\mathfrak{m}'_i \subset R'_i$. Therefore we have
\[
 \calO(1)_{p_{\ii_i}} = \calO_{p_{\ii_i}} f^\theta_i.
\]
\end{proof}
\end{lemma}

\begin{corollary}
 \label{cor:2}
 For $i=1$, $\dots$, $l$, we have 
\[
 \dim_{q,t} \calO(1)(p_{\ii_i}) \big|_{t=q^{-1}}
 = q^{d_{\ii_i}^\theta}.
\]
\begin{proof}
 First we calculate the case of $i=1$, then we have
\begin{align*}
 f^\theta_1 &= \prod_{j=1}^{l-1} (t_{\ii_{j+1}} \cdots t_{\ii_l})^{b_j^\theta} \\
 &= t_{\ii_1+1}^{\theta_{\ii_1}} 
 t_{\ii_1 + 2}^{\theta_{\ii_1} + \theta_{\ii_1+1}} \cdots
 t_{\ii_1 + l - 1}^{\theta_{\ii_1} + \theta_{\ii_1+1} + \dots +
 \theta_{\ii_1 + l - 2}}.
\end{align*} 
Consider the degree given by $\deg t_i = 1$ and $\deg \xi_i = -1$ for
$i=0$, $\dots$, $l-1$.
The degree of $f^\theta_1$ is 
\[
 \deg f^\theta_1 = (l-1) \theta_{\ii_1} + (l-2) \theta_{\ii_1+1}
 + \dots + \theta_{\ii_1+l-2} 
 = d_{\ii_1}^\theta.
\]
Thus the statement of the corollary is valid for $i=1$. On the 
other hand, we have
$d_{\ii_i}^\theta - d_{\ii_{i+1}}^\theta = l b_i^\theta$
and $\deg f^\theta_i - \deg f^\theta_{i+1} = l b_i^\theta$. Therefore we have
$\deg f^\theta_i = d_{\ii_i}^\theta$ for $i=2$, $\dots$, $l-1$ by 
induction on $i$. 
\end{proof}
\end{corollary}

Next we consider the fibers of $\widetilde{\calP}_\theta$ at the fixed points.
At the fixed  point $p_i = [a_j, b_j]_{j=0, \dots, l-1}$ defined by
\refeq{eq:34}, we consider the stalk 
$(\widetilde{\calP}_\theta)_{p_i} 
= \bigoplus_{k=0}^{l-1} \C E_{0k} \otimes (\widetilde{\calL}_k)_{p_i}$.
For $k = 0$, $\dots$, $l-1$,
let $v_k$ be the germ of $(\widetilde{\calP}_\theta)_{p_i}$
defined as follows.
\begin{align*}
 v_k &= E_{0k} \otimes \nu_1 \cdots \nu_{k}  \qquad \text{if $k \leq i-1$}, \\
 v_k &= E_{0k} \otimes \nu'_{k+1} \cdots \nu'_{l-1} \nu'_{0}  \qquad
\text{if $k \geq i$}
\end{align*}
where
\begin{align*}
 \nu_j &= \left\{
  \begin{array}{ll}
   t_j^{-1} & (\text{if $a_j \ne 0$}) \\
   \xi_j & (\text{if $b_j \ne 0$}) \\
  \end{array}\right. , \\
 \nu'_j &= \left\{
  \begin{array}{ll}
   t_j & (\text{if $a_j \ne 0$}) \\
   \xi_j^{-1} & (\text{if $b_j \ne 0$}) \\
  \end{array}
 \right. .
\end{align*}
Note that
the group $GL(\delta)$ acts on $v_k$ by the character $\chi_{\tau_k}$,
and $v_k$ is not zero at $p_i$.

The stalk $(\widetilde{\calP}_\theta)_{p_i}$ is
a free $\calO_{p_i}$-module of
rank $l$. As
the $v_k$ above for $k=0$, $\dots$, $l-1$ are clearly
linearly independent, we have the following lemma.

 \begin{lemma}
  For $i=0$, $\dots$, $l-1$,
  the stalk $(\widetilde{\calP}_\theta)_{p_i}$ has the 
  $\calO_{p_i}$-basis $\{v_k\}_{k=0}^{l-1}$.
 \end{lemma}

\begin{corollary}
 \label{cor:3}
 For $i=0$, $\dots$, $l-1$, we have
 \[
  \dim_{q,t} \widetilde{\calP}_\theta(p_i) \Bigr|_{t=q^{-1}} = 
 q^{l-i} \frac{1-q^{-l}}{1-q^{-1}}.
 \]
\end{corollary}

Finally, the cotangent space of $\M_\theta(\delta)$ at the
fixed points has the following well-known structure.

\begin{lemma}[\cite{Mu}]
 \label{lemma:4}
 For $i=1$, $\dots$, $l$, the cotangent space of $\M_\theta(\delta)$
 at $p_{\ii_i}$ has $\bbT$-action with weights 
 \begin{align*}
  (v_1, v_2) &= (l-i, -i), \\
  (w_1, w_2) &= (-l+i+1, i+1).
 \end{align*}
\end{lemma}

Now we apply \refthm{thm:1} for $\calF = \widetilde{\calP}_\theta \otimes \calO(m)$
with $m \in \Z_{>0}$ to prove \refthm{thm:3}. 

By \reflemma{lemma:10} and \refcor{cor:1}, we have the modules
$H^p(\M_\theta(\delta), \widetilde{\calP}_\theta \otimes \calO(m))$ are
$\C[t_0 \cdots t_{l-1}]$-free. Thus we have
\begin{equation}
\label{eq:33} 
 \dim_{q,t} \C \otimes_{\C[t_0 \cdots t_{l-1}]}
 H^p(\M_\theta(\delta), \widetilde{\calP}_\theta \otimes \calO(m))
 = (1 - q^l) 
 \dim_{q,t}
 H^p(\M_\theta(\delta), \widetilde{\calP}_\theta \otimes \calO(m)).
\end{equation}

By \refcor{cor:1}, 
\begin{multline*}
 \dim_{q,t} \C \otimes_{\C[t_0 \cdots t_{l-1}]}
 e_0 M_l(\C[\mu^{-1}(0)])^{GL(\delta), \chi_\theta^m} \\
 = \sum_p (-1)^p \dim_{q,t}
 \C \otimes_{\C[t_0 \cdots t_{l-1}]}
 H^p(\M_\theta(\delta), \widetilde{\calP}_\theta \otimes
 \calO(m)).
\end{multline*}
By \refthm{thm:1} together with \refcor{cor:2}, \refcor{cor:3},
\reflemma{lemma:4} and \refeq{eq:33}, we have
\begin{align*}
\lefteqn{\sum_p (-1)^p \dim_{q,t}
 \C \otimes_{\C[t_0 \cdots t_{l-1}]}
 H^p(\M_\theta(\delta), \widetilde{\calP}_\theta \otimes
 \calO(m)) \Bigr|_{t=q^{-1}}} &
 \\
&= \sum_{i=1}^l (1-q^l) 
\frac{\dim_{q,t} \widetilde{\calP}_\theta(p_{\ii_i}) \dim_{q,t}
 \calO(m)(p_{\ii_i})}
{(1-v_{p_{\ii_i}}(q,t))(1-w_{p_{\ii_i}}(q,t))} \Bigr|_{t=q^{-1}} \\
&= \sum_{i=1}^l (1-q^l) 
\frac{q^{d_{\ii_i}^{m\theta}} q^{l-\ii_i} (1-q^{-l})/(1-q^{-1}) }
{(1-q^l)(1-q^{-l})} \\
&= \sum_{i=1}^l q^{d_{i}^{m\theta} + l-i} \frac{1}{1-q^{-1}}.
\end{align*}
Thus we have \refeq{eq:22} of \refthm{thm:3}. 

\section{Rational Cherednik algebras}
\subsection{Rational Cherednik algebras}
\label{sec:rati-cher-algebr}

First we introduce the rational Cherednik algebra 
$H_{\kappa} = H_{\kappa}(\Z_l)$ for the group $\Z_l = \Z / l \Z$
with a parameter $\kappa = (\kappa_0, \dots, \kappa_{l-1}) \in \R^l$.
As a vector space, $H_\kappa$ is given by
\[
 H_\kappa = \C[x] \otimes \C\Z_l \otimes \C[y].
\]
The relations of $H_\kappa$ are as follows:
\begin{align*}
 \gamma x \gamma^{-1} &= \zeta^{-1} x, \\
 \gamma y \gamma^{-1} &= \zeta y, \\
 [y, x] &= 1 + l \sum_{i=0}^{l-1} (\kappa_{i+1} - \kappa_i) \bar{e}_i.
\end{align*}
Here we set $\kappa_l = \kappa_0$ and $\bar{e}_i$ is an idempotent 
defined in \refsec{sec:basic-notations}.
Note that $H_\kappa$ depends only on $\kappa_{i+1} - \kappa_i$
for $i=0$, $\dots$, $l-1$. We have
\[
 x \bar{e}_i = \bar{e}_{i+1} x, \qquad
 y \bar{e}_i = \bar{e}_{i-1} y.
\]

The polynomial algebras $\C[x]$ and $\C[y]$ are subalgebras of $H_\kappa$.
Moreover, the smash products $\C[x] \smashprod \Z_l$ and 
$\C[y] \smashprod \Z_l$ are subalgebras of $H_\kappa$.
We also define the 
spherical subalgebra $U_\kappa$ of $H_\kappa$ as 
$U_\kappa = \bar{e}_0 H_\kappa \bar{e}_0$.

In \cite{DO}, the following homomorphism of algebras 
from $H_\kappa$ into the algebra 
$D(\C^{*}) \smashprod \Z_l$ was defined,
\begin{align*}
 H_\kappa &\longrightarrow D(\C^*) \smashprod \Z_l, \\
 x &\mapsto x, \\
 \gamma &\mapsto \gamma, \\
 y &\mapsto D_y = 
 \frac{d}{dx} + \frac{l}{x} \sum_{i=0}^{l-1} \kappa_i \bar{e}_i.
\end{align*}
This homomorphism is injective. This map is called the Dunkl-Cherednik
embedding, and the operator $D_y$ is called the Dunkl operator.
Let $\catO_\kappa$ be the subcategory of $\fmod{H_\kappa}$ 
such that $y \in H_\kappa$ acts locally nilpotently on 
objects of $\catO_\kappa$. By \cite{DO} and \cite{GGOR}, $\catO_\kappa$ is a 
highest weight category with index poset 
$\Lambda = \{0, 1, \dots, l-1\}$ in the sense of \cite{CPS}.

We define the standard modules of $\catO_\kappa$. For 
$i = 0$, $\dots$, $l-1$, we have the irreducible $\Z_l$-modules
$L_i = \C \hwvec_i$. Let $y$ act trivially on $L_i$,
so this induces an action of the algebra $\C[y] \smashprod \Z_l$. 
The algebra $\C[y] \smashprod \Z_l$ is a subalgebra of $H_\kappa$,
thus we define
the standard module $\Delta_\kappa(i)$ as the induced module
\[
 \Delta_\kappa(i) = H_\kappa \otimes_{\C[y] \smashprod \Z_l}
 L_i.
\]
The following proposition is due to \cite{DO} and \cite{GGOR}.

\begin{proposition}
\label{prop:9}
(1) For each $i=0$, $\dots$, $l-1$,
the standard module $\Delta_{\kappa}(i)$ has a unique 
simple quotient which we denote by $L_\kappa(i)$.

(2) For any simple object $L \in \catO_\kappa$, we have
 an isomorphism $L \simeq L_\kappa(i)$  for some
$i=0$, $\dots$, $l-1$.
\end{proposition}

\subsection{Deformed preprojective algebras}
\label{sec:deform-prepr-algebr}

Deformed preprojective algebras were first introduced by
\cite{CBH}. We use another equivalent definition which was defined by
\cite{Ho}.

As in \refsec{sec:quivers} and \refsec{sec:defin-quiv-vari}, 
let $Q = (I, E)$ be  the quiver of type $A_{l-1}^{(1)}$. 

We consider the space of representations
\[
 \Rep(Q, \delta) = \{ (a_i)_{i=0, \dots, l-1} \;|\; a_i \in \C \}
 \simeq \C^l
\]
with the dimension vector $\delta = (1, \dots, 1)$. Consider
the algebra $D(\Rep(Q, \delta))$ of algebraic differential operators.
Let $t_0$, $\dots$, $t_{l-1} \in \C[\Rep(Q, \delta)]$ be 
the coordinate functions such that $t_i((a_j)_{j=0, \dots, l-1}) = a_i$.
Set $\partial_i = \partial / \partial t_i$. 
The algebra $D(\Rep(Q, \delta))$ is
generated by $t_0$, $\dots$, $t_{l-1}$, $\partial_0$, $\dots$,
$\partial_{l-1}$.

The group $GL(\delta)$ acts on $D(\Rep(Q, \delta))$ 
by
\begin{align*}
 g \cdot t_i &= g_i^{-1} g_{i-1} t_i, \\
 g \cdot \partial_i &= g_i g_{i-1}^{-1} \partial_i 
\end{align*}
for $i=0$, $\dots$, $l-1$ and
$g = (g_k)_{k=0, \dots, l-1} \in GL(\delta)$.
The 
action of $GL(\delta)$ on $\Rep(Q, \delta)$ induces a 
homomorphism of Lie algebras
\[
 \iota : \gl(\delta) \longrightarrow D(\Rep(Q, \delta))^{GL(\delta)}.
\]
As a Lie algebra,  
$\gl(\delta) = \bigoplus_{i=0}^{l-1} \gl_1(\C) \simeq 
\bigoplus_{i=0}^{l-1} \C e^{(i)}$
where $e^{(i)}$ is a natural basis of the $i$-th component. 
Then, we have
\[
 \iota(e^{(i)}) = t_{i+1} \partial_{i+1} - t_{i} \partial_{i}.
\]

Consider the $l \times l$ matrix algebra $M_l(D(\Rep(Q, \delta)))$.
We have an isomorphism 
\begin{equation}
\label{eq:43} 
 M_l(D(\Rep(Q, \delta))) \simeq M_l(\C) \otimes_{\C} 
D(\Rep(Q, \delta)).
\end{equation}
The group $GL(\delta)$ acts on $M_l(\C)$ by \refeq{eq:42}.
It also acts on $D(\Rep(Q, \delta))$.
Thus $GL(\delta)$ acts diagonally on $M_l(D(\Rep(Q, \delta)))$ through
the isomorphism \refeq{eq:43}. 

We have the following
homomorphism of Lie algebras:
\begin{gather*}
 \tau : \gl(\delta) \longrightarrow M_l(D(\Rep(Q, \delta)))^{GL(\delta)}, \\
 \tau = \varpi \otimes 1 + 1 \otimes \iota,
\end{gather*}
where $\varpi : \gl(\delta) \rightarrow M_l(\C)$ is given by
$\varpi(e^{(i)}) = E_{ii}$.

For a parameter $\lambda = (\lambda_i)_{i=0, \dots, l-1} \in \R^l_1$,
we define the deformed preprojective algebra $\T_\lambda$ as
\[
 \T_\lambda = M_l(D(\Rep(Q, \delta)))^{GL(\delta)}
 \Bigm/ \sum_{i=0}^{l-1} M_l(D(\Rep(Q, \delta)))^{GL(\delta)}
 (\tau(e^{(i)}) - \lambda_i),
\]
and define the spherical subalgebra $\A_\lambda$ of $\T_\lambda$
as $\A_\lambda = e_0 \T_\lambda e_0$ where $e_i = E_{ii}$ for
$i=0$, $\dots$, $l-1$. It is clear that
\begin{equation}
\label{eq:50} 
 \A_\lambda = D(\Rep(Q, \delta))^{GL(\delta)}
 \Bigm/ \sum_{i=0}^{l-1} D(\Rep(Q, \delta))^{GL(\delta)}
 (\iota(e^{(i)}) -  \bar{\lambda}_i)
\end{equation}
where $\bar{\lambda}_0 = \lambda_0 - 1$ and $\bar{\lambda}_i = \lambda_i$
for $i \ne 0$.

By the proposition of \cite{Ho} together with \cite{CBH},
we have the following isomorphisms of algebras.

\begin{proposition}[\cite{CBH}, \cite{Ho}]
 \label{prop:6}
For $\lambda = (\lambda_i)_{i=0, \dots, l-1}$, 
$\lambda_i = \kappa_{i+1} - \kappa_i + 1/l$, we have
isomorphisms of algebras:
\begin{equation}
\label{eq:49} 
 H_\kappa \simeq \T_\lambda, \qquad U_\kappa \simeq \A_\lambda.
\end{equation}
This isomorphisms are given by
\begin{align*}
\T_\lambda &\longrightarrow H_\kappa, \\
A_i = E_{i, i-1} \otimes t_i &\mapsto \bar{e}_{i} x \bar{e}_{i-1}, \\
A^*_i = E_{i-1, i} \otimes
 \partial_i &\mapsto \bar{e}_{i-1} y \bar{e}_{i}, \\
 e_i &\mapsto \bar{e}_i.
\end{align*}
\end{proposition}
Set $A = A_0 + A_1 + \dots + A_{l-1}$, and $A^*= A^*_0 +
A^*_1 + \dots + A^*_{l-1}$. They correspond to $x$,
$y \in H_\kappa$ under the above isomorphism. We have
the following triangular decomposition of $\T_\lambda$.
\begin{equation}
\label{eq:23}
\T_\lambda = \C[A]  \otimes_{\C} \left(\bigoplus_{i=0}^{l-1} \C e_i \right)
\otimes_{\C} \C[A^*].
\end{equation}

By the isomorphism \refeq{eq:49}, we identify the rational Cherednik 
algebra $H_\kappa$ and the deformed preprojective algebra 
$\T_\lambda$ with $\lambda_i = \kappa_{i+1} - \kappa_i + (1/l)$.
Thus we regard category $\catO_\kappa$ 
as a subcategory of $\fmod{\T_\lambda}$, and denote
it by $\catO_\lambda$. Then the category $\catO_\lambda$ is
the subcategory of $\fmod{\T_\lambda}$ such that the operator
$A^*$ acts locally nilpotently on each object of $\catO_\lambda$.
We also regard the standard modules of $H_\kappa$ as 
$\T_\lambda$-modules. Denote them by $\Delta_\lambda(i)$ for $i=0$, 
$\dots$, $l-1$. As $\T_\lambda$-modules, we have the natural 
description of the standard modules
\begin{equation}
 \label{eq:24}
 \Delta_\lambda(i) = (\T_\lambda / \T_\lambda A^*) e_i.
\end{equation}
By \refeq{eq:23}, we have 
\[
 \Delta_\lambda(i) = \C[A] \hwvec_i
\]
as a vector space.
By \refprop{prop:9}, $\Delta_\lambda(i)$ has a unique simple 
quotient which we denote by $L_\lambda(i)$.

\begin{lemma}
 \label{lemma:20}
 We have $\Hom_{\T_\lambda}(\Delta_\lambda(j), \Delta_\lambda(i)) \ne 0$
 if and only if $\lambda_{i} + \dots + \lambda_{j-1} \in \Z_{\leq 0}$.
 Moreover in this case,
 $\Hom_{\T_\lambda}(\Delta_\lambda(j), \Delta_\lambda(i))$ is 
 one-dimensional and any non-zero homomorphism from
 $\Delta_\lambda(j)$ to $\Delta_\lambda(i)$ is
 injective. 
 \begin{proof}
  To construct a homomorphism from $\Delta_\lambda(j)$ to 
  $\Delta_\lambda(i)$,
  it is enough to find a vector
  $v \in \Delta_\lambda(i)$ such 
  that $A^* v = 0$ and 
  $e_k v = \delta_{kj} v$ for $k=0$, $\dots$, $l-1$. 
  Indeed, if we have a
  non-zero homomorphism $\phi \in \Hom(\Delta_\lambda(j), \Delta_\lambda(i))$, 
  the vector 
  $v = \phi(\hwvec_j)$ satisfies $A^* v = 0$ and $e_k v = \delta_{kj} v$.
  Conversely, assume there is a vector $v$ 
  such 
  that $A^* v = 0$ and 
  $e_k v = \delta_{kj} v$ for $k=0$, $\dots$, $l-1$.
  Then we can define a homomorphism $\phi : \Delta_\lambda(j) 
  \rightarrow \Delta_\lambda(i)$ as 
  $\phi(A^m \hwvec_j) = A^m v$ for $m \in \Z_{\geq 0}$.
  Moreover, $v$ is equal to $A^p \hwvec_i$ and $i+p$ is equivalent to $j$
  modulo $l$. 
  
  Assume the above $v = A^p \hwvec_i \in \Delta_\lambda(i)$ exists.
  By the relation $[A^*, A] = \sum_{k=0}^{l-1} \lambda_k e_k$, we have
\begin{equation}
\label{eq:51} 
   0 = A^* v = A^* (A^p \hwvec_i) = [A^*, A^p] \hwvec_i
   =  l \sum_{k=0}^{p-1} \lambda_{i+k} \hwvec_i.
\end{equation}
  By $\lambda_0 + \dots + \lambda_{l-1} = 1$ and $i+p \equiv j$ modulo
  $l$, we have 
\[
\sum_{k=0}^{p-1} \lambda_{i+k} =
  n (\lambda_i + \dots + \lambda_{i-1}) + 
  (\lambda_i + \dots + \lambda_{j-1}) = n + (\lambda_i + \dots + \lambda_{j-1}).
\]
where $n = (p-j+i)/l \in \Z_{\geq 0}$. 
Then, we have
$\lambda_i + \dots + \lambda_{j-1} = -n \in \Z_{\leq 0}$.
Conversely, when $-n = \lambda_i + \dots + \lambda_{j-1} \in \Z_{\leq 0}$,
the vector $v = A^{nl + j - i} \hwvec_{i}$ satisfies 
$A^* v = 0$ and $e_k v = \delta_{jk} v$.
  Moreover, since such $v$ is uniquely determined by $i$, $j$ and $n$, we have
  $\Hom_{\T_\lambda}(\Delta_\lambda(j), \Delta_\lambda(i)) = \C$.
  Obviously we have $\T_\lambda v = \C[A] v$, thus this map is injective.
 \end{proof}
\end{lemma}

For $i \ne j$ such that 
$\Hom_{\T_\lambda}(\Delta(j), \Delta(i)) \ne 0$. Let
$L_\lambda(i,j)$ be the quotient 
\[
 0 \rightarrow \Delta_\lambda(j) {\rightarrow}
 \Delta_\lambda(i) \rightarrow L_\lambda(i,j) \rightarrow 0.
\]
By the above lemma $L_\lambda(i,j)$ is uniquely determined.
By the proof of \reflemma{lemma:20}, we have
$L_\lambda(i,j) \simeq \C[A] \hwvec_i / \C[A] A^{nl+j-i} \hwvec_i$ 
for some $n \in \Z_{\geq 0}$. Therefore, we have
$\langle \underline{\dim} L_\lambda(i,j), \lambda \rangle = 
\sum_{k=0}^{nl+j-i-1} \lambda_{i+k} = 0$ by \refeq{eq:51}.

Consider the following functors between the categories of modules:
\begin{align*}
 E_\lambda : \modu{\T_\lambda} &\longrightarrow \modu{\A_\lambda}, \\
 M &\mapsto e_0 M, \\
 F_\lambda : \modu{\A_\lambda} &\longrightarrow \modu{\T_\lambda}, \\
 N &\mapsto \T_\lambda e_0 \otimes_{\A_\lambda} N.
\end{align*}
Restricting the functors $E_\lambda$ and $F_\lambda$, we 
have functors between $\fmod{\T_\lambda}$ and 
$\fmod{\A_\lambda}$. We also denote it by the same symbols $E_\lambda$
and $F_\lambda$.

\begin{proposition}
 \label{prop:11}
 If $\langle \lambda, \beta \rangle \ne 0$ for all Dynkin
 roots $\beta \in \Z^l$, $E_\lambda$ is an equivalence of
 categories with quasi-inverse $F_\lambda$.
\begin{proof}
 The following proof is essentially the same as the argument in
 the proof of Theorem 3.3 of \cite{GS1}.

 To prove the equivalence, we show that 
 $\T_\lambda e_0 \otimes_{\A_\lambda} e_o \T_\lambda \simeq
 \T_\lambda e_0 \T_\lambda = \T_\lambda$ and
 $e_0 \T_\lambda \otimes_{\T_\lambda} \T_\lambda e_0 \simeq \A_\lambda$.
 It is clear that $e_0 \T_\lambda \otimes_{\T_\lambda} \T_\lambda e_0
 \simeq e_0 \T_\lambda e_0 = \A_\lambda$.
 Assume that $\T_\lambda e_0 \T_\lambda \ne \T_\lambda$, so then
 $\T_\lambda e_0 \T_\lambda$ is proper two-sided ideal.
 By the generalized Duflo theorem proved by \cite{Gi},
 $\T_\lambda e_0 \T_\lambda$ annihilates the irreducible module
 $L_\lambda(i)$ for some $i=0$, $\dots$, $l-1$. Namely, there
 is an $i=0$, $\dots$, $l-1$ such that $e_0 L_\lambda(i) = 0$.
 If $\Delta_\lambda(i) = L_\lambda(i)$, then 
 $e_0 L_\lambda(i) = e_0 \Delta_\lambda(i) \ne 0$ because
 we have $\Delta_\lambda(i) = \C[x] \hwvec_i$. Thus, it
 contradicts the assumption $e_0 L_\lambda(i) = 0$. 
 Assume $\Delta_\lambda(i) \ne L_\lambda(i)$, then there is
 an exact sequence 
\[
 0 \rightarrow \Delta_\lambda(j) \rightarrow
 \Delta_\lambda(i) \rightarrow L_\lambda(i) \rightarrow 0
\]
for some $j \ne i$. Let $\alpha = \underline{\Dim} L_\lambda(i)$,
then we have $\langle \lambda, \alpha \rangle = 0$ and
$\alpha \in \Z^l$ is a root. Moreover, 
by the assumption $e_0 L_\lambda = 0$, $\alpha$ is a Dynkin
root. This is a contradiction.
\end{proof}
\end{proposition}

\subsection{Parameters and orderings}
\label{sec:parameters-orderings}

In the next subsection, we  define
a functor $\calS_\lambda^\theta : \modu{\A_\lambda} \rightarrow
\modu{\A_{\lambda+\theta}}$ called the shift functor.
The shift functor $\calS_\lambda^\theta$ depends on the parameter
$\lambda \in \R^l_1$ and $\theta \in \Z^l_0$. 
 In this paper, we concentrate 
our attention on the case where $\calS_\lambda^\theta$ is an
 equivalences of categories. 
In this section we define our space of parameters.

For $\lambda = (\lambda_i)_{i=0, \dots, l-1} \in \R^l_1$, we 
have the highest weight category $\catO_\lambda$. We have the
ordering $\repgeq{\lambda}$ on the index set 
$\Lambda = \{0, \dots, l-1\}$ which
arises from the 
structure of the highest weight category $\catO_\lambda$.
Namely,
\begin{equation}
\label{eq:3}
 i \repgeq{\lambda} j \Leftrightarrow
 \Hom_{\T_\lambda}(\Delta_\lambda(j), \Delta_\lambda(i)) \ne 0 
\Leftrightarrow 
 \lambda_{i} + \dots + \lambda_{j-1} \in \Z_{\leq 0}
\end{equation}
as we proved in \reflemma{lemma:20}. 

Define 
\begin{equation}
\label{eq:7} 
 \R^l_{reg} = \{\lambda = (\lambda_i)_{i=0, \dots, l-1}
\in \R^l_1\; | \; \bar{\lambda}_{i} + \dots + \bar{\lambda}_{j-1} \ne 0
\quad \text{for all $i \ne j$}\}.
\end{equation}
where $\bar{\lambda}_i = \lambda_i - \delta_{i0}$. 

Fix $\lambda = (\lambda_i)_{i=0, \dots, l-1} \in \R^l_{reg}$.
Then we have $\bar{\lambda}_{i} + \dots + \bar{\lambda}_{j-1} < 0$
for all $i \repgo{\lambda} j$. 

The set of parameters $\Z^l_{reg}$ defined in \refeq{eq:16} is
separated into $(l-1)!$ alcoves by the hyperplanes 
$\theta_i + \dots + \theta_{j-1} = 0$ for $i \ne j$.
Set 
\begin{equation}
\label{eq:6} 
 \Z^l_{\lambda} = \{ \theta = (\theta_i)_{i=0, \dots, l-1}
\in \Z^l_{reg}\; | \;\theta_{i} + \dots + \theta_{j-1} < 0
\quad \text{if $\lambda_i + \dots + \lambda_{j-1} \in \Z_{\leq 0}$}\}.
\end{equation}
The set $\Z_\lambda^l$ is a union of alcoves in $\Z^l_{reg}$
depending on $\lambda$. 
If $\lambda \in \R^l_{reg}$ is generic, we have
 $\Z^l_\lambda = \Z^l_{reg}$. If $\lambda$ belongs to 
$\R^l_{reg} \cap \Z^l$,
$\Z^l_\lambda$ is one of $(l-1)!$ alcoves in $\Z^l_{reg}$.
For
$\lambda \in \R^l_{reg}$, $\theta \in \Z^l_\lambda$ and
$m \in \Z_{\geq 0}$, we have
$\Z^l_{\lambda + m \theta} = \Z^l_{\lambda}$ and
$\repgeq{\lambda + m \theta}$ is equal to $\repgeq{\lambda}$.

By the result of \cite{Ro}, we have the following theorem and
corollary. 
\begin{theorem}[\cite{Ro}]
 \label{thm:4}
 If $\catO_\lambda$ and $\catO_{\lambda'}$ have the same
 ordering defined by \refeq{eq:3}, then there exists
 an equivalence of categories
 $\catO_\lambda \simeq \catO_{\lambda'}$.
\end{theorem}

\begin{corollary}
 \label{cor:6}
 For $\lambda \in \R^l_{reg}$ and $\theta \in \Z^l_\lambda$, 
 there is an equivalence of categories between 
 $\catO_\lambda$ and $\catO_{\lambda+\theta}$.
\end{corollary}

Fix $\lambda = (\lambda_i)_{i=0, \dots, l-1} \in \R^l_{reg}$
and $\theta = (\theta_i)_{i=0, \dots, l-1} \in \Z^l_\lambda$.
Define a new ordering $\unrhd_\theta$ of $\Lambda$ by
\begin{equation}
 \label{eq:4}
  i \unrhd_\theta j
  \Leftrightarrow \theta_{i} + \dots + \theta_{j-1} \leq 0.
\end{equation}
It is a total ordering and it refines the ordering $\repgeq{\lambda}$,
i.e., $i \repgo{\lambda} j$ implies $i \unrhd_{\theta} j$.
If $\theta$ and $\theta'$ belong to the same alcove, $\unrhd_\theta$
is equal to $\unrhd_{\theta'}$. If we take $\theta$ and 
$\theta'$ from different alcoves, $\unrhd_\theta$ is different from
$\unrhd_{\theta'}$. 

By \refeq{eq:35}, the ordering $\unrhd_{\theta}$ is exactly
same as the ordering $\geogeq{\theta}$ defined in \refsec{sec:defin-quiv-vari}.
Therefore we have 
\begin{equation}
\label{eq:5} 
 \ii_{l} \rhd_\theta \ii_{l-1} \rhd_\theta \dots
 \rhd_\theta \ii_1.
\end{equation}

The following lemma will be repeatedly used in the next 
subsection. 
\begin{lemma}
 \label{lemma:7}
 For any $i > j$, we have
 $\bar{\lambda}_{\ii_j} + \dots + \bar{\lambda}_{\ii_i-1} > 0$
 or
 $\bar{\lambda}_{\ii_j} + \dots + \bar{\lambda}_{\ii_i-1} \not\in \Z$.
\begin{proof}
 Assume $\bar{\lambda}_{\ii_j} + \dots + \bar{\lambda}_{\ii_i-1} \in
 \Z$. By \refeq{eq:3}, we have $\ii_i \repgo{\lambda} \ii_j$ or
 $\ii_j \repgo{\lambda} \ii_i$. Since $\rhd_\theta$ refines 
 $\repgo{\lambda}$ and we have $\ii_i \rhd_\theta \ii_j$, 
 the case $\ii_j \repgo{\lambda} \ii_i$ cannot occur. 
 Thus we have $\ii_{i} \repgo{\lambda} \ii_j$.
 Therefore we have $\bar{\lambda}_{\ii_j} + \dots + \bar{\lambda}_{\ii_i-1}
 >0$ by \refeq{eq:7}.
\end{proof}
\end{lemma}

\subsection{Shift functors}
\label{sec:shift-functors}

As in \cite{Bo}, we define a functor called the shift functor between the two
categories of modules of the 
rational Cherednik algebras with different parameters.
Moreover we prove that it gives the equivalence of categories that we
 discussed in
\refthm{thm:4} and \refcor{cor:6}. 

Fix a parameter $\lambda = (\lambda_i)_{i=0, \dots, l-1} \in \R^l_1$ of 
the rational Cherednik algebra. We take another parameter 
$\theta = (\theta_i)_{i=0, \dots, l-1} \in \Z^l_0$. 
Define
\begin{equation}
\label{eq:52} 
 \B_\lambda^\theta = \left[D(\Rep(Q, \delta))
 \Bigm/ \sum_{i=0}^{l-1} D(\Rep(Q, \delta))
 (\iota(e^{(i)}) -  \bar{\lambda}_i)\right]^{GL(\delta), \chi_\theta}
\end{equation}
where $\chi_\theta$ is the character of $GL(\delta)$ defined in
\refsec{sec:defin-quiv-vari}. It is easy to see that
$\B_\lambda^\theta$ has an $(\A_{\lambda+\theta}, \A_\lambda)$-bimodule
structure. 

\begin{definition}
 We define the functor
\begin{align*}
 \calS_\lambda^{\theta}: \modu{\A_\lambda} &\longrightarrow
 \modu{\A_{\lambda+\theta}}, \\
 M &\mapsto \B_\lambda^\theta \otimes_{\A_\lambda} M.
\end{align*}
The functor $\calS_\lambda^{\theta}$ is called a shift functor.
Restricting $\calS_\lambda^\theta$ to the subcategory 
$\fmod{\A_\lambda}$, we have the functor
\[
\calS_\lambda^\theta : \fmod{\A_\lambda} \longrightarrow
 \fmod{\A_{\lambda+\theta}}.
\]
\end{definition}

In \refsec{sec:deform-prepr-algebr}, we defined the functors
$E_\lambda$, $F_\lambda$ between $\modu{\T_\lambda}$ and
$\modu{\A_\lambda}$. Using $E_\lambda : \modu{\T_\lambda} \rightarrow
\modu{\A_\lambda}$ and $F_{\lambda+\theta} : \modu{\A_{\lambda+\theta}}
\rightarrow \modu{\T_{\lambda+\theta}}$,
we define the shift functor 
\[
 \widehat{\calS}_\lambda^\theta = F_{\lambda+\theta} \circ 
 \calS_\lambda^\theta \circ E_\lambda :
 \modu{\T_\lambda} \longrightarrow \modu{\T_{\lambda+\theta}}.
\]
We also denote the restricted functor
\[
 \widehat{\calS}_\lambda^\theta :
 \fmod{\T_\lambda} \longrightarrow \fmod{\T_{\lambda+\theta}}
\]
by the same symbol $\widehat{\calS}_\lambda^\theta$.
\begin{lemma}
 \label{lemma:8}
 The functor 
 $\widehat{\calS}_\lambda^\theta$
 restricts to a functor 
 \[
  \widehat{\calS}_\lambda^\theta : \catO_\lambda \longrightarrow
 \catO_{\lambda+\theta}.
 \]
\begin{proof}
 Fix $M \in \catO_\lambda$. Since we have 
 $(A^*)^l = (e_0 + \dots + e_{l-1}) \otimes \partial_0 \cdots \partial_{l-1}$,
to prove the lemma, we only need to
 show $\partial_0 \cdots \partial_{l-1}$ acts locally nilpotently
on $\calS_\lambda^\theta(e_0 M) = \B_\lambda^\theta 
\otimes_{\A_\lambda} e_0 M$.
Fix $b \in \B_\lambda^\theta$ and $m \in e_0 M$.
We have $(\partial_0 \cdots \partial_{l-1})^p m = 0$ for a sufficiently 
large $p \in \Z_{\geq 0}$. 
Consider a filtration $\{F_k \B_\Lambda^\theta\}_{k}$ defined by
the degree
\[
\deg t_i = 1, \quad \deg \partial_i = 0 \qquad (i=0, \dots, l-1).
\]
For $b \in F_k \B_\lambda^\theta \backslash F_{k-1}\B_\lambda^\theta$, 
let $q$ be an integer greater
than $k+p$. Then, we have 
\[
 [(\partial_0 \cdots \partial_{l-1})^q, b] = b' (\partial_0 \cdots \partial_{l-1})^p
\]
Here $b' = \sum_{j=1}^{k} \mathrm{ad}(\partial_0 
\cdots \partial_{l-1})^j (b') \cdot (\partial_0 
\cdots \partial_{l-1})^{q-p-j} \in \B_\lambda^\theta$. Thus, we have
\begin{align*}
(\partial_0 \cdots \partial_{l-1})^q b \otimes m &=
 b (\partial_0 \cdots \partial_{l-1})^q \otimes m +
 [(\partial_0 \cdots \partial_{l-1})^q, b] \otimes m \\
 &= b \otimes (\partial_0 \cdots \partial_{l-1})^q m 
+ b' \otimes (\partial_0 \cdots \partial_{l-1})^p m = 0.
\end{align*}
Therefore $\partial_0 \cdots \partial_{l-1}$ acts locally nilpotently
on $\calS_\lambda^\theta(e_0 M)$.
\end{proof} 
\end{lemma}

Let $\bar{\A}_\lambda$ be the subalgebra of $\A_\lambda$
generated by the elements $t_0 \partial_0$ and $t_0 \cdots t_{l-1}$.

\begin{lemma}
 \label{lemma:9}
 For $k=1$, $\dots$, $l-1$,
 let $b^\theta_k$ be the non-negative integer defined by
 \refeq{eq:20} with $\theta' = \theta$.
 For $k=1$, $\dots$, $l-1$, $n = 0$, $\dots$, $b_k^\theta - 1$, define
 \[
 \tilde{g}_{k}(n) = (t_{\ii_{k+1}} \cdots t_{\ii_l})^{b^\theta_{k}-n}
 \prod_{j=k+1}^{l-1} (t_{\ii_{j+1}} \cdots t_{\ii_l})^{b_j^\theta}
 (\partial_{\ii_1} \cdots \partial_{\ii_k})^{n}
 \prod_{j=1}^{k-1} (\partial_{\ii_1} \cdots \partial_{\ii_j})^{b^\theta_j}
\]
 and, for $n \in \Z_{\geq 0}$, define
 \[
 \tilde{g}_l(n)
 = (\partial_0 \cdots \partial_{l-1})^n
 \prod_{j=1}^{l-1} (\partial_{\ii_1} \cdots \partial_{\ii_j})^{b^\theta_j}.
 \]
Then $\{\tilde{g}_k(n)\}_{k,n}$ generates $\B_\lambda^\theta$ as a
left $\bar{\A}_\lambda$-module.
\begin{proof}
 Consider the filtration of $\B_\lambda^\theta$ defined by 
 the order of differential operators in $D(\Rep(Q, \delta))$. 
By \cite[(5.1)]{Bo}, the
 associated graded module is
 \[
  \gr \B_\lambda^\theta = \C[\mu^{-1}(0)]^{GL(\delta), \chi_\theta}.
 \]
 Thus the statement of the lemma follows from \refprop{prop:3} and
 \reflemma{lemma:2}. 
\end{proof}
\end{lemma}

\begin{proposition}
 \label{prop:7}
 For $\lambda = (\lambda_i)_{i=0, \dots, l-1} \in \R^l_{reg}$
 and $\theta = (\theta_i)_{i=0, \dots, l-1} \in \Z^l_\lambda$,
 we have 
 \[
  \calS_\lambda^\theta(e_0 \Delta_\lambda(i))
 \simeq e_0 \Delta_{\lambda+\theta}(i) \quad \text{and} \quad
 \widehat{\calS}_\lambda^\theta(\Delta_\lambda(i)) 
 \simeq \Delta_{\lambda+\theta}(i)
 \]
 for all $i=0$, $\dots$, $l-1$.
\begin{proof}
 We show that $\calS_\lambda^\theta(e_0 \Delta_\lambda(\ii_i))$ is
 isomorphic to $e_0 \Delta_{\lambda+\theta}(\ii_i)$ for all 
 $i=1$, $\dots$, $l$. To prove this, we see the structure of 
 $\calS_\lambda^\theta(e_0 \Delta_\lambda(\ii_i))$ with the help of
 the geometric information which we studied in
 \refsec{sec:line-bundles}. As a result of it, 
 we can construct the isomorphism
 $\calS_\lambda^\theta(e_0 \Delta_\lambda(\ii_i)) \simeq
 e_0 \Delta_{\lambda+\theta}(\ii_i)$ explicitly. 

 Let $w_k(n) = \tilde{g}_k(n) \otimes e_0 t_0 t_{l-1} \dots t_{\ii_i+1}
 \hwvec_{\ii_i}$ be an element of
 $\calS_\lambda^\theta(e_0 \Delta_\lambda(\ii_i)) =
\B_\lambda^\theta \otimes_{\A_\lambda} e_0 \Delta_\lambda(\ii_i)$.
By \reflemma{lemma:9}, $\{w_k(n)\}_{k,n}$ span the module
$\B_\lambda^\theta \otimes_{\A_\lambda} e_0 \Delta_\lambda(\ii_i)$.
We show that the vector $w_i(0)$ generates 
$\B_\lambda^\theta \otimes_{\A_\lambda} e_0 \Delta_\lambda(\ii_i)$
and $(\partial_0 \cdots \partial_{l-1}) w_k(0) = 0$.

 First we show $w_k(n)$ is non-zero when $k < i$ or $k=i$ and $n=0$. 
 We identify $w_k(0) = w_{k-1}(b_{k-1}-1)$. Then, by a straightforward 
 calculation, we have
\begin{align*}
 \lefteqn{(t_0 \cdots t_{l-1}) w_k(n)} &\\
 &=
 (t_0 \cdots t_{l-1})
(t_{\ii_{k+1}} \cdots t_{\ii_l})^{b^\theta_{k}-n}
 \prod_{j=k+1}^{l-1} (t_{\ii_{j+1}} \cdots t_{\ii_l})^{b_j^\theta} \\
 & \qquad
 (\partial_{\ii_1} \cdots \partial_{\ii_k})^{n}
 \prod_{j=1}^{k-1} (\partial_{\ii_1} \cdots
 \partial_{\ii_j})^{b^\theta_j}
 \otimes e_0 t_0 t_{l-1} \cdots t_{\ii_i} \hwvec_{\ii_i} \\
 &= 
(t_{\ii_{k+1}} \cdots t_{\ii_l})^{b^\theta_{k}-n+1}
 \prod_{j=k+1}^{l-1} (t_{\ii_{j+1}} \cdots t_{\ii_l})^{b_j^\theta}
 (t_{\ii_1} \partial_{\ii_1} \cdots t_{\ii_k} \partial_{\ii_k}) \\
 & \qquad
 (\partial_{\ii_1} \cdots \partial_{\ii_k})^{n-1}
 \prod_{j=1}^{k-1} (\partial_{\ii_1} \cdots
 \partial_{\ii_j})^{b^\theta_j}
 \otimes e_0 t_0 t_{l-1} \cdots t_{\ii_i} \hwvec_{\ii_i}.
\end{align*}
For $1 \leq p \leq k$, we have
\begin{align*}
\lefteqn{ (t_{\ii_p} \partial_{\ii_p}) 
 (\partial_{\ii_1} \cdots \partial_{\ii_k})^{n-1}
 \prod_{j=1}^{k-1} (\partial_{\ii_1} \cdots
 \partial_{\ii_j})^{b^\theta_j}
 \otimes e_0 t_0 t_{l-1} \cdots t_{\ii_i} \hwvec_{\ii_i}} & \\
 &= - 
\Bigl(\bar{\lambda}_{\ii_p} + \dots +
 \bar{\lambda}_{\ii_i-1} + \sum_{q=p}^{k-1} b_q^\theta + n - 1\Bigr) \\
 & \qquad
  (\partial_{\ii_1} \cdots \partial_{\ii_k})^{n-1}
  \prod_{j=1}^{k-1} (\partial_{\ii_1} \cdots
 \partial_{\ii_j})^{b^\theta_j}
 \otimes e_0 t_0 t_{l-1} \cdots t_{\ii_i} \hwvec_{\ii_i}.
\end{align*}
Thus we have 
\[
 (t_0 \cdots t_{l-1}) w_k(n)
 \left\{\prod_{p=1}^{k} - \Bigl(\bar{\lambda}_{\ii_p} + \dots +
 \bar{\lambda}_{\ii_i-1} + \sum_{q=p}^{k-1} b_q^\theta + n - 1\Bigr)
\right\}
 w_k(n-1).
\]
 For any
 $p=1$, $\dots$, $i-1$ we have $\bar{\lambda}_{\ii_p} + \dots
 + \bar{\lambda}_{\ii_i-1} > 0$ or 
 $\bar{\lambda}_{\ii_p} + \dots
 + \bar{\lambda}_{\ii_i-1} \not\in \Z$ by \reflemma{lemma:7}. 
 Thus the coefficient of the right hand side of this equation
 is non-zero. 

 Therefore we have
 \[
  (t_0 \cdots t_{l-1})^{\sum_{j=1}^{i-1} b_j^\theta} w_i(0)
= C w_1(0) =
 C \prod_{j=1}^{l-1} (t_{\ii_{j+1}} \cdots t_{\ii_l})^{b_j^\theta}
 \otimes e_0 t_0 t_{l-1} \cdots t_{\ii_i+1} \hwvec_{\ii_i}
 \]
 where $C \in \C \backslash \{0\}$. Since 
 the right hand side of this equation is non-zero, so is
 $w_i(0)$. Moreover, for $k < i$, $w_k(n)$ is non-zero and it 
 belongs to $\C[t_0 \cdots t_{l-1}] w_i(0)$. Then, we have
\begin{equation}
\label{eq:27} 
\C[t_0 \cdots t_{l-1}] w_i(0) = \C w_i(0) \oplus 
 \bigoplus_{k < i, n} \C w_k(n),
\end{equation}
and it is $\C[t_0 \cdots t_{l-1}]$-free.

 Next we show $w_k(n) = 0$ when $k > i$ or $k=i$  and $n \geq 1$.
 Inserting $\prod_{j=k+1}^l t_{\ii_j} \partial_{\ii_j}$ into 
 the factors of $w_k(n)$, we have
\begin{align}
 \lefteqn{(t_{\ii_{k+1}} \cdots t_{\ii_l})^{b_k^\theta - n}
 \prod_{j=k+1}^{l-1} (t_{\ii_{j+1}} \cdots t_{\ii_l})^{b^\theta_j}
 \Bigl(\prod_{j=k+1}^{l} t_{\ii_j} \partial_{\ii_j}\Bigr)} & 
\nonumber\\
 \lefteqn{\qquad \times(\partial_{\ii_1} \cdots \partial_{\ii_k})^{n}
\prod_{j=1}^{k-1} (\partial_{\ii_1} \cdots \partial_{\ii_j})^{b^\theta_j}
 \otimes e_0 t_0 t_{l-1} \cdots t_{\ii_i+1} \hwvec_{\ii_i}}  & 
 \nonumber\\
 &=  
(t_{\ii_{k+1}} \cdots t_{\ii_l})^{b_k^\theta - n + 1}
 \prod_{j=k+1}^{l-1} (t_{\ii_{j+1}} \cdots t_{\ii_l})^{b^\theta_j}
 (\partial_0 \cdots \partial_{l-1})
 \nonumber\\
 &\qquad \times
 (\partial_{\ii_1} \cdots \partial_{\ii_k})^{n - 1}
\prod_{j=1}^{k-1} (\partial_{\ii_1} \cdots \partial_{\ii_j})^{b^\theta_j}
 \otimes e_0 t_0 t_{l-1} \cdots t_{\ii_i+1} \hwvec_{\ii_i} \nonumber\\
 \label{eq:26}
 &= 0. 
\end{align}
On the other hand, we have
\begin{align}
 \lefteqn{(t_{\ii_{k+1}} \cdots t_{\ii_l})^{b_k^\theta - n}
 \prod_{j=k+1}^{l-1} (t_{\ii_{j+1}} \cdots t_{\ii_l})^{b^\theta_j}
 \Bigl(\prod_{j=k+1}^{l} t_{\ii_j} \partial_{\ii_j}\Bigr)} & \nonumber\\
 \lefteqn{\qquad \times(\partial_{\ii_1} \cdots \partial_{\ii_k})^{n}
\prod_{j=1}^{k-1} (\partial_{\ii_1} \cdots \partial_{\ii_j})^{b^\theta_j}
 \otimes e_0 t_0 t_{l-1} \cdots t_{\ii_i+1} \hwvec_{\ii_i} } &
 \nonumber\\
 \label{eq:25}
&= \left\{ \prod_{j=k+1}^{l} - 
(\bar{\lambda}_{\ii_j} + \dots + \bar{\lambda}_{\ii_i-1}) \right\}
 w_k(n) 
\end{align}
 For $j \ne i$, we have
 $\bar{\lambda}_{\ii_j} + \dots + \bar{\lambda}_{\ii_i-1} \ne 0$.
 By \refeq{eq:26} and \refeq{eq:25}, we have
\begin{equation}
\label{eq:28}
 w_k(n) = 0
\end{equation}
 for $k > i$ or $k=i$ and $n \geq 1$. 

By \reflemma{lemma:9}, we have
 \[
  \B_\lambda^\theta \otimes_{\A_\lambda} e_0 \Delta_\lambda(\ii_i)
= \sum_{k, n} \C[t_0 \cdots t_{l-1}] w_k(n).
 \]
 Then, by \refeq{eq:28}, we have
\begin{equation}
\label{eq:30} 
  \B_\lambda^\theta \otimes_{\A_\lambda} e_0 \Delta_\lambda(\ii_i)
 = \C[t_0 \cdots t_{l-1}] w_i(0) + 
 \sum_{k < i, n} \C[t_0 \cdots t_{l-1}] w_k(n).
\end{equation}
 By \refeq{eq:25}, we have $w_k(n) \in \C[t_0 \cdots t_{l-1}] w_i(0)$
 for $k < i$. Thus we have,
\begin{equation}
\label{eq:29} 
 \sum_{k < i, n} \C[t_0 \cdots t_{l-1}] w_k(n) =
 \C[t_0 \cdots t_{l-1}] w_i(0).
\end{equation}
By \refeq{eq:30} and \refeq{eq:29}, we have
\[
  \B_\lambda^\theta \otimes_{\A_\lambda} e_0 \Delta_\lambda(\ii_i)
 = \C[t_0 \cdots t_{l-1}] w_i(0).
\]

Therefore we have $A^* w_i(0) = 0$ and $w_i(0)$ generates the module
$\B_\lambda^\theta \otimes_{\A_\lambda} e_0 \Delta_\lambda(\ii_i)$.
Thus we have $\calS_{\lambda}^\theta(e_0 \Delta_\lambda(\ii_i)) \simeq
e_0 \Delta_{\lambda+\theta}(\ii_i)$.

By \refprop{prop:11}, $F_{\lambda+\theta}$ is an equivalence of
categories. Thus, we have \break
$F_{\lambda+\theta}(e_0 \Delta_{\lambda+\theta}(i)) \simeq 
\Delta_{\lambda+\theta}(i)$. Therefore we have
\[
 \widehat{\calS}_\lambda^\theta(\Delta_\lambda(i))
= F_{\lambda+\theta} \circ \calS_\lambda^\theta \circ
E_\lambda (\Delta_\lambda(i))
\simeq F_{\lambda+\theta}(e_0 \Delta_{\lambda+\theta}(i)) \simeq 
\Delta_{\lambda+\theta}(i).
\]
\end{proof}
\end{proposition}

Next, we show that the shift
functor $\calS_\lambda^\theta$ is an equivalence of categories
between $\modu{\A_\lambda}$ and $\modu{\A_{\lambda+\theta}}$.

\begin{lemma}
 \label{lemma:11}
 For $\lambda \in \R^l_{reg}$, $\theta \in \Z^l_\lambda$ and
 $i$, $j=0$, $\dots$, $l-1$ such that $i \repgo{\lambda} j$,
 the shift functor $\widehat{\calS}_\lambda^\theta$ sends the exact 
 sequence in $\catO_\lambda$,
\[
 0 \rightarrow \Delta_\lambda(j) \overset{\varphi}{\rightarrow}
 \Delta_\lambda(i) \rightarrow L_\lambda(i,j) \rightarrow 0
 \] 
to the exact sequence in $\catO_{\lambda+\theta}$,
\[
 0 \rightarrow \Delta_{\lambda+\theta}(j) \rightarrow 
 \Delta_{\lambda+\theta}(i) \rightarrow L_{\lambda+\theta}(i,j) \rightarrow 0.
\] 
\begin{proof}
 By \refprop{prop:7}, we have 
 $\widehat{\calS}_\lambda^\theta(\Delta_\lambda(k)) = 
 \Delta_{\lambda+\theta}(k)$ for $k=i$, $j$. Then 
 $\widehat{\calS}_\lambda^\theta(\varphi)$ is a homomorphism
 \[
  \widehat{\calS}_\lambda^\theta(\varphi) : 
 \Delta_{\lambda+\theta}(j) \longrightarrow
 \Delta_{\lambda+\theta}(i).
 \]
By \reflemma{lemma:20}, $\widehat{\calS}_\lambda^\theta(\varphi)$ is injective
and its quotient is $L_{\lambda+\theta}(i,j)$. Since 
 $\widehat{\calS}_\lambda^\theta$ is a right exact functor,
 it implies $\widehat{\calS}_\lambda^\theta(L_\lambda(i,j)) \simeq
 L_{\lambda+\theta}(i,j)$.
\end{proof}
\end{lemma}

\begin{corollary}
 \label{cor:4}
 For $\lambda \in \R^l_{reg}$, $\theta \in \Z^l_\lambda$
 and $i=0$, $\dots$, $l-1$, we have
 \[
  \widehat{\calS}_\lambda^\theta(L_\lambda(i)) \simeq L_{\lambda+\theta}(i).
 \]
\end{corollary}

\begin{proposition}
 \label{prop:8}
 For $\lambda \in \R^l_{reg}$ and $\theta \in \Z^l_\lambda$,
 the functor $\widehat{\calS}_\lambda^\theta$ is an exact functor 
 from $\catO_\lambda$ to $\catO_{\lambda+\theta}$.
\begin{proof}
 Since $\widehat{\calS}_\lambda^\theta$ is right exact, to prove the
 exactness it is enough to show that 
 $\widehat{\calS}_\lambda^\theta$ sends 
 injective homomorphisms to injective homomorphisms. 
 Assume there is 
 a non-zero module $M \in \catO_\lambda$ such that 
 $\widehat{\calS}_\lambda^\theta(M) = 0$.
 Without loss of generalities, we can suppose that $M$ is 
 irreducible. By \refprop{prop:9}, 
 $M$ is isomorphic to $L_\lambda(i)$ for some 
 $i=0$, $\dots$, $l-1$. On the other hand, we have
 $\widehat{\calS}_\lambda^\theta(L_\lambda(i)) \simeq L_{\lambda+\theta}(i)$
 by \refcor{cor:4}. 
 This contradicts the assumption $\widehat{\calS}_\lambda^\theta(M) = 0$.
\end{proof}
\end{proposition}

The following proposition is a result of the general theory
of highest weight categories. 
R.~Rouquier suggested it to the author as an approach to proving
that $\widehat{\calS}_\lambda^\theta$ is an equivalence. 
using it to prove the equivalence of $\calS_\lambda^\theta$.
The following proof of the proposition is
given by S.~Ariki.

\begin{proposition}
 \label{prop:10}
 Assume there are two highest weight categories 
 $(\catO, \Lambda)$, $(\catO', \Lambda')$ which are equivalent
to each other. If an exact functor $F: \catO \longrightarrow \catO'$
 preserves the partial orderings of $\Lambda$ and $\Lambda'$, and
 $F$ sends the standard modules of $\catO$ to the standard modules
 of $\catO'$, then $F$ is an equivalence of categories.
\begin{proof}
We denote the partial ordering of index poset $\Lambda$ by
$\rhd$. 
We also denote the standard modules of $\catO$ by $\Delta(i)$ and
 the simple modules of $\catO$ by $L(i)$ for $i \in \Lambda$. For
$i \in \Lambda$, let $P(i)$ be the projective cover of $L(i)$.

Let $G : \catO' \longrightarrow \catO$ be an equivalence of
categories. Consider the exact functor 
$F' = G \circ F : \catO \longrightarrow \catO$. Since $G$ and
$F$ preserve the partial orderings of $\Lambda$ and $\Lambda'$,
so does $F'$. Since $G$ is the equivalence, $F'$ is an
equivalence of categories if and only if $F$ is an equivalence of
categories. Therefore we assume 
$\catO' = \catO$ and
$F : \catO \longrightarrow \catO$ is an exact functor such that
$F(\Delta(i)) \simeq \Delta(i)$ for any $i \in \Lambda$.

First we show that $F(L(i)) \simeq L(i)$ for any $i \in \Lambda$
by induction on $i$. If $i$ is minimal in $\Lambda$, we 
have $L(i) = \Delta(i)$. Thus we have $F(L(i)) \simeq L(i)$.
Assume $F(L(j)) \simeq L(j)$ for all $j \lhd i$. Consider
the exact sequence 
\[
 0 \rightarrow N(i) \rightarrow \Delta(i) \rightarrow
 L(i) \rightarrow 0.
\]
In the composition factors of $N(i)$, only $L(j)$ with $j \lhd i$ appears.
Thus $F(N(i))$ and $N(i)$ has the same composition factors 
by the hypothesis of the induction. Therefore we have
$F(L(i)) \simeq L(i)$. 

Second, we have
\begin{equation}
\label{eq:46}
\Ext^n(M, N) \simeq \Ext^n(F(M), F(N))
\end{equation}
 for any $M$, $N$ and $n$ by inductions on the length of
$M$ and $N$. In particular, $F$ is fully faithful.

By \refeq{eq:46}, we have $\Ext^1(F(P(i)), L(j)) = 0$ for
any $i$, $j \in \Lambda$. Moreover, we have
\[
 \Ext^1(F(P(i)), M) = 0
\]
for any $i \in \Lambda$ and $M$ by induction on the length of $M$.
Thus $F(P(i))$ is a projective object in $\catO$. 
Since, $\End(F(P(i))) \simeq \End(P(i))$ is a local ring, 
$F(P(i))$ is indecomposable. On the other hand, $F(P(i))$ has
$F(L(i)) \simeq L(i)$ as its quotient. Therefore, we have
$F(P(i)) \simeq P(i)$ for all $i \in \Lambda$. 

Let $A$ be a finite dimensional algebra such that $\catO \simeq \fmod{A}$.
Since $F(P(i)) \simeq P(i)$ for all $i \in \Lambda$, we have
$F(A) \simeq A$. 

Therefore we have
\[
 F(M) \simeq F(A \otimes_A M) \simeq F(A) \otimes_A M
 \simeq M.
\]
for any $M \in \catO$. Therefore $F$ is an equivalence of
categories.
\end{proof}
\end{proposition}

\begin{remark}
 Since we proved 
 $\widehat{\calS}_\lambda^\theta(L_\lambda(i)) \simeq
 L_{\lambda+\theta}(i)$ in \refcor{cor:4}, we actually do not need
 the first part of the above proof.
\end{remark}

\begin{theorem}
 \label{thm:5}
 For $\lambda \in \R^l_{reg}$ and $\theta \in \Z^l_\lambda$, 
 the shift functor $\widehat{\calS}_\lambda^\theta: \catO_{\lambda} 
 \longrightarrow \catO_{\lambda+\theta}$ is an 
 equivalence of categories.
\begin{proof}
 By \refcor{cor:6}, we have an equivalence 
 $\catO_\lambda \simeq \catO_{\lambda+\theta}$. By \refprop{prop:8}
 $\widehat{\calS}_\lambda^\theta$ is an exact functor. The assumption
 of \refprop{prop:10} is satisfied for $F = \widehat{\calS}_\lambda^\theta$ by
 \refprop{prop:7}. Then $\widehat{\calS}_\lambda^\theta$ is an 
 equivalence of categories.
\end{proof}
\end{theorem}

\begin{corollary}
 \label{cor:5}
 The shift functor $\calS_\lambda^\theta$ is an equivalence of categories
 between $\modu{\A_{\lambda}}$ and $\modu{\A_{\lambda+\theta}}$,
 and $\widehat{\calS}_\lambda^\theta$ is
 an equivalence of categories between $\modu{\T_\lambda}$ and
 $\modu{\T_{\lambda+\theta}}$.
\begin{proof}
 This proof is essentially same as the proof of Theorem~3.3 in
 \cite{GS1}.

 To prove the equivalence, we show that 
 $\B_{\lambda+\theta}^{-\theta} \otimes_{\A_{\lambda+\theta}}
 \B_{\lambda}^{\theta} \simeq \A_{\lambda}$ and
 $\B_\lambda^\theta \otimes_{\A_{\lambda}}
 \B_{\lambda+\theta}^{-\theta} \simeq \A_{\lambda+\theta}$. 
 Assume that $I := \B_{\lambda+\theta}^{-\theta} \otimes_{\A_{\lambda+\theta}}
 \B_{\lambda}^{\theta} \simeq \B_{\lambda+\theta}^{-\theta} 
 \B_{\lambda}^{\theta} \ne \A_{\lambda}$. Then $I$ is a 
 proper two-sided ideal of $\A_\lambda$. By the generalized Duflo theorem
 proved in \cite{Gi}, $I$ annihilates a irreducible module
 $e_0 L_\lambda(i)$ for some $i=0$, $\dots$, $l-1$.
 However, by \refthm{thm:5}, we have 
 \[
  I e_0 L_\lambda(i) \simeq 
 \B_{\lambda+\theta}^{-\theta} \otimes_{\A_{\lambda+\theta}}
 \B_{\lambda}^{\theta} \otimes_{\A_\lambda} e_0 L_\lambda(i)
 \simeq e_0 L_\lambda(i).
 \]
 Therefore $I = \A_\lambda$. The second isomorphism can be
 proved similarly. 
\end{proof}
\end{corollary}

In the rest of this subsection, we consider the
$(\A_{\lambda+\theta}, \T_\lambda)$-bimodule
$\calS_\lambda^\theta(e_0 \T_\lambda) =
\B_\lambda^\theta \otimes_{\A_\lambda} e_0 \T_\lambda$.

First we define a space
\[
 \calM_\lambda = M_l(D(\Rep(Q, \delta))) \Bigm/
 \sum_{i=0}^{l-1} M_l(D(\Rep(Q, \delta))) (\tau(e^{(i)}) - \lambda_i).
\]
Then, consider the following natural maps:
\begin{align*}
 \phi_1 : e_0 \T_\lambda = e_0 \calM_\lambda^{GL(\delta)} 
 &\longrightarrow e_0 \calM_\lambda, \\
 \phi_2 : e_0 \A_\lambda = e_0 \T_\lambda e_0 
 &\longrightarrow e_0 \calM_\lambda, \\
 \phi_3 : \B_\lambda^\theta 
 = e_0 \calM_\lambda^{GL(\delta), \chi_\theta} e_0
 &\longrightarrow e_0 \calM_\lambda. 
\end{align*}
The above maps $\phi_1$, $\phi_2$ and $\phi_3$ are
clearly injective.
Then, these injective maps induce a 
map
\[
 \B_\lambda^\theta \otimes_{\A_\lambda} e_0 \T_\lambda
 \longrightarrow \B_\lambda^\theta e_0 \T_\lambda \subset 
 e_0 \calM_\lambda.
\]
Clearly the image of this map is inside the
subspace $e_0 \calM_\lambda^{GL(\delta), \chi_\theta}$,
thus we have the following map:
\begin{align}
 \label{eq:54}
 \Theta : \B_\lambda^\theta \otimes_{\A_\lambda} e_0 \T_\lambda
 &\longrightarrow e_0 \calM_\lambda^{GL(\delta), \chi_\theta}, \\
 b \otimes m &\mapsto \phi_3(b) \phi_1(m). 
\nonumber
\end{align}
It is a homomorphism of $(\A_{\lambda+\theta}, \T_\lambda)$-bimodules.

\begin{lemma}[\cite{Bo}, Lemma~6.8]
 \label{lemma:12}
 Let $B$ be a left Ore domain and $P$ an 
 $(A, B)$-bimodule which yields Morita equivalence between
 $A$ and $B$. If $P'$ is torsion free $A$-module, then
 every surjective homomorphism $P \longrightarrow P'$ is
 isomorphism.
\end{lemma}

Set
\begin{equation}
\label{eq:53} 
 \widetilde{\R}^l_{reg} = \{(\lambda_i)_{i=0, \dots, l-1} \in \R^l_{reg}
\; | \; \lambda_i + \dots + \lambda_{j-1} \ne 0 \quad
\text{for all $i \ne j$}\}.
\end{equation}

\begin{lemma}
 \label{lemma:13}
 For $\lambda \in \widetilde{\R}^l_{reg}$, $\theta \in \Z^l_\lambda$,
 the homomorphism $\Theta$ is injective.
\begin{proof}
 Set $\A_{\lambda}^{(i)} = e_i \T_\lambda e_i$ for each 
$i=0$, $\dots$, $l-1$. Each $\A_\lambda^{(i)}$ is a left
Ore domain. The $(\T_\lambda, \A_\lambda^{(i)})$-bimodule
$\T_\lambda e_i$ yields a Morita equivalence when 
$\lambda \in \widetilde{\R}^l_{reg}$
by the same argument as in \refprop{prop:11} for 
the $(\A_{\lambda}, \T_\lambda)$-bimodule $e_0 \T_\lambda$.
Therefore $\B_\lambda^\theta \otimes_{\A_\lambda} e_0 \T_\lambda e_i
= \B_\lambda^\theta \otimes_{\A_\lambda} e_0 \T_\lambda
 \otimes_{\T_\lambda} T_\lambda e_i$ yields a Morita equivalence 
between $\A_{\lambda+\theta}$ and $\A_\lambda^{(i)}$.
The module $\B_\lambda^\theta e_0 \T e_i$ is a torsion free 
 $\A_{\lambda+\theta}$-module. Applying \reflemma{lemma:12}
to $P = \B_\lambda^\theta \otimes_{\A_\lambda} e_0 \T_\lambda e_i$ and
$P' = \B_\lambda^\theta e_0 \T_\lambda e_i$,
 we have $\B_\lambda^\theta \otimes_{\A_\lambda} e_0
 \T_\lambda e_i \simeq \B_\lambda^\theta e_0 \T_\lambda e_i
\subseteq e_0 \calM_{\lambda}^{GL(\delta), \chi_\theta} e_i$
for $i=0$, $\dots$, $l-1$. Therefore
 we have $\B_\lambda^\theta \otimes_{\A_\lambda} e_0
 \T_\lambda \simeq \B_\lambda^\theta e_0 \T_\lambda 
\subseteq e_0 \calM_{\lambda}^{GL(\delta), \chi_\theta}$.
\end{proof}
\end{lemma}

\subsection{$q$-dimension of representations}
\label{sec:q-dimens-repr}

In this subsection we calculate the $q$-dimension of the module
$\C \otimes_{\C[t_0 \cdots t_{l-1}]} \B_\lambda^\theta
\otimes_{\A_\lambda} e_0 \T_\lambda$. This result is used in
\refsec{sec:proof-main-theorem} to prove our main theorem,
\refthm{thm:2}.

Consider the Euler operator
\[
 \eu_\lambda = \sum_{i=0}^{l-1} A_i A_i^* - \sum_{i=0}^{l-1}
 c_i(\lambda) e_i  \in \T_\lambda
\]
where $c_i(\lambda) \in \R$ such that
$c_{i+1}(\lambda) - c_i(\lambda) = l \kappa_{i+1} - l \kappa_i = l \lambda_i - 1$
and $c_0(\lambda) + c_1(\lambda) + \dots + c_{l-1}(\lambda) = 1$.

For $\hwvec_i \in \Delta_\lambda(i)$, $\eu_\lambda$ acts as follows:
\begin{equation}
 \label{eq:44}
  \eu_\lambda \hwvec_i = -c_i(\lambda) \hwvec_i.
\end{equation}

The following lemma is proved by a straightforward calculation.
\begin{lemma}
 \label{lemma:1}
 \begin{enumerate}
  \item $[\eu_\lambda, A_i] = A_i$.
  \item $[\eu_\lambda, A^*_i] = -A^*_i$.
  \item $[\eu_\lambda, e_i] = 0$.
 \end{enumerate}
\end{lemma}

A $\T_\lambda$-module (or $\A_\lambda$-module) $M$ is 
called a graded module if $M$ has a vector space 
decomposition $M = \bigoplus_m M_m$ such that 
$A_i M_m \subseteq M_{m+1}$, $A^*_i M_m \subseteq M_{m-1}$
and $e_i M_m \subseteq M_m$ for all $i$.

For a graded module $M = \bigoplus_m M_m$ and $k \in \R$,
let $M[k]$ be the graded module shifted degree by $k$, 
i.e.
\[
 (M[k])_m = M_{m+k}.
\]

For a module $M \in \catO_\lambda$, $M$ has a
vector space decomposition $M = \bigoplus_m M_m$ where
$M_m$ is the generalized eigenspace for an eigenvalue $m$ with respect
to $\eu_\lambda$. By \reflemma{lemma:1}, this decomposition
makes $M$ a graded module. This grading is called the canonical
grading. 

For a standard module $\Delta_\lambda(i)$ ($i=0$, $\dots$, $l-1$),
let $\widetilde{\Delta}_\lambda(i)$ be a graded module which is 
isomorphic to $\Delta_\lambda(i)$ as an ungraded module, and
\[
 \widetilde{\Delta}_\lambda(i) = \bigoplus_{m \in \Z_{\geq 0}}
 \left(\widetilde{\Delta}_\lambda(i)\right)_m, \qquad
  \left(\widetilde{\Delta}_\lambda(i)\right)_m = \C A^m \hwvec_i.
\]
Considering $\Delta_\lambda(i)$ to be a graded $\T_\lambda$-module
with the canonical grading, we have
\begin{equation}
\label{eq:31} 
 \Delta_\lambda(i) \simeq \widetilde{\Delta}_\lambda(i)[-c_i(\lambda)]
\end{equation}
as a graded $\T_\lambda$-module by \refeq{eq:44}.

For a $(\T_{\lambda'}, \T_\lambda)$-bimodule $M$, let
${}_{\lambda'}\eu_\lambda$ be the operator on $M$:
\[
{}_{\lambda'}\eu_\lambda (m) =  \eu_{\lambda'} \cdot m - m \cdot
\eu_\lambda .
\]
for $m \in M$. We have the decomposition of $M$,
\[
 M = \bigoplus_{n} M_n
\]
where $M_n$ is the generalized eigenspace for
an eigenvalue $n$ with respect to the operator 
${}_{\lambda'}\eu_{\lambda}$. By this decomposition,
$M$ is a graded module. This grading is called the adjoint
grading. 

The $(\T_\lambda, \T_\lambda)$-bimodule $\T_\lambda$ has
the decomposition
\[
\T_\lambda = \bigoplus_{n \in \Z} (\T_\lambda)_n
\]
where $(\T_\lambda)_n$ is the eigenspace for an eigenvalue
$n \in \Z$ with respect to the operator ${}_\lambda \eu_\lambda$.
This grading coincides with the grading given by the
degree,
\[
 \deg E_{ij} \otimes t_k = 1, \quad
 \deg E_{ij} \otimes \partial_k = -1.
\]

For an $(\A_{\lambda'}, \A_\lambda)$-bimodule $M$, let
${}_{\lambda'}\eu_\lambda$ be the operator on $M$,
\[
{}_{\lambda'}\eu_\lambda (m) = e_0 \eu_{\lambda'} e_0 \cdot m - m \cdot
e_0 \eu_\lambda e_0.
\]
Then, the operator ${}_{\lambda'}\eu_\lambda$ gives $M$ 
the structure of a graded module. We also call this grading
the adjoint grading. 

 Consider the following decomposition of $\B_\lambda^\theta$,
 \[
  \B_\lambda^\theta = \bigoplus_{n \in \Z}
 \left(\B_\lambda^\theta\right)_n
 \]
where $\left(\B_\lambda^\theta\right)_n$ is the
eigenspace 
for an eigenvalue $n$ with respect to the operator
${}_{\lambda+\theta}\eu_{\lambda}$. By the above decomposition,
$\B_\lambda^\theta$ is a graded module. 
This grading coincides with the grading given by the 
 degree,
\[
 \deg t_k = 1, \quad \deg \partial_k = -1.
\]

In the rest of this subsection, we assume
$\lambda \in \widetilde{\R}^l_{reg}$ and $\theta \in \Z^l_\lambda$.

\begin{lemma}
 For $i=0$, $\dots$, $l-1$, we consider the grading of the module
\[
 \calS_\lambda^\theta(\widetilde{\Delta}_\lambda(i)) =
\B_\lambda^\theta \otimes_{\A_\lambda} e_0 \widetilde{\Delta}_\lambda(i)
\]
 induced from the adjoint grading of $\B_\lambda^\theta$ and 
 the grading of $\widetilde{\Delta}_\lambda(i)$. Then we have
 \[
 \B_\lambda^\theta \otimes_{\A_\lambda} e_0 \widetilde{\Delta}_\lambda(i)
  \simeq e_0 \widetilde{\Delta}_{\lambda+\theta}(i)[d_i^\theta]
 \]
 where $d_i^{\theta}$ are the integers defined on \refeq{eq:21}.
\begin{proof}
 We have the two gradings on the standard modules $\Delta_\lambda(i)$
 and $\Delta_{\lambda+\theta}(i)$. Let $\deg_{can}$ be the degree
 defined by the canonical grading and let $\deg_\lambda$ be 
 the degree defined by the grading of $\widetilde{\Delta}_\lambda(i)$.
 For the other graded modules, let $\deg$ be the degree of
 the grading of each module.

 We have 
 \[
 \B_\lambda^\theta \otimes_{\A_\lambda} e_0 \Delta_\lambda(i)
  \simeq e_0 {\Delta}_{\lambda+\theta}(i)
 \]
 as ungraded modules by \refprop{prop:7}. Considering the
 canonical grading, we have 
 \[
  \deg_{can}( b \otimes v ) = \deg b + \deg_{can} v
 \]
 for $b \in \B_\lambda^\theta$ and 
 $v \in e_0 \widetilde{\Delta}_\lambda(i)$. 

 By \refeq{eq:31}, we have
 \begin{align*}
  \deg_{can}(v) &= \deg_\lambda(v) - c_i(\lambda) \\
  \deg_{can}(v') &= \deg_{\lambda+\theta}(v') - c_i(\lambda+\theta)
 \end{align*}
 for $v \in \Delta_\lambda(i)$ and $v' \in \Delta_{\lambda+\theta}(i)$.
 Therefore we have
\begin{align*}
 \deg_{\lambda+\theta}(b \otimes v) &= 
 \deg_{can}(b \otimes v) + c_i(\lambda+\theta) \\
 &= \deg(b) + \deg_{can}(v) + c_i(\lambda+\theta) \\
 &= \deg(b) + \deg_{\lambda}(v) + c_i(\lambda+\theta) -
 c_i(\lambda) \\ 
  &= \deg(b \otimes v) + c_i(\lambda+\theta) -  c_i(\lambda).
\end{align*}
 On the other hand, we can easily obtain
 \begin{gather*}
  d_{i+1}^\theta - d_i^\theta = l \theta_i = 
  (c_{i+1}(\lambda+\theta) - c_i(\lambda+\theta)) - 
  (c_{i+1}(\lambda) - c_i(\lambda)), \\
  d_0^\theta + d_1^\theta + \dots + d_{l-1}^\theta = 0.
 \end{gather*}
 Thus we obtain
 \[
  c_i(\lambda+\theta) - c_i(\lambda) = d_i^\theta.
 \]
\end{proof}
\end{lemma}

For a graded module $M = \bigoplus_m M_m$, define the
$q$-dimension of $M$ as a formal series
\[
 \Dim_q M = \sum_{m} (\Dim M_m) q^m.
\]
By the above lemma, we have the $q$-dimension of 
$\calS_\lambda^\theta(e_0 \widetilde{\Delta}_\lambda(i))
= \B_\lambda^\theta \otimes_{\A_\lambda} e_0 \widetilde{\Delta}_\lambda(i)$,
\[
 \dim_q \B_\lambda^\theta \otimes_{\A_\lambda} e_0
 \widetilde{\Delta}_\lambda(i) = 
  q^{d_i^\theta + l-i} \frac{1}{1-q^l}
\]
for $i=0$, $\dots$, $l-1$. Here we set $d_l^\theta = d_0^\theta$.

The adjoint gradings of $\B_\lambda^\theta$ and $\T_\lambda$ induce
the gradings of the $(\A_{\lambda+\theta}, \T_\lambda)$-bimodule
$\B_\lambda^\theta \otimes_{\A_\lambda} \T_\lambda$. They also induce
the grading of the left $\A_{\lambda+\theta}$-module
$\B_\lambda^\theta \otimes_{\A_\lambda} \T_\lambda \otimes_{\C[A^*]}
\C$
and the grading of the right $\T_\lambda$-module
$\C \otimes_{\C[t_0 \cdots t_{l-1}]}
\B_\lambda^\theta \otimes_{\A_\lambda} \T_\lambda$.

By \refeq{eq:24}, 
we have the following natural
isomorphism as graded $\A_{\lambda+\theta}$-modules
\[
 \B_\lambda^\theta \otimes_{\A_\lambda} e_0 \T_\lambda
 \otimes_{\C[A^*]} \C
 \simeq 
 \bigoplus_{i=0}^{l-1} \B_\lambda^\theta \otimes_{\A_\lambda} 
 e_0 \widetilde{\Delta}_\lambda(i)
\]
By the above equations, we have
\[
 \dim_q  \B_\lambda^\theta \otimes_{\A_\lambda} e_0 \T_\lambda
 \otimes_{\C[A^*]} \C
= q^{d_0^\theta} \frac{1}{1-q^l} + 
\sum_{i=1}^{l-1}   q^{d_i^\theta + l-i} \frac{1}{1-q^l}.
\]

\begin{lemma}[\cite{GS1}, Theorem A.1]
 \label{lemma:6}
 Let $R$ be a connected $\Z_{\geq 0}$-graded $\C$-algebra.
 Let $P$ be an $R$-module that is both graded and
 projective. 
Then $P$ is a graded-free $R$-module in the sense that P has a 
free basis of 
homogeneous elements. 
\end{lemma}

\begin{lemma}
\label{lemma:19}
 \begin{enumerate}
  \item The module $\calS_\lambda^\theta(e_0 \T_\lambda) = 
	\B_\lambda^\theta \otimes_{\A_\lambda}
e_0 \T_\lambda$ is graded-free as a left $\C[t_0 \cdots t_{l-1}]$-module and
graded-free as a right $\C[A^*]$-module.
  \item $\B_\lambda^\theta \otimes_{\A_\lambda}
e_0 \T_\lambda \otimes_{\C[A^*]} \C$ is a finitely generated,
graded-free $\C[t_0 \cdots t_{l-1}]$-module.
  \item $\C \otimes_{\C[t_0 \cdots t_{l-1}]}
	\B_\lambda^\theta \otimes_{\A_\lambda}
	e_0 \T_\lambda$ is a finitely generated,
	graded-free right $\C[A^*]$-module
 \end{enumerate}
 \begin{proof}
  The following proof is essential the same as the proof
  of \cite[Lemma 6.11]{GS1}.

  (1) By \refprop{prop:11} and \refcor{cor:5}, 
  $\B_\lambda^\theta \otimes_{\A_\lambda} e_0 \T_\lambda$ is
  projective as a left $\A_{\lambda+\theta}$-module and
  a right $\T_\lambda$-module. By the structure of the graded
  module $\B_\lambda^\theta \otimes_{\A_\lambda} e_0 \T_\lambda$
  which is defined above,
  $\B_\lambda^\theta \otimes_{\A_\lambda} e_0 \T_\lambda$ is
  graded as a left $\C[t_0 \cdots t_{l-1}]$-module and
  a right $\C[A^*]$-module. By \reflemma{lemma:6}, 
  $\B_\lambda^\theta \otimes_{\A_\lambda} e_0 \T_\lambda$ is
  graded-free as a $\C[t_0 \cdots t_{l-1}]$-module and
  a right $\C[A^*]$-module.

  (2) By \refeq{eq:24} and \refprop{prop:7}, we have
  \[
  \B_\lambda^\theta \otimes_{\A_\lambda} e_0 \T_\lambda 
  \otimes_{\C[A^*]} \C
  \simeq \bigoplus_{i=0}^{l-1} e_0 \Delta_{\lambda+\theta}(i)
  \simeq \C[t_0 \cdots t_{l-1}] \otimes \C \Z_l.
  \]
  Therefore, $\B_\lambda^\theta \otimes_{\A_\lambda} e_0 \T_\lambda 
  \otimes_{\C[A^*]} \C$ is graded-free.

  (3) First, we show that
  $\B_\lambda^\theta \otimes_{\A_\lambda} e_0 \T_\lambda$
  is a finitely generated right module over 
  $R = (\C[t_0 \cdots t_{l-1}])^{\mathrm{op}} \otimes_{\C} \C[A^*]$.
  By \reflemma{lemma:13}, 
  $\B_\lambda^\theta \otimes_{\A_\lambda} e_0 \T_\lambda \subset
  e_0 \calM_\lambda^{GL(\delta), \chi_\theta}$. Thus
  $\gr \B_\lambda^\theta \otimes_{\A_\lambda} e_0 \T_\lambda \subset
  \gr e_0 \calM_\lambda^{GL(\delta), \chi_\theta} =
  e_0 M_l(\C[\mu^{-1}(0)])^{GL(\delta), \chi_\theta}$, which is
  certainly a noetherian $\C[t_0 \cdots t_{l-1}] \otimes \C[A^*]$-module.
  The $\C[t_0 \cdots t_{l-1}] \otimes \C[A^*]$-module structure of
  $\gr \B_\lambda^\theta \otimes_{\A_\lambda} e_0 \T_\lambda$ is
  the one induced from the $R$-module structure of
  $\B_\lambda^\theta \otimes_{\A_\lambda} e_0 \T_\lambda$. 
  Therefore, $\B_\lambda^\theta \otimes_{\A_\lambda} e_0 \T_\lambda$
  is finitely generated. 

  By (2), $\Sigma = \{ A^*, (t_0 \cdots t_l-1) \}$ is a regular
  sequence for the right $R$-module 
  $\B_\lambda^\theta \otimes_{\A_\lambda} e_0 \T_\lambda$. 
  In particular, if $\mathfrak{n} = A^* R + (t_0 \cdots t_{l-1}) R$,
  then $\Sigma$ is a regular sequence for the $R_{\mathfrak{n}}$-module
  $(\B_\lambda^\theta \otimes_{\A_\lambda} e_0
  \T_\lambda)_{\mathfrak{n}}$.
  By the Auslander-Buchsbaum formula \cite[Ex.~4, p.114]{Ma}, 
  $(\B_\lambda^\theta \otimes_{\A_\lambda} e_0
  \T_\lambda)_{\mathfrak{n}}$ is free as a $R_{\mathfrak{n}}$-module.
  
  Finally, $\C \otimes_{\C[t_0 \cdots t_{l-1}]} \B_\lambda^\theta 
  \otimes_{\A_\lambda} e_0 \T_\lambda$ is a finitely generated,
  graded $\C[A^*]$-module and so corresponds to a $\C^*$-equivariant
  coherent sheaf on $\C$. Therefore, the locus where 
  $\C \otimes_{\C[t_0 \cdots t_{l-1}]} \B_\lambda^\theta 
  \otimes_{\A_\lambda} e_0 \T_\lambda$ is not free is a 
  $\C^*$-stable closed subvariety of $\C$. If this locus is
  non-empty, it must contain $0 \in \C$. By the conclusion of the 
  last paragraph, the stalk at $0 \in \C$ of 
  $\C \otimes_{\C[t_0 \cdots t_{l-1}]} 
  \B_\lambda^\theta \otimes_{\A_\lambda} e_0 \T_\lambda$
  is free. Therefore, $\C \otimes_{\C[t_0 \cdots t_{l-1}]} \B_\lambda^\theta 
  \otimes_{\A_\lambda} e_0 \T_\lambda$ must be free.
  Since $\C \otimes_{\C[t_0 \cdots t_{l-1}]} \B_\lambda^\theta 
  \otimes_{\A_\lambda} e_0 \T_\lambda$ is graded module, 
  it is graded free by \reflemma{lemma:6}.
 \end{proof}
\end{lemma}

By \reflemma{lemma:19}, a homogeneous $\C[A^*]$-basis of 
$\C \otimes_{\C[t_0 \cdots t_{l-1}]} \otimes B_\lambda^\theta
\otimes_{\A_\lambda} e_0 \T_\lambda$ is given by a homogeneous
$\C$-basis of 
$\C \otimes_{\C[t_0 \cdots t_{l-1}]} \otimes B_\lambda^\theta
\otimes_{\A_\lambda} e_0 \T_\lambda 
\otimes_{\C[A^*]} \C$. Therefore we have
the following proposition.

\begin{proposition}
 \label{prop:12}
 We have 
 \begin{equation}
 \Dim_q 
 \C \otimes_{\C[t_0 \cdots t_{l-1}]} \otimes B_\lambda^\theta
 \otimes_{\A_\lambda} e_0 \T_\lambda = 
 \sum_{i=1}^{l} q^{d_i^\theta + (l-i)} \frac{1}{1-q^{-1}}.
 \end{equation}
\end{proposition}

\section{Gordon-Stafford functors}
\label{sec:gord-staff-funct}
\subsection{$\Z$-algebras}

In this section, we define the functor $\widehat{\Phi}_\lambda^\theta$
of \refeq{eq:13} as Boyarchenko defined it in \cite{Bo}.
First we review the definition and basic
properties of $\Z$-algebras.

\begin{definition}
 A lower-triangular $\Z$-algebra $B$ is an algebra such that
\begin{enumerate}
 \item
 $B$ is bigraded by $\Z$ in the following way:
 $B = \bigoplus_{i \geq j \geq 0} B_{ij}$. 
 \item The multiplication
 of $B$ is defined in matrix fashion, i.e., $B$ satisfies
 $B_{ij} B_{jk} \subseteq B_{ik}$ for $i \geq j \geq k \geq 0$ 
 but $B_{ij} B_{lk} = 0$ if $j \ne l$.
 \item $B_{ii}$ is an unital subalgebra for all $i \in \Z_{\geq 0}$.
\end{enumerate}
\end{definition}

We also define graded modules of lower-triangular $\Z$-algebras.
Let $B$ be a lower-triangular $\Z$-algebra. A graded $B$-module
is $\Z_{\geq 0}$-graded left $B$-module $M = \bigoplus_{i \in \Z_{\geq
0}} M_i$, such that $B_{ij} M_j \subseteq M_i$ for all 
$i \geq j \geq 0$ and $B_{ij} M_k = 0$ if $j \ne k$.
Homomorphisms of graded $B$-modules are defined to be graded
homomorphisms of degree zero.

We denote the category of graded $B$-modules by 
$\Grmod{B}$, and denote the subcategory of finitely generated
graded $B$-modules by $\fgrmod{B}$.

A graded module $M = \bigoplus_{i \in \Z_{\geq 0}} M_i \in
\Grmod{B}$  is bounded if $M_i = 0$ all but finitely many
$i \in \Z_{\geq 0}$, and torsion if it is a direct limit
of bounded modules. We denote the subcategory of torsion
modules in $\Grmod{B}$ by $\Tor{B}$, and the subcategory
of bounded modules in $\fgrmod{B}$ by $\ftor{B}$. The 
corresponding quotient categories are written 
$\Qgr{B} = \Grmod{B} / \Tor{B}$ and 
$\fqgr{B} = \fgrmod{B} / \ftor{B}$.

For a graded commutative algebra $S = \bigoplus_{m \in \Z_{\geq 0}} S_m$,
we define a lower-triangular $\Z$-module $\widehat{S} = 
\bigoplus_{i \geq j \geq 0} \widehat{S}_{ij}$ where
$\widehat{S}_{ij} = S_{i - j}$ for $i \geq j \geq 0$.
Define the categories $\Grmod{S}$, $\dots$, $\fqgr{S}$ in the
usual manner. Then, as in Section 5.3 of \cite{GS1}, we have
equivalences of categories:
\begin{equation}
 \label{eq:14}
  \begin{split}
   \Qgr{S} &\longrightarrow \Qgr{\widehat{S}}, \\
   \fqgr{S} &\longrightarrow \fqgr{\widehat{S}}, \\
   M = \bigoplus_{i \in \Z_{\geq 0}} M_i &\mapsto 
   M = \bigoplus_{i \in \Z_{\geq 0}} M_i.
  \end{split}
\end{equation}

We define another example of $\Z$-algebra called a Morita
$\Z$-algebra.
Suppose we have countably many Morita equivalent algebras
$\{B_i\}_{i \in \Z_{\geq 0}}$ and $B_{ij}$ is a
$(B_i, B_j)$-bimodule which yields Morita equivalences for
$i > j \geq 0$. Moreover, suppose we have an isomorphism
$B_{ij} \otimes_{B_j} B_{jk} \simeq B_{ik}$ for $i \geq j \geq k \geq
0$. Set $B_{ii} = B_i$ and define the Morita $\Z$-algebra
$B$ to be $B = \bigoplus_{i \geq j \geq 0} B_{ij}$. Note that
our definition of Morita $\Z$-algebras is same as one of \cite{GS1},
and it requires a stronger condition than one of \cite{Bo}.

\begin{lemma}[GS1, Lemma 5.5]
 \label{lemma:17}
 Assume $B_0$ is noetherian, then
\begin{enumerate}
 \item Each finitely generated graded left $B$-module
       is noetherian.
 \item The association $\phi : 
       M \mapsto \bigoplus_{i \in \Z_{\geq 0}} B_{i0} \otimes_{B_0} M$
       induces an equivalence of categories between $\fmod{B_0}$
       and $\fqgr{B}$.
\end{enumerate}
\end{lemma}

\subsection{Construction of the functor}

For a stability parameter $\theta \in \Z^l_{reg}$, set
\[
S = \bigoplus_{m \in \Z_{\geq 0}} S_m, \qquad
S_m =
\C[\mu^{-1}(0)]^{GL(\delta), \chi_\theta^m} 
\]
as defined in
\refsec{sec:defin-quiv-vari}. By \refprop{prop:2}, we have
an isomorphism $\M_\theta(\delta) \simeq \Proj S$.

Let $\widehat{S} = \bigoplus_{i \geq j \geq 0} \widehat{S}_{ij}$
be the lower-triangular $\Z$-algebra obtained from the above 
graded algebra $S$. By \refeq{eq:14}, we have an 
equivalence of categories
$\fqgr{\widehat{S}} \simeq \Coh(\M_\theta(\delta))$.

Fix $\lambda \in \R^l_{reg}$ and $\theta \in \Z^l_\lambda$.
Recall the algebra $\A_{\lambda}$ of \refeq{eq:50} and
the bimodule $\B_\lambda^\theta$ of \refeq{eq:52}.
Set $B_i = \A_{\lambda + i\theta}$ for $i \in \Z_{\geq 0}$
and set $B_{ij} = \B_{\lambda + j\theta}^{(i-j) \theta}$ for
$i > j \geq 0$. By \refcor{cor:5}, $B_{ij}$ is a
$(B_i, B_j)$-bimodule which yields a Morita equivalence.

\begin{proposition}
 \label{prop:17}
 For $\lambda \in \R^l_{reg}$ and $\theta \in \Z^l_\lambda$, we
 have the isomorphism
\[
 \B_{\lambda + j\theta}^{(i-j) \theta} \otimes_{\A_{\lambda + j\theta}}
 \B_{\lambda + k \theta}^{(j-k) \theta} \simeq 
 \B_{\lambda + k\theta}^{(i-k) \theta}.
\]
\begin{proof}
 We apply \reflemma{lemma:12} for the algebras
 $A = \A_{\lambda+i\theta}$,
 $B = \A_{\lambda+k\theta}$ and the modules
$P =  \B_{\lambda + j\theta}^{(i-j) \theta} \otimes_{\A_{\lambda + j\theta}}
 \B_{\lambda + k \theta}^{(j-k) \theta}$,
 $P'= \B_{\lambda + k\theta}^{(i-k) \theta}$. It is clear that
 $B$ is a left Ore domain, $P'$ is torsion free and
 $P$ yields an equivalence of categories by \refcor{cor:5}.
 We have the
 surjective homomorphism
\begin{align*}
  P = \B_{\lambda + j\theta}^{(i-j) \theta} \otimes_{\A_{\lambda + j\theta}}
 \B_{\lambda + k \theta}^{(j-k) \theta}
 &\longrightarrow 
 P' = \B_{\lambda + k\theta}^{(i-k) \theta} \\
 b_1 \otimes b_2 &\mapsto b_1 b_2
\end{align*}
 Thus, by the above lemma, it is an isomorphism. 
\end{proof}
\end{proposition}
By the above proposition, we have the Morita $\Z$-algebra
$B = \bigoplus_{i \geq j \geq 0} B_{ij}$ where 
$B_{ii} = B_i = \A_{\lambda+i\theta}$. By \reflemma{lemma:17},
we have an equivalence of categories
\begin{align*}
 \fmod{\A_\lambda} &\longrightarrow \fqgr{B}, \\
 M &\mapsto \widetilde{M} = \bigoplus_{m \in \Z_{\geq 0}}
 \B_\lambda^{m \theta} \otimes_{\A_\lambda} M.
\end{align*}

The algebra $\A_\lambda$ and the bimodule $\B_\lambda^\theta$ 
are filtered by the order of differential operators in 
$D(\Rep(Q, \delta))$, and we have
 \begin{align*}
  \gr \A_\lambda &\simeq \C[\mu^{-1}(0)]^{GL(\delta)} = S_0, \\
  \gr \B_\lambda^{m \theta} &\simeq 
  \C[\mu^{-1}(0)]^{GL(\delta), \chi_\theta^m} = S_m.
 \end{align*}
The filtration induces the filtration on the Morita
$\Z$-algebra $B = \bigoplus_{i \geq j \geq 0} B_{ij}$.
Thus we have the following theorem as in \cite{GS1} and
\cite{Bo}. 

\begin{theorem}
 For $\lambda \in \R^l_{reg}$, $\theta \in \Z^l_\lambda$, define
 $\Z$-algebras $B$ and $\widehat{S}$ as above. Then:
\begin{enumerate}
 \item There is an equivalence of categories
       $\fmod{\A_\lambda} \simeq \fqgr{B}$.
 \item There is an isomorphism of lower-triangular $\Z$-algebras
       $\gr B \simeq \widehat{S}$.
 \item We have an equivalence of categories
       $\fqgr{\widehat{S}} \simeq \Coh(\M_\theta(\delta))$.
\end{enumerate}
\end{theorem}

Let $\Filt{\T_\lambda}$, $\Filt{\A_\lambda}$ be categories
of filtered modules. Given $(M, \Lambda) \in \Filt{\T_\lambda}$
(resp. $\in \Filt{\A_\lambda}$), we say $\Lambda$ is a good
filtration on $M$ if $\gr_{\Lambda} M$ is a finitely generated
$\gr \T_\lambda$-module (resp. $\gr \A_\lambda$-module). We
denote the subcategories of good filtered modules by 
$\filt{\T_\lambda}$ and $\filt{\A_\lambda}$. 

Let $(M, \Lambda) \in \Filt{\A_\lambda}$. Then each module
$M(i) = \B_\lambda^{i \theta} \otimes_{\A_\lambda} M$ is 
filtered by the tensor product filtration 
\[
 \Lambda^k M(i) = \sum_{l \in \Z} F^l \B_\lambda^{i\theta}
 \otimes \Lambda^{k-l} M
\]
where $F$ is the filtration of $\B_\lambda^{i \theta}$.
Therefore, the $B$-module $\widetilde{M} = \bigoplus_{i \in \Z_{\geq 0}}
M(i)$ is filtered and we have the graded 
$\widehat{S}$-module $\gr \widetilde{M} = \bigoplus_{i \in \Z_{\geq 0}}
\gr M(i)$ associated to $\widetilde{M}$. For a graded
$B$-module with filtration $(\widetilde{M}, \Lambda)$,
we call $\Lambda$ is a good filtration if 
$\gr_{\Lambda} \widetilde{M}$ is a finitely generated 
$\widehat{S}$-module. 
The following lemma is due to \cite{GS2}. 

\begin{lemma}[\cite{GS2}, Lemma 2.5]
 If $\Lambda$ is a good filtration of $M$, then the
 induced filtration $\Lambda$ on $\widetilde{M}$ is also good.
\end{lemma}

As in \cite{Bo}, we define the functor
\begin{align*}
 \Phi_\lambda^\theta : 
 \Filt{\A_\lambda} &\longrightarrow \Qcoh(\M_\theta(\delta)), \\
 (M, \Lambda) &\mapsto \gr_{\Lambda} \widetilde{M}.
\end{align*}
By restricting $\Phi_\lambda^\theta$ to $\filt{\A_\lambda}$, 
we have a functor
\[
  \Phi_\lambda^\theta : 
  \filt{\A_\lambda} \longrightarrow \Coh(\M_\theta(\delta)).
\]
We also define the functor from $\Filt{\T_\lambda}$ and
$\filt{\T_\lambda}$
\begin{align*}
 \widehat{\Phi}_\lambda^\theta = \Phi_\lambda^\theta
 \circ E_\lambda : \Filt{\T_\lambda} &\longrightarrow
 \Qcoh(\M_\theta(\delta)), \\
 \widehat{\Phi}_\lambda^\theta = \Phi_\lambda^\theta
 \circ E_\lambda : \filt{\T_\lambda} &\longrightarrow
 \Coh(\M_\theta(\delta)). 
\end{align*}
We call the above functors $\Phi_\lambda^\theta$ and 
$\widehat{\Phi}_\lambda^\theta$ the Gordon-Stafford functors.

\section{Construction of a tautological bundle}
\subsection{Main theorem}
\label{sec:main-theorem}
Now we consider our main theorem. We determine the image of 
the module $e_0 \T_\lambda$ by the functor $\Phi_\lambda^\theta$.

Recall the set of parameter $\widetilde{\R}^l_{reg}$ of \refeq{eq:53}.
Fix the parameters $\lambda = (\lambda_i)_{i=0, \dots, l-1} \in
\widetilde{\R}^l_{reg}$ and 
$\theta = (\theta_i)_{i=0, \dots, l-1} \in \Z^l_\lambda$.
Consider an $\A_\lambda$-module $e_0 \T_\lambda$. Considering 
the filtration by the order of differential operators in 
$D(\Rep(Q, \delta))$, $e_0 T_\lambda$ is a filtered 
$A_\lambda$-module with good filtration.
The following theorem is the main result of this paper.

\begin{theorem}
 \label{thm:2}
 For $\lambda \in \widetilde{\R}^l_{reg}$, $\theta \in \Z^l_\lambda$
 and $m \in \Z_{\geq 0}$,
we have the 
 isomorphism
 \[
 \gr \B_\lambda^{m \theta} \otimes_{\A_\lambda} e_0 \T_\lambda
 \simeq e_0 M_l(\C[\mu^{-1}(0)])^{GL(\delta), \chi_\theta^m}.
 \]
 Therefore we have
 \[
  \Phi_\lambda^{\theta}(e_0 \T_\lambda) \simeq \widetilde{\calP}_\theta
 \]
 as coherent sheaves on $\M_\theta(\delta)$.
\end{theorem}

In the next subsection, we will give the proof of
 the above theorem.

\subsection{Proof of the main theorem}
\label{sec:proof-main-theorem}

In this subsection, we complete the proof of our main 
theorem, \refthm{thm:2}. 

First, we construct the homomorphism of \refthm{thm:2}.
The module $\B_\lambda^\theta \otimes_{\A_\lambda} e_0 \T_\lambda$
is filtered by the tensor product filtration. The module 
$e_0 \calM_\lambda^{GL(\delta), \chi_\theta}$ is filtered by the order of
differential operators in $D(\Rep(Q, \delta))$. 
Clearly, the homomorphism $\Theta$ of \refeq{eq:54} 
is a homomorphism of filtered 
modules. Thus we have the associated homomorphism between
the associated graded modules
\[
 \gr \Theta : 
 \gr \left(\B_\lambda^\theta \otimes_{\A_\lambda} e_0 \T_\lambda \right)
 \longrightarrow \gr e_0 \calM_\lambda^{GL(\delta), \chi_\theta}
 = e_0 M_l(\C[\mu^{-1}(0)])^{GL(\delta), \chi_\theta}.
\]
This homomorphism was what we wanted to construct.

The $(\A_{\lambda+\theta}, \T_\lambda)$-bimodule 
$\B_\lambda^\theta \otimes_{\A_\lambda} e_0 \T_\lambda$ is graded
by the grading induced from the adjoint gradings of $\B_\lambda^\theta$
and $e_0 \T_\lambda$. The $(\A_{\lambda+\theta}, \T_\lambda)$-bimodule 
$e_0 \calM_\lambda^{GL(\delta), \chi_\theta}$ is graded by the
adjoint grading with respect to the operator 
${}_{\lambda+\theta} \eu_\lambda$. 

\begin{lemma}
 \label{lemma:14}
 The homomorphism $\Theta$ is
 homogeneous with respect to the gradings of 
$\B_\lambda^\theta \otimes_{\A_\lambda} e_0 \T_\lambda$ and
$e_0 \calM_\lambda^{GL(\delta), \chi_\theta}$. 
\end{lemma}

By \reflemma{lemma:13} we have
the injective homomorphism
\[
 \Theta : \B_\lambda^{m\theta} \otimes_{\A_\lambda} e_0 \T_\lambda
 \simeq \B_\lambda^{m\theta} e_0 \T_\lambda 
 \longrightarrow e_0 \calM_\lambda^{GL(\delta), \chi^m_\theta}.
\]
To complete the proof of theorem, we need to prove 
the equality of the inclusion
$\B_\lambda^\theta e_0 \T_\lambda \subseteq e_0 
\calM_\lambda^{GL(\delta), \chi_\theta^m}$. The following proof
is essentially the same as the proof \cite[Section 6.17]{GS1}.

\begin{lemma}
 \label{lemma:16}
 The modules $\B_\lambda^{m \theta} e_0 \T_\lambda$ and
 $e_0 \calM_{\lambda}^{GL(\delta), \chi_\theta^m}$ are free 
 as left
 $\C[t_0 \cdots t_{l-1}]$ modules. 
\end{lemma}

\begin{lemma}
 \label{lemma:15}
 We have an equality of localized spaces
\[
 (\B_\lambda^{m \theta} e_0 \T_\lambda)[(t_0 \cdots t_{l-1})^{-1}]
 = e_0 \calM_\lambda^{GL(\delta), \chi_\theta^m}[(t_0 \cdots t_{l-1})^{-1}].
\]
\begin{proof}
 It is clear that 
 \[
 (\B_\lambda^{m \theta} e_0 \T_\lambda)[(t_0 \cdots t_{l-1})^{-1}]
 \subseteq e_0 \calM_\lambda^{GL(\delta), \chi^m_\theta}[(t_0 \cdots t_{l-1})^{-1}].
 \]
 Fix an arbitrary element $\sum_{i=0}^{l-1}  f_i \otimes E_{0i} 
\in e_0 \calM_{\lambda}^{GL(\delta), \chi^m_\theta}$ where
$f_i \in \B_\lambda^{m \theta + \tau_i}$. We show by induction that,
for any $f \otimes E_{0i} 
\in F_n e_0 \calM_\lambda^{GL(\delta), \chi^m_\theta} e_i$,
there exists $p \in \Z_{\geq 0}$ such that
\begin{equation}
 \label{eq:45}
  (t_0 \cdots t_{l-1})^p f \otimes E_{0i} 
\in \B_\lambda^{m \theta} e_0 \T_\lambda.
\end{equation}

For $n=0$, we have 
\[
 F_0 e_0 \calM_\lambda^{GL(\delta), \chi^m_\theta} =
 e_0 M_l(\C[t_0, \dots, t_{l-1}])^{GL(\delta), \chi_\theta^m}.
\]
Thus we have
\[
 (t_0 \cdots t_{l-1}) f = (t_{i+1} \cdots t_{l-1} t_0 f) (t_1 \cdots t_{i}).
\]
We have $(t_{i+1} \cdots t_{l-1} t_0 f) \in \B_\lambda^{m \theta}$ and
$(t_1 \cdots t_i) \otimes E_{0i} \in e_0 \T_\lambda e_i$. Therefore we have
\refeq{eq:45} for $n=0$. 

Assume we have \refeq{eq:45} for $n < n_0$. For 
$f \otimes E_{0i} 
\in F_{n_0} e_0 \calM_\lambda^{GL(\delta), \chi_\theta^m} e_i$,
we have
\[
 (t_0 \cdots t_{l-1}) f \otimes E_{0i} = f (t_0 \cdots t_{l-1}) \otimes E_{0i}
 + [t_0 \cdots t_{l-1}, f] \otimes E_{0i}.
\]
The second term $[t_0 \cdots t_{l-1}, f] \otimes E_{0i}$ belongs to
$F_{n_0-1} e_0 \calM_\lambda^{GL(\delta), \chi_\theta^m} e_i$. By
the hypothesis of the induction, there exists $p \in \Z_{\geq 0}$
such that
\[
 (t_0 \cdots t_{l-1})^{p} [t_0 \cdots t_{l-1}, f]
 \otimes E_{0i}\in \B_\lambda^{m \theta} e_0 \T_\lambda.
\]
We have $(t_0 \cdots t_{l-1})^{p} (f t_{i+1} \cdots t_{l-1} t_0)
\in \B_\lambda^{m \theta}$ and
$(t_1 \cdots t_i) \otimes E_{0i} \in e_0 \T_\lambda e_i$. 
Therefore, we have
\begin{align*}
 \lefteqn{(t_0 \cdots t_{l-1})^{p+1} f \otimes E_{0i}} &\\
&= 
 (t_0 \cdots t_{l-1})^p f (t_{i+1} \cdots t_{l-1} t_0) (t_1 \cdots t_i)
 \otimes E_{0i} \\
&\qquad + (t_0 \cdots t_{l-1})^p [t_0 \cdots t_{l-1}, f] \otimes E_{0i}
 \in \B_\lambda^{m
 \theta} e_0 \T_\lambda.
\end{align*}
Therefore we have \refeq{eq:45} for $n=n_0$. 
\end{proof}
\end{lemma}

Let $\{a_{gp}\}_{g,p}$ be a $\C[t_0 \cdots t_{l-1}]$-free basis of
$\B_\lambda^{m \theta} e_0 \T_\lambda$ such that 
$a_{gp}$ is a homogeneous vector of degree $g$. Let 
$\{b_{gq}\}_{g,q}$ be a $\C[t_0 \cdots t_{l-1}]$-free basis of 
$e_0 \calM_\lambda^{GL(\delta), \chi_\theta^m}$ such that 
$b_{gq}$ is a homogeneous vector of degree $g$. 

By \refthm{thm:3} and \refprop{prop:12}, we have the equality
\begin{equation}
 \label{eq:8}
\begin{split}
 \dim_q \C \otimes_{\C[t_0 \cdots t_{l-1}]} 
 \B_\lambda^{m \theta} e_0 \T_\lambda &=
 \sum_{i=1}^{l} q^{d_{i+1}^{m \theta} + (l-i-1)} \frac{1}{1-q^{-1}} 
 \\
 &= \dim_{q} \C \otimes_{\C[t_0 \cdots t_{l-1}]}
 e_0 \calM_\lambda^{GL(\delta), \chi_\theta^m}.
\end{split}
\end{equation}
By \refeq{eq:8}, we have:
\renewcommand{\theenumi}{($\dagger$\arabic{enumi})}
\begin{enumerate}
 \item \label{item:1}
       For any $g \in \R$, $\{a_{gp}\}_p$ and $\{b_{gq}\}_q$ have
       finite cardinality.
 \item \label{item:2}
       There is $T \in \R$ such that there is no nonzero
       $a_{gp}$ and $b_{gq}$ when $g > T$.
 \item \label{item:3}
       For any $g \in \R$ $\#\{a_{gp}\}_{p}$ is equal to 
       $\#\{b_{gq}\}_q$. 
\end{enumerate}
\renewcommand{\theenumi}{(\arabic{enumi})}
We show that we can adjust the basis $\{b_{gq}\}_q$ to
be equal to the basis $\{a_{gp}\}_p$ by a downward induction
on $g$. By \ref{item:3}, we have 
$\{a_{gp}\}_p = \{b_{gq}\}_q = \emptyset$ for $g > T$.

Let $-\infty < G \leq T$, suppose that 
$\{b_{gq}\}_q = \{a_{gp}\}_p$ for all $g > G$ by induction.
Suppose that there exists an element $b_{Gq_0}$ which 
does not belong to $\{a_{gp}\}_p$. By \reflemma{lemma:15},
there exists an integer $n \in \Z_{\geq 0}$ such that
we have a homogeneous equation
\begin{equation}
 \label{eq:10}
  (t_0 \cdots t_{l-1})^n b_{Gq_0} = 
  \sum_{g < G, p} c_{gp} a_{gp} + \sum_{p} c_{Gp} a_{Gp} +
  \sum_{g > G, q} c'_{gq} b_{gq}
\end{equation}
where each $c_{gp}$, $c'_{gq} \in \C[t_0 \cdots t_{l-1}]$.
Note that we use the hypothesis of the induction 
$\{a_{gp}\}_p = \{b_{gq}\}_q$ for $g > G$. Since 
$\B_\lambda^{m \theta} e_0 \T_\lambda \subseteq e_0 
\calM_\lambda^{GL(\delta), \chi_\theta^m}$, we can write 
each $a_{gp} = \sum_{h, q} d_{gp}^{hq} b_{hq}$ for some
$d_{gp}^{hq} \in \C[t_0 \cdots t_{l-1}]$. Thus we obtain
a homogeneous equation
\begin{equation}
 \label{eq:11}
  (t_0 \cdots t_{l-1})^n b_{Gq_0} = 
  \sum_{g < G, p, h, q} c_{gp} d_{gp}^{hq} b_{hq} 
  + \sum_{p, h, q} c_{Gp} d_{Gp}^{hq} b_{hq} 
  \sum_{g > G, q} c'_{gq} b_{gq}
\end{equation}
The above equations \refeq{eq:10}, \refeq{eq:11} are homogeneous
of degree $G+ln$. By \refeq{eq:10}, $\deg c_{gp} \geq ln$ for each
$g$ and $p$. Thus the $b_{gq}$ in the first two terms of the right 
hand side of \refeq{eq:11} has degree $\leq G$. Since 
$\{b_{gq}\}_{g,q}$ are $\C[t_0 \cdots t_{l-1}]$-basis, the third 
term $\sum_{g > G, q} c'_{gq} b_{gq}$ is actually zero in
\refeq{eq:10}, \refeq{eq:11}. 

Now consider where $b_{Gq_0}$ appears on the right hand side of
\refeq{eq:11}. For $g < G$, \refeq{eq:10} implies $\deg c_{gp}
> ln$ for each $p$, thus there is no $b_{Gq_0}$ in the first 
term of the right hand side of \refeq{eq:11}. Thus $b_{Gq_0}$
appears only in the second term of the right hand side of
\refeq{eq:11}.
Since $\{b_{gq}\}_{g,q}$ is $\C[t_0 \cdots t_{l-1}]$-basis, there
is nonzero $c_{Gp} d_{Gp}^{Gq_0} b_{Gq_0}$ in the second term of
the right hand side of \refeq{eq:11}. In this case by \refeq{eq:10}
we have $\deg c_{Gp} = nl$. Hence $d_{Gp}^{Gq_0} \in \C \backslash \{0\}$,
and we have 
\[
 a_{Gp} = d_{Gp}^{Gq_0} b_{Gq_0} + \sum_{(h, q) \ne (G, q_0)}
 d_{Gp}^{hq} b_{hq}.
\]
Thus we can replace $b_{Gq_0}$ by $a_{Gp}$ in our basis of
$e_0 \calM_\lambda^{GL(\delta), \chi_\theta^m}$. By \ref{item:3}, 
$\{a_{Gp}\}_p$ and $\{b_{Gq}\}_q$ have the same cardinality. 
After a finite number of steps, we have $\{b_{Gq}\}_q \subseteq
\{a_{Gp}\}_p$ and hence $\{b_{Gq}\}_q = \{a_{Gp}\}_p$. This
completes the induction on $g$. Thus, we have
$\{b_{gq}\}_q = \{a_{gp}\}_p$ for all $g$. It implies the 
equality $\B_\lambda^{m\theta} e_0 \T_\lambda = e_0
\calM_\lambda^{GL(\delta), \chi_\theta^m}$.

\subsection{Characteristic cycles}
\label{sec:char-cycl}

In this section, we determine the characteristic cycle of
the standard module $\Delta_\lambda(i)$ in $\M_\theta(\delta)$
for $i=0$, $\dots$, $l-1$.

Fix the parameters $\lambda = (\lambda_i)_{i=0, \dots, l-1}
 \in \widetilde{\R}^l_{reg}$ and $\theta = (\theta_i)_{i=0, \dots, l-1}
\in \Z^l_\lambda$. For $M \in \filt{A_\lambda}$, we consider
$\widetilde{M} = \bigoplus_{m \in \Z_{\geq 0}} \B_\lambda^{m \theta}
\otimes_{\A_\lambda} M$ and the coherent sheaf 
$\Phi_\lambda^\theta(M) = (\gr \widetilde{M})\sptilde$. 
Let the associated radical ideal 
$N_{\widetilde{M}}$ to be the radical ideal $N_{\widetilde{M}}
= \sqrt{\ann_{S} \gr \widetilde{M}}$. The characteristic
variety of $M$ is
\[
 \Char(M) = \supp(\Phi_\lambda^\theta(M)) = \calV(N_{\widetilde{M}})
\]
where $\calV(I)$ is closed subvariety associated to ideal $I$.
Note that $\Char(M) \subseteq \pi_\theta^{-1}(\{y=0\}) = 
 \U_0 \cup \dots \cup \U_{l-1}$
for $M \in \catO_\lambda$. 

Let $\Min \widetilde{M}$ be a set of prime ideals minimal over
$N_{\widetilde{M}}$. If $P \in \Min \widetilde{M}$, let 
$n_{\widetilde{M}, P}$ be the length of the $S_{P}$-module 
$(\gr \widetilde{M})_{P}$. Define the characteristic cycle 
$\Ch(M)$ of $M$ to be 
\[
 \Ch(M) = \sum_{P \in \Min \widetilde{M}} n_{\widetilde{M}, P} 
 \calV(P).
\]
Let $\Min' \widetilde{M}$ be a subset of $\Min \widetilde{M}$
consisting of those prime ideals $P$ for which $\dim \calV(P)
= \dim \calV(N_{\widetilde{M}})$. Define the restricted 
characteristic cycle $\rCh(M)$ to be
\[
 \rCh(M) = \sum_{P \in \Min' \widetilde{M}}
 n_{\widetilde{M}, P} \calV(P).
\]

The following lemmas are proved in \cite{GS2}. 

\begin{lemma}
The characteristic varieties, the characteristic cycles and
the restricted characteristic cycles are independent of
the choice of good filtrations. 
\end{lemma}

\begin{lemma}
 \label{lemma:18}
 Let $0 \rightarrow A \rightarrow B \rightarrow C \rightarrow 0$
be a short exact sequence of finitely generated 
$\A_\lambda$-modules. Then one of the following cases occurs:
\begin{enumerate}
 \item $\dim \Char(A) < \dim \Char(B) = \dim \Char(C)$ and
       $\rCh(B) = \rCh(C)$.
 \item $\dim \Char(A) = \dim \Char(B) > \dim \Char(C)$ and
       $\rCh(B) = \rCh(A)$.
 \item $\dim \Char(A) = \dim \Char(B) = \dim \Char(C)$ and
       $\rCh(B) = \rCh(A) + \rCh(C)$.
\end{enumerate}
\end{lemma}

By \refthm{thm:2}, we have the isomorphism
\[
 \Phi_\lambda^\theta(e_0 \T_\lambda) \simeq \widetilde{\calP}_\theta.
\]
By the isomorphism of $\A_\lambda$-modules
\[
 e_0 \T_\lambda / e_0 \T_\lambda A^* \simeq
 \bigoplus_{i=0}^{l-1} e_0 \Delta_\lambda(i),
\]
we have
\begin{equation}
 \label{eq:9}
  \widetilde{\calP}_\theta / \widetilde{\calP}_\theta \bar{A}^* \simeq
  \Phi_\lambda^\theta( e_0 \T_\lambda / e_o \T_\lambda A^*)
  \simeq \bigoplus_{i=0}^{l-1} \Phi_\lambda^\theta(
  e_0 \Delta_\lambda(i))
\end{equation}
as coherent sheaves on $\M_\theta(\delta)$. Thus we can
determine the characteristic cycles from \refeq{eq:9}. 

\begin{proposition}
 \label{prop:13}
 For $\lambda \in \widetilde{\R}^l_{reg}$ and $\theta \in \Z^l_\lambda$,
 we have an isomorphism of coherent sheaves
\[
 \Phi_\lambda^\theta(e_0 \Delta_\lambda(i))
 \simeq (\widetilde{\calP}_\theta / \widetilde{\calP}_\theta \bar{A}^*) e_i.
\]
\end{proposition}

\begin{proposition}
 \label{prop:14}
 For $\lambda \in \widetilde{\R}^l_{reg}$, $\theta \in \Z^l_\lambda$
 and $i=0$, $\dots$, $l-1$,
 we have
 \[
 \Ch_\theta(e_0 \Delta_\lambda(i))  
 \simeq [\:\U_{i}\: ] + \sum_{\theta_{i} + \dots + \theta_{j-1} < 0} 
 [\: \U_{j} \:] = \sum_{i \unrhd_{\theta} j} [\: \U_j \:]
 = \sum_{i \geogeq{\theta} j} [\: \U_j \:].
 \]
\begin{proof}
 Since 
 the sheaf $\widetilde{\calP}_\theta$ is the tautological
 bundle on $\M_\theta(\delta)$ defined by \refeq{eq:47}, 
 The action of $\bar{A}^*$ 
on a fiber 
 $\widetilde{\calP}_\theta([a_k, b_k]_{k=0, \dots, l-1}) \simeq 
\bigoplus_{i=0}^{l-1} \C E_{0i}$ is given by
 the right-multiplication of a matrix 
 $\sum_{k=0}^{l-1} b_k E_{k-1, k}$. By the definition of 
 $\U_j$, we have $b_i = 0$ if and only if 
$\theta_{i} + \dots + \theta_{j-1} < 0$
 at $[a_k, b_k]_{k=0, \dots, l-1} \in \U_j$. Therefore
 \[
  (\widetilde{\calP}_\theta / \widetilde{\calP}_\theta \bar{A}^*)([a_k, b_k]_{k=0, \dots, l-1})
 = \bigoplus_{\theta_{i} + \dots + \theta_{j-1} < 0} \C E_{0,i}
 \]
 at any point $[a_k, b_k]_{k=0, \dots, l-1} \in \U_j$. 

 By \refprop{prop:13}, we have
\[
 \Phi_\lambda^\theta(e_0 \Delta(\lambda)(i)) \simeq
  (\widetilde{\calP}_\theta / \widetilde{\calP}_\theta \bar{A}^*) e_{i}.
\]
 Therefore we have
 \[
  \Ch_\theta(e_0 \Delta(\lambda)(i)) \simeq
 [\: \U_i \:] + \sum_{\theta_{i} + \dots + \theta_{j-1} < 0} [\: \U_j \:].
 \]
 By the definition of $\unrhd_\theta$ and the definition of 
 $\geogeq{\theta}$, we have
\[
 i \rhd_\theta j \Leftrightarrow i \geogo{\theta} j \Leftrightarrow
 \theta_{i} + \dots + \theta_{j-1} < 0.
\]
Thus we have the statement of the proposition.
\end{proof}
\end{proposition}

Note that \refprop{prop:14} is an affirmative answer to
Question~10.2 of \cite{Go} in the case of $G = \Z / l \Z$.

\begin{corollary}
 Taking the reduced characteristic cycle $\rCh_\theta(M)$ of
 a module $M \in \catO_\lambda$ induces an isomorphism of
 vector spaces
 \[
  \rCh_\theta : K(\catO_\lambda) \otimes_{\Z} \C 
 \longrightarrow H^0(\pi_\theta^{-1}(\{y = 0\}), \C).
 \]
\end{corollary}

\begin{corollary}
 If $\ii_i \repgo{\lambda} \ii_j$ and there is no $i > k > j$ such that
 $\ii_i \repgo{\lambda} \ii_k$, then
 \[
  \rCh_\theta(e_0 L_{\lambda}(\ii_i))
 = [\: \U_{\ii_i} \:] + [\: \U_{\ii_{i-1}} \:] + \dots + [\: \U_{\ii_{j-1}} \:].
 \]
\begin{proof}
 By the assumption, we have an exact sequence 
\[
 0 \rightarrow e_0 \Delta_\lambda(\ii_{j}) \rightarrow
 e_0 \Delta_\lambda(\ii_i) \rightarrow
 e_0 L_\lambda(\ii_i) \rightarrow 0
\]
 of $\A_\lambda$-modules. By \reflemma{lemma:18}, we have
\[
 \rCh_\theta(e_0 \Delta_\lambda(\ii_i))
 = \rCh_\theta(e_0 \Delta_\lambda(\ii_j))
 + \rCh_\theta(e_0 L_\lambda(\ii_i)).
\]
Therefore we have the statement of the corollary.
\end{proof}
\end{corollary}

\end{document}